  \documentclass[journal]{IEEEtran}
\usepackage{cite}
\usepackage{color,yfonts}
\usepackage[pdftex]{graphicx}
\graphicspath{{fig/}{jpeg/}}
\usepackage[cmex10]{amsmath}
\usepackage{amssymb}
\usepackage{epstopdf}
\usepackage{multirow}
\usepackage{algorithm}
\usepackage{tikz}
\usetikzlibrary{shapes,arrows,matrix}
\usetikzlibrary{decorations.pathmorphing} 
\usetikzlibrary{fit}                    
\usetikzlibrary{backgrounds}    
\usepackage{verbatim}
\usepackage{algpseudocode}
\usepackage{amsmath,amssymb,lipsum}
\usepackage{hyperref}
\usepackage{comment}
\usepackage{graphicx}
\usepackage[font=small]{caption}
\usepackage{subcaption}
\usepackage{mathtools,enumerate}
\usepackage{multicol,lipsum}

\providecommand{\Ex}[1]{\mathbb{E}\left[#1\right]}
\usepackage{stix}

\input{mysymbol.sty}
\renewcommand{\theenumi}{\roman{enumi}}%

\newcommand{\ayche}{\mathscr{h}}

\newcommand{\overbar}[2][3]{{}\mkern#1mu\overline{\mkern-#1mu#2}}
\newcommand{\INDSTATE}[1][1]{\State\hspace{3mm}}
\newcommand{\INDSTATED}[1][1]{\State\hspace{6mm}}

\usepackage{theorem}
\newtheorem{theorem}{Theorem}

\newtheorem{lemma}{Lemma}
\newtheorem{proposition}{Proposition}

\newtheorem{example}{\bf Example}
\theoremstyle{definition}
\newtheorem{assumption}{Assumption}

\def \EE {{\mathbb{E}}}

\def \x {{\mathbf{x}}}

\def \y {{\mathbf{y}}}
\def \h {{\mathbf{h}}}

\newcommand{\qed}{\hfill\blacksquare}
 \providecommand{\Ex}[1]{\mathbb{E}\left[#1\right]}
 
 \providecommand{\norm}[1]{\left\|#1\right\|}
 \providecommand{\ip}[1]{\boldsymbol{\langle}#1\boldsymbol{\rangle}}

    \newcommand{\Et}[1]{\ensuremath{\EE[#1~|~\mathcal{F}_t]}}  
  \tikzstyle{agent}=[circle,
  thick,
  minimum size=.5cm,
  draw=blue!80,
  fill=blue!20]

  \tikzstyle{neighbor}=[circle,
  thick,
  minimum size=.5cm,
  draw=cyan!85!black,
  fill=cyan!40,
  decorate,
  decoration={random steps,
      segment length=2pt,
      amplitude=2pt}]
  
  \tikzstyle{local_nat}=[rectangle,
  thick,
  minimum size=.5cm,
  draw=red!80,
  fill=red!20]
  
  \tikzstyle{glob_nat}=[rectangle,
  thick,
  minimum size=.5cm,
  draw=red!100,
  fill=red!40]      
  
  \tikzstyle{background}=[rectangle,
  fill=gray!10,
  inner sep=0.2cm,
  rounded corners=5mm]
  
  \tikzstyle{background2}=[rectangle,
  fill=green!20,
  inner sep=0.2cm,
  rounded corners=5mm]
{\tiny }
\title{\vspace{-0cm}{Nonparametric Compositional Stochastic Optimization for Risk-Sensitive Kernel Learning}}

\author{Amrit~Singh~Bedi,
    Alec~Koppel,
    Ketan~Rajawat, and Panchajanya Sanyal
    \thanks{
        A.S. Bedi and A. Koppel contributed equally to this work. They both are with the U.S. Army Research Laboratory, Adelphi, MD, USA. (e-mail: amrit0714@gmail.com, akoppel@seas.upenn.edu.). K. Rajawat is with the Department of Electrical Engineering,
        Indian Institute of Technology Kanpur, Kanpur 208016, India (e-mail:
        ketan@iitk.ac.in). A part of this work was presented in American Control Conference (ACC), Philadelphia, USA, 2019 \cite{BEdi_ACC} and was spotlighted in \cite{9084339}.}
    \vspace{-0cm}    }
\begin{document}
\maketitle

\begin{abstract}%
In this work, we address optimization problems where the objective function is a nonlinear function of an expected value, i.e., compositional stochastic {programs}. We consider the case where the decision variable is not vector-valued but instead belongs to a Reproducing Kernel Hilbert Space (RKHS), motivated by risk-aware formulations of supervised learning.
  We develop the first memory-efficient stochastic algorithm for this setting, which we call Compositional Online Learning with Kernels (COLK). COLK, at its core a two time-scale stochastic approximation method, addresses the facts that (i) compositions of expected value problems cannot be addressed by stochastic gradient method due to the presence of an inner expectation; and (ii) the RKHS-induced parameterization has complexity which is proportional to the iteration index which is mitigated through greedily constructed subspace projections. \blue{We provide, for the first time, a non-asymptotic tradeoff between the complexity of a function parameterization and its required convergence accuracy for both strongly convex and non-convex objectives under constant step-sizes.} Experiments with risk-sensitive supervised learning demonstrate that COLK consistently converges and performs reliably even when data is full of outliers, and thus marks a step towards overfitting. \blue{Specifically, we observe a favorable tradeoff between model complexity, consistent convergence, and statistical accuracy for data associated with heavy-tailed distributions.}

\end{abstract}

\section{Introduction}\label{sec:intro}
In this work, we focus on compositional stochastic programming, a setting where the objective function is an expectation over a set of random convex functions, each of which depends on the expected value of a different random convex function. This problem setting has received recent attention in operations research \cite{wang2017stochastic,lian2016finite} and machine learning \cite{dai2017learning} due to its ability to gracefully address technicalities that arise
%
%
in the theory of Markov Decision Problems (MDPs) \cite{sutton2009convergent} and bias-variance issues in supervised learning \cite{dentcheva2017statistical}. Our goal is to solve this class of problems when the decision variable is not vector-valued, as in \cite{wang2017stochastic}, but is instead itself a function. This setting arises intrinsically in MDPs defined over the continuous state and action spaces \cite{Tolstaya_ACC2018} or when accounting for risk \cite{ahmed2006convexity} in supervised learning with nonlinear interpolators \cite{Kivinen2004}. 

The theory of optimization in function space began with variational calculus \cite{gelfand2000calculus} and Hamilton's Principle \cite{bailey1982hamilton}. However, in modern applications, we require solutions to such problems in situations where classical methods no longer apply. Two different issues arise: (1) how to evaluate the expectation (integral) and (2) how to parameterize the function so that tractable updates may be obtained. Setting aside (1) for now, to address (2), i.e., to handle the intractability of general functional optimization, one must restrict the function we seek to not only yield a computationally tractable formulation, but also one be rich enough to address common experimental settings. In learning theory, for instance, we typically restrict the function to be a neural network \cite{haykin1994neural} or a nonparametric basis expansion in terms of data \cite{shawe2004kernel}, whereas in control systems, polynomial interpolation \cite{jarvis2005control} and kriging \cite{rasmussen2004gaussian} are popular. In this work, we address the case where the function class is a Reproducing Kernel Hilbert Space (RKHS), motivated by a recently developed memory-efficient parameterization of a function that is infinite dimensional \cite{koppel2016parsimonious}. This approach subsumes polynomial interpolation \cite{jarvis2005control}, avoids the memory explosion associated with large sample-size kriging \cite{rasmussen2004gaussian}, and preserves convexity, thus avoiding convergence to poor stationary points rampant in neural network training \cite{safran2017spurious}. 

With the function class specified, we turn to discuss how to solve the associated functional stochastic program: doing so requires iterative stochastic methods \cite{Robbins1951,shapiro2014lectures}, since deterministic approaches \cite{Boyd2004} require computing gradients that depend on infinitely many realizations of a random variable, thus exhibiting prohibitive complexity. Unfortunately, standard stochastic gradient descent (SGD) is inapplicable to the compositional setting, because, for a single stochastic descent direction, one requires the evaluation of an additional inner expectation, an observation that was popularized in reinforcement learning as ``the double sampling problem" \cite{sutton2009convergent}. To ameliorate this issue, we develop a functional nonparametric extension of stochastic quasi-gradient (SQG) method \cite{ermoliev1983stochastic,wang2017stochastic}, which uses two time-scale stochastic approximation: one uses a quasi-stationary estimate of the inner expectation, whereas the other executes stochastic descent \cite{konda2004convergence}.

Different from the estimator for the inner expectation in \cite{wang2017stochastic}, in this work we consider a momentum (gradient tracking) scheme that uses a difference of instantaneous costs, which permits us to tighten existing rate analyses for compositional problems with minimal additional computational overhead. {Gradient tracking has been investigated extensively in the distributed optimization literature to have agents' track the gradient of the global objective -- see for instance, \cite{olfati2005consensus,yang2007distributed,zhu2010discrete} and more recently \cite{di2016next,qu2017harnessing,scutari2019distributed,pu2020distributed}. Here we adopt the spirit of using the previous iterate, but rather than track the gradient itself, we use history information to track the evaluation of the inner objective that appears in compositional problems, which to the best of our knowledge, is the first time it has been employed in this context, and permits us to refine existing rates for compositional stochastic programming. This contrasts with its use for improving the vanilla gradient estimation error that arises in non-convex optimization as studied by \cite{cutkosky2019momentum}. }

However, with the choice of $\ccalH$ as an RKHS, and sequential application of the Representer Theorem \cite{scholkopfgeneralized}, means that the function parameterization grows with the iteration index \cite{Kivinen2004}, and thus becomes untenable for expected value problems. While memory-reduction methods for RKHS optimization exist using, for instance, random feature approximations \cite{rahimi2008random}, forgetting factors \cite{Kivinen2004}, random dropping \cite{conf/icml/HoiWZJW12}, or projections onto subspaces of fixed size \cite{wang2012breaking}, \blue{none of the aforementioned approaches can ensure that the stochastic descent properties of the algorithm are nearly preserved with probability $1$}, and hence fail to (approximately) preserve almost sure convergence of the optimization method to which they are applied. Thus, it is an open challenge how to generalize the results of \cite{wang2017stochastic} to optimization over RKHS. 
%
%
Thus, our main contributions are to:
\begin{itemize}
	\item extend SQG to RKHS, which we compress with matching pursuit \cite{Pati1993} (Sec. \ref{sec:algorithm}). We tailor the compression to ensure  valid descent \cite{koppel2016parsimonious,pkgtd}. We call this method Compositional Online Learning with Kernels (COLK). 
	\item  \blue{characterize the number of iterations and the memory size requirement mathematically to achieve $\delta$-suboptimal solutions for strongly convex  (Theorem \ref{theorem:constant_stepsize_convergence0}) and non-convex (Theorem \ref{nonconvex_theorem}) compositional objectives in terms of mean-square error to the optimizer and the expected gradient norm, respectively. }
	%
	\item \blue{establish that the worst-case complexity of the function sequence is finite (Theorem \ref{model_order}), with an explicit dependence on the parameter dimension, algorithm step-size, and compression budget. This result establishes a trades-off between convergence accuracy and model complexity.}
	\item experimentally (Sec. \ref{sec:experiments}) validate this method on a problem instantiation defined by robust supervised learning. Doing so yields nonlinear statistical models whose bias \emph{and} variance are small, first on a synthetic data \texttt{regression outliers} which has a heavier tailed distribution, i.e., more outliers are present, and then on benchmark data: \texttt{lidar} \cite{ruppert2009semiparametric}. We observe that COLK yields \emph{consistently} accurate performance across training realizations, meaning that it does not overfit, in contrast to other methods that cannot minimize risk functionals \cite{dentcheva2017statistical}.
\end{itemize} 
\blue{Overall, then, this work differs from \cite{wang2017stochastic} in the following ways: (i) \cite{wang2017stochastic}  considers vector-valued decision variables, whereas in this case, we focus on the RKHS setting. This more general setting implies that the function class to which statistical models belong is universal, but this universality comes at the cost of complexity. In the vector setting, stochastic compositional gradients may be used for updates without any memory bottleneck, whereas for RKHS one must use approximate gradients to ensure tractable implementation, but this approximation incurs errors in the gradients, which we analyze. (ii) We consider a modification of stochastic quasi-gradient methods that employ a momentum scheme in the faster time-scale in order to improve the overall convergence rates. We further note that all analysis in \cite{wang2017stochastic} is for attenuating step-sizes, whereas our focus is on constant step-sizes. (iii) We further characterize the dependency of the algorithm's non-asymptotic parameterization complexity on the underlying dimension of the feature space under constant step-sizes (Theorem \ref{model_order}), an attribute that is unique to nonparametric statistics which does not appear in \cite{wang2017stochastic}. Moreover, we illuminate the differences between these works experimentally in Section \ref{sec:experiments}.}

\textbf{Notation:} All the scalars are denoted by letters in regular font $s$, vectors are bold $\bbs$, and matrices are capitalized $\bbS$. The notation $\norm{\cdot}_{\mathcal{H}}$ represents the norm in the RKHS symbolized by $\mathcal{H}$. Subspaces of the RKHS are denoted by $\ccalH_{\bbS}$. The inner product operator between two functions $f\in\mathcal{H}
$ and $f'\in\mathcal{H}$ is denoted by $\ip{f,f'}_{\mathcal{H}}$.  The expectation operator is symbolized by $\mathbb{E}$. The notation $\circ$ describes the compositional operator, for instance, $f\circ g(x)$ denotes the composition as $f(g(x))$. A kernel function is denote by $\kappa(\cdot,\cdot)$ which takes two vectors as arguments. Stacking of kernel evaluations is denoted as the empirical kernel map $\bbkappa_{\bbU}(\cdot)$ and the kernel matrix $\bbK_{\bbU, \bbV}$ is the matrix whose $(i,j)$ entry is given by $\kappa(\bbu_i, \bbv_j)$. $\ccalP_\ccalC$ denotes the orthogonal projection operator onto a convex set $\ccalC$ .  Gradients are denoted in the usual way, and partial derivative, for instance, of scalar valued function $\ell(\bbu)$ with respect to $u_i$ at $\bbv$ , is denoted by, $\frac{\partial \ell(\bbu)}{\partial u_i}|_{\bbu=\bbv}$.

\section{Compositional Stochastic Programming in RKHS}\label{sec:problem}
In this work, we focus on solving \emph{functional} optimization problems whose objective is a nonlinear function of an expected value. That is, the objective function is a composition of two functions, each of which is an expected value over a set of functions parameterized by a pair of random variables. More specifically, there are two sets of {random variables  $\{\bbxi_t\}\subset \reals^p$ and $\{\bbtheta_t\}\subset\mathbb{R}^p$. In general both the random variables are allowed to be dependent, but for the ease of analysis and understanding, we assume that $\bbxi \in \bbXi \subset \reals^p$,  $\bbtheta \in \bbTheta \subset \reals^p$ and $\bbxi$, $\bbtheta$ are independent of each other. Considering these random pairs, the cost takes the form $J(f):= (L\circ \bbH)(f)$, where  $\bbH(f)=\mathbb{E}_{\bbxi}\left[ \boldsymbol{\ayche}_{\bbxi}( f(\bbxi)) \right]$ is a map $\bbH: \ccalH \rightarrow \reals^m$ that is an expectation over a set of random functions $\boldsymbol{\ayche}_{\bbxi}( f(\bbxi))$. {Specific instances of $\ayche$ will be discussed in Examples \ref{eg_risk}-\ref{eg_mdps} shortly. Often in applications, $\ayche$ quantifies statistics of some loss function, possibly $\ell$.} Similarly,  $L(\bbu)=\mathbb{E}_{\bbtheta}\left[ \ell_{\bbtheta}(\bbu) \right] $ is a map {$L: \mathbb{R}^m \rightarrow \reals$} that is an expected value over a random collection variable. Further, $\ccalH$ is a function space to be subsequently specified. In this work, we focus on the functional compositional stochastic program:} 
\begin{align}\label{eq:main_problem}
	\min_{f\in\ccalH} \mathbb{E}_{\bbtheta}\left[ \ell_{\bbtheta}\left( \mathbb{E}_{\bbxi}\left[ \boldsymbol{\ayche}_{\bbxi}( f(\bbxi)) \right]\right) \right] + \frac{\lambda}{2}\|f \|^2_{\ccalH} \; .
\end{align}
\blue{We consider both the cases that $J(f)$ is convex with respect to function $f$} with Tikhonov regularizer $\frac{\lambda}{2}\|f \|_{\ccalH}$ to ensure strong convexity, defining the regularized loss \cite{friedman2001elements}
\begin{align}\label{regularized_loss}
	f^\star=\argmin_{f\in\ccalH}\ R(f) := J(f) + \frac{\lambda}{2}\|f \|^2_{\ccalH}\; ,
\end{align}
\blue{and the case where $R(f)$ is non-convex\footnote{In which case the regularizer is inconsequential, and may be set as $\lambda=0$.}}.
The feasible set  $\ccalH$ of \eqref{eq:main_problem}, the domain of $H$, and hence $J$, is not Euclidean space $\reals^p$, as in \cite{wang2017stochastic}, but instead is a Hilbert space equipped with a unique \emph{kernel function}, $\kappa: \ccalU \times \ccalU \rightarrow \reals$, such that:
\begin{align} \label{eq:rkhs_properties}
	& (i) \  \langle \!f, \kappa(\bbu,\!\cdot) \rangle _{\ccalH} \!=\! f(\bbu) \ \  \text{for all } \bbu \in \ccalU,\nonumber
	\\
	&(ii) \ \ccalH = \overbar{\text{span}\{ \kappa(\bbu , \cdot) \}} \quad\text{for all } \bbu \in \ccalU\; .
\end{align}
where $\langle \cdot , \cdot \rangle_{\ccalH}$ denotes the Hilbert inner product for $\ccalH$ and $\ccalU := \bbXi \cup \bbTheta$ denotes the union of data domains $\bbXi$ and $\bbTheta$, whose elements $\bbu$ are random variables $\bbxi$ or $\bbtheta$. We assume that the kernel is nonnegative, i.e. $\kappa(\bbu, \bbu') \geq 0$ for all $\bbu, \bbu' \in \ccalU$ so that it is a Mercer kernel. Function spaces of this type are called reproducing kernel Hilbert spaces (RKHS) \cite{kimeldorf1971some}. 

In \eqref{eq:rkhs_properties}, property (i) is called the reproducing property of the kernel and comes from the Riesz Representation Theorem \cite{wheeden1977measure}. Replacing $f$ by $\kappa(\bbu' , \cdot) $  in \eqref{eq:rkhs_properties} (i) yields the expression $ \langle \kappa(\bbu', \cdot) , \kappa(\bbu, \cdot) \rangle_{\ccalH} = \kappa(\bbu, \bbu')$, which is why $\kappa$ is called ``reproducing."  This property provides a practical means by which to access a nonlinear transformation of the input space $\ccalU$.  Specifically, denote by $\phi(\cdot)$ a nonlinear map of the feature space that assigns to each $\bbu$ the kernel function $\kappa(\cdot, \bbu)$. Then the reproducing property of the kernel allows us to write the inner product of the image of distinct feature vectors $\bbu$ and $\bbu'$ under the map $\phi$ in terms of kernel evaluations only: $\langle \phi(\bbu), \phi(\bbu') \rangle_{\ccalH} =\kappa(\bbu, \bbu')$. This is the \emph{kernel trick}, and it provides a principled method for function estimation.

Moreover, property \eqref{eq:rkhs_properties} (ii) states that any function $f\in \ccalH$ may be written as a linear combination of kernel evaluations. For kernelized and regularized empirical risk minimization (i.e., the sample average approximation of \eqref{eq:main_problem} for some fixed $N$ realizations of $\bbxi$ and $\bbtheta$), the Representer Theorem \cite{kimeldorf1971some,scholkopfgeneralized} establishes that the optimal $f$ in function class $\ccalH$ may be written as an expansion of kernel evaluations \emph{only} at elements of the training set as\footnote{\blue{In the non-convex setting, this holds when $R$ is either quasi-convex \cite{boyer2019representer} or strongly convex in a neighborhood of a stationary point.}}
\begin{equation}\label{eq:kernel_expansion}
	f(\bbu) = \sum_{n=1}^{N} w_n \kappa(\bbxi_n, \bbu) \; .
\end{equation}
where $\bbw = [w_1, \cdots, w_N]^T \in \reals^N$ denotes the weight vector. The upper summand index $N$ in \eqref{eq:kernel_expansion} is henceforth referred to as the model order. Common choices $\kappa$ include the polynomial kernel and the radial basis kernel, i.e., $\kappa(\bbu,\bbu') = \left(\bbu^T\bbu'+b\right)^c $ and $\kappa(\bbu,\bbu') = \exp\left\{ -\frac{\lVert \bbu - \bbu' \rVert_2^2}{2c^2} \right\}$, respectively, with $\bbu, \bbu' \in \ccalU$. 

Then, one may use the Representer Theorem \eqref{eq:kernel_expansion} to transform the sample average approximation of \eqref{eq:main_problem} over all of $\ccalH$ into the parametric problem of two  $N$-dimensional weight vector $\bbw$. However, as $N\rightarrow \infty$, the function representation becomes infinite as well. In this work, we seek to find functions that are close-to-optimal solutions to \eqref{eq:main_problem} but also admit a finite-memory representation.
Before turning to develop an algorithmic tool which does so, we note that the problem setting \eqref{eq:main_problem} arises in diverse applications. Here we mention two, the first of which is the focus of this work. {In both instances, we focus on the case where samples from the data distribution are revealed incrementally, the input-ouput relationship between features and targets is not necessarily linearly separable, and hence feature selection must done on the fly during training. This setting makes RKHS parameterizations advantageous in terms of balancing descriptive richness with preserving convexity, and hence defining training schemes that may find the optimal representation within the defined function class.}
\vspace{-0mm}
\begin{example}[\bf Robust Supervised Learning]\label{eg_risk}\normalfont
	Consider a random pair $(\bbx, \bby)\in\ccalX\times\ccalY$, realizations of which are training examples $(\bbx_n, \bby_n)$, and $\ccalX \subset \reals^p$, the $p$-dimensional Euclidean space. In comparison to formulation in \eqref{eq:main_problem}, we have $\bbtheta=\bbxi=\bbx$ and $\bby^{\bbtheta}=\bby^{\bbxi}=\bby$ for this example {which represents the corresponding target values.} In the case of classification with $C$ classes, $\ccalY=\{1,\dots,C\}$, whereas in the case of regression $\ccalY\subset\reals^q$. In supervised learning, we learn an estimator $f(\bbx)$ according to its ability to minimize a loss function {$l : \ccalH \times \ccalX \times \ccalY \rightarrow \reals$} averaged over data:
	{\begin{equation}\label{eq:supervised_loss}
			f^\star=\argmin_{f\in\ccalH} \mathbb{E}[l(f(\bbx_{\bbtheta}),\bby_{\bbtheta})] \; ,
	\end{equation}}
	where we define { $L(f) =\mathbb{E}[l(f(\bbx_{\bbtheta}),\bby_{\bbtheta})] $,} and ignore the regularizer for the moment. The loss {$l$} quantifies the merit of the estimator $f(\bbx)$ with respect to its target $\bby$. However, as it is well known in statistics \cite{friedman2001elements}, solving \eqref{eq:supervised_loss} is really only an approximation of the Bayes optimal estimator
	%
	\begin{equation}\label{eq:bayes_risk}
		\hat{\bby}^\star=\argmin_{\hat{\bby}\in\ccalY^{\ccalX}} \mathbb{E}[l(\hat{\bby}(\bbx_{\bbtheta}),\bby_{\bbtheta})] \; ,
	\end{equation}
	where $\ccalY^{\ccalX}$ denotes the space of all functions  $\hat{\bby}:\ccalX\rightarrow \ccalY$ that map data $\bbx$ to target variables $\bby$. Suppose we obtain some estimate $\hat{f}$ by approximately minimizing \eqref{eq:supervised_loss}. Then the performance difference associated with $\hat{f}$ and the Bayes optimal $\hat{\bby}^\star$  \eqref{eq:bayes_risk}  is given as
	\begin{align}\label{eq:bias_variance}
		&\underbrace{\mathbb{E}[l(\hat{f}(\bbx_{\bbtheta}),\bby_{\bbtheta})]  - \min_{f\in\ccalH} \mathbb{E}[l(f(\bbx_{\bbtheta}),\bby_{\bbtheta})]}_{\text{\normalfont Estimation Error}}\nonumber\\
		&\quad\quad + 
		\underbrace{ \min_{f\in\ccalH} \mathbb{E}[l(f(\bbx_{\bbtheta}),\bby_{\bbtheta})] - \min_{\hat{\bby}\in\ccalY^{\ccalX}} \mathbb{E}[l(\hat{\bby}(\bbx_{\bbtheta}),\bby_{\bbtheta})]}_{\text{\normalfont Approximation Error}}\; ,
	\end{align}
	where we add and subtract the optimal supervised cost $\min_{f\in\ccalH} \mathbb{E}[l(f(\bbx_{\bbtheta}),\bby_{\bbtheta})]$ to obtain that this discrepancy decomposes into two terms: the estimation error, and approximation error. In most cases, the estimation error may be identified with bias, and the approximation error may be identified with variance plus noise -- we subsequently identify the decomposition \eqref{eq:bias_variance} with bias and variance, as bias and estimation error are both defined as subsampling errors, whereas variance and approximation error are both defined by quality of model fit across  training runs. See \cite{poggio2003mathematics}[Sections 3.1 - 3.2].
	In supervised learning, typically we try to make the model bias small as the number of data points goes to infinity, resigning ourselves to the fact that the variance is an intrinsic penalty we suffer for selecting a particular modeling hypothesis: the function class to which $f$ belongs.
	However, authors in operations research \cite{ahmed2006convexity} and applied probability \cite{ruszczynski2006optimization} have proposed to optimize both the expected loss over all data plus a measure of the \emph{dispersion} of the estimate with respect to its target variable as a way of accounting for the unknown approximation error of the modeling hypothesis in \eqref{eq:bias_variance}. Many measures of dispersion are possible, but one which yields a convex formulation is the semivariance.
	\begin{equation}\label{eq:semivariance}
		\!\!\widetilde{\text{Var}}[l(f(\bbx_{\bbtheta}),\bby_{\bbtheta})]= \mathbb{E}\Big\{ \!\big(l(f(\bbx_{\bbtheta}),\bby_{\bbtheta})- {\mathbb{E}[l(f(\bbx_{\bbxi}),\bby_{{\bbxi}})]} \big)_+^2\Big\}\; ,
	\end{equation}
	where $a_+=\max(a,0)$ denotes the positive projection. Note that when we omit the positive projection, \eqref{eq:semivariance} reduces to the variance of the instantaneous loss $l(f(\bbx_{\bbtheta}),\bby_{\bbtheta})$. To see that variance is a composite function observe that when $\bbtheta$ and $\bbxi$ are independent, we have that
	\begin{align}
		\ell_{\bbtheta}(\bbu) &= (u^1-u^2)^2 & L(\bbu) &= \Ex{(u^1-u^2)^2}
		\\
		\!\!\!\![\boldsymbol{\ayche}_{\bbxi}(f({\bbtheta}))]_1 &= l(y_{\bbtheta} - f({\bbtheta})) & [\bbH(f)]_1 &=  l(y_{\bbtheta} - f({\bbtheta}))
		\\
		\!\!\!\![\boldsymbol{\ayche}_{\bbxi}(f({\bbxi}))]_2 &= l(y_{\bbxi}- f({\bbxi})) & [\bbH(f)]_2 & = \Ex{l(y_{\bbxi} - f({\bbxi}))}
	\end{align}
	where note that $u_1$ and $u_2$ are also random quantities so the first expectation cannot be dropped.

	However, the subtraction of the second moment of $l(f(\bbx_{\bbtheta}),\bby_{\bbtheta})$ without positive projection makes the problem non-convex. Thus, using positive projections yields semivariance rather than true variance \cite{ahmed2006convexity}.  With the measure of dispersion in \eqref{eq:semivariance}, one robust formulation of supervised learning over an RKHS $\ccalH$ which accounts for approximation error is
	\begin{equation}\label{eq:robust_supervised_loss}
		\!f^\star\!\!=\argmin_{f\in\ccalH} \mathbb{E}[l(f(\bbx_{\bbtheta}),\bby_{\bbtheta})]  + \eta  \widetilde{\text{Var}}[l(f(\bbx_{\bbtheta}),\bby_{\bbtheta})]\! + \frac{\lambda}{2}\|f \|_{\ccalH}^2 .
	\end{equation}
	where $\eta$ is a scaling parameter that tunes the importance of estimation or approximation error. Note that we have added the regularizer back into \eqref{eq:robust_supervised_loss}. Solutions of \eqref{eq:robust_supervised_loss}, as compared with \eqref{eq:supervised_loss}, are better attuned to data points associated with high variance objective evaluation, which practically may be interpreted as outliers or in the case of classification, situations where training examples possess characteristics corresponding to multiple classes. Due to the fact that our analysis requires Lipschitz gradients (Section \ref{sec:convergence}), numerically we approximate the positive projection in \eqref{eq:semivariance} with the softmax in Section \ref{sec:experiments} -- see \cite{Boyd2004}. 
	An alternative risk measure popular in finance is the conditional value-at-risk (CVaR) \cite{uryasev2000conditional}, which quantifies the loss function at different quantiles of its distribution. 
\end{example}

\begin{example}[\bf  Policy Evaluation in MDPs]\normalfont\label{eg_mdps}
	
	%
	Another instantiation of \eqref{eq:main_problem} is the task of policy evaluation in a continuous Markov Decision Problem (MDP) \cite{Bellman:1957}, a quintuple $(\ccalX, \ccalA, \mathbb{P}, r, \gamma)$, where $\mathbb{P}$ is the action-dependent transition probability of the process: when the agent starts in state $\bbx_t \in \ccalX \subset \reals^p$ at time $t$ and takes an action $\bba_t \in \ccalA$, a transition to next state $\bby_t\in \ccalX$ is distributed according to 
	%
	$\bby_t \sim \mathbb{P}(\cdot \given \bbx_t, \bba_t).$
	%
	After the agent transitions to a particular $\bby_t$, the MDP provides to it an instantaneous reward $r(\bbx_t, \bba_t, \bby_t)$, where the reward function is a map $r :\ccalX \times \ccalA \times \ccalX \rightarrow \reals$.
	
	In \emph{policy evaluation}, control decisions $\bba_t$ are chosen according to a stochastic stationary policy $\pi: \ccalX \rightarrow \rho(\ccalA)$, where $\rho(\ccalA)$ denotes the set of probability distributions over $\ccalA$, and one seeks to compute the \emph{value} of a policy when starting in state $\bbx$, quantified by the discounted expected sum of rewards, or value function $V^{\pi}(\bbx)$:
	\begin{equation}\label{eq:value_function}
		\!\!\!\!V^{\pi}(  \bbx ) = \mathbb{E}_{\bby}\Big[\sum_{t=0}^\infty \gamma^t r(\bbx_t, \bba_t, \bby_t ) \given \bbx_0 = \bbx,\{ \bba_t = \pi(\bbx_t)\}_{t=0}^\infty \Big]\; .
	\end{equation}
	For a single trajectory through the state space $\ccalX$, $\bby_t = \bbx_{t+1}$. The discount factor $\gamma \in (0,1)$ determines the agent's farsightedness. 
	%
	%
	From the definition of the value function in \eqref{eq:value_function}, one may derive the Bellman evaluation equation \cite{Bellman:1957}:
	\begin{equation}\label{eq:bellman_value}
		V^\pi( \bbx) = \int_{\ccalX} [ r(\bbx,\pi(\bbx),\bby) + \gamma V(\bby)] \mathbb{P}(d\bby \given \bbx, \pi(\bbx)), \;
	\end{equation}
	for all $\bbx \in \ccalX$. The functional fixed point problem \eqref{eq:bellman_value} defined by Bellman's equation may be reformulated as a nested stochastic program. To do so, rewrite the integral as an expectation, subtract the value function $V^{\pi}(\bbx)$ that satisfies the fixed point relation from both sides, and then pull it inside the expectation. Then, to solve \eqref{eq:bellman_value} in an {initialization-independent} manner, integrate out $\bbx$, the starting point of the trajectory defining the value function \eqref{eq:value_function}, as well as policy $\pi(\bbx)$. Doing so, \blue{i.e., considering the \emph{Galerkin relaxation} of \eqref{eq:bellman_value} followed by squaring the resulting expression,} yields the compositional objective
	\begin{align}\label{eq:main_prob_0}
		J\!(V)\!\!:=\! \mathbb{E}_{\bbx, \pi(\bbx)}\!\Big\{\!\frac{1}{2}\big( \mathbb{E}_{\bby}\big[ r(\bbx,\pi(\bbx),\bby) \!+\! \gamma V(\bby) \!-  \!V (\bbx)   \given \bbx, \pi(\bbx)\big]\big)^2\!\Big\}.
	\end{align}
	However, since it is intractable to optimize over all bounded functions $\ccalB(\ccalX)$, one may restrict the minimization of $J(V)$ to an RKHS $\ccalH$. This hypothesis, however, requires the introduction of regularization. Assuming that the Bellman fixed point $V^{\pi}$ is a continuous function, the RKHS approximation may be made close to the true $V^{\pi}$ when used with a universal kernel \cite{micchelli2006universal} -- see \cite{pkgtd}. We note that the problem formulation in  \eqref{eq:main_prob_0} is slight generalization of the problem formulation in \eqref{eq:robust_supervised_loss}. In particular, observe that the inner random variable $\y$ in \eqref{eq:main_prob_0} depends on the outer random variables $\x$ and $\pi(\x)$, necessitating the use of the conditional expectation operator. However, this generalization can be readily handled as detailed in \cite{wang2017stochastic} and the subsequent analysis considers \eqref{eq:robust_supervised_loss} for brevity. 
	
	With the problem setting clarified, we next shift focus to developing an iterative numerical method to solve \eqref{eq:main_problem}.

\end{example}

\section{Algorithm Development}
\label{sec:algorithm}
Now we turn to solving the stochastic compositional optimization \eqref{eq:main_problem} over the RKHS $\ccalH$. We focus on the development of stochastic approximation methods such that we may minimize $R$ over $\ccalH$ with only sequentially revealed independent and identically distributed realizations of $\bbtheta$ and $\bbxi$ without knowledge of their underlying  distributions. Related ideas are developed for the specialized objective of Example \ref{eg_mdps} in \cite{pkgtd} and in the vector-valued case in \cite{wang2017stochastic}. \blue{A key distinguishing feature of our method is the use of a gradient tracking scheme to reduce the estimation error of the inner expectation -- this is the first time momentum-style updates have been employed for compositional problems, and their use permits us to improve upon existing rate analyses.} The fundamental building blocks of our proposed algorithm are a functional generalization of the stochastic quasi-gradient method (Section \ref{subsec:sqd}) operating in tandem with low-dimensional subspace projections that are greedily constructed using matching pursuit (Section \ref{subsec:projection}), which we detail next.

\subsection{Functional Stochastic Quasi-gradient Descent}\label{subsec:sqd}

Note that the functional gradient of the objective function $J(f)$ is given by 
\begin{align}
	\langle\mathbb{E}\left[\nabla _f\boldsymbol{\ayche}_{\bbxi}(f(\bbxi))\right], \mathbb{E}\left[\nabla \ell_{\bbtheta}(\mathbb{E}[\boldsymbol{\ayche}_{\bbxi}(f(\bbxi))])\right]\rangle
\end{align} where $\nabla _f\boldsymbol{\ayche}_{\bbxi}(f(\bbxi)) \in \mathcal{H}\times \mathbb{R}^m$ and $\nabla \ell_{\bbtheta}\in \mathbb{R}^m$ which is defined as $\left[\nabla \ell_{\bbtheta}(\bbu)\right]_i=\frac{\partial \ell_{\bbtheta}(\bbu)}{\partial u_i}$.  Observe that the stochastic version of the functional gradient is obtained by dropping the outer expectation as follows
\begin{align}\label{eq:stochastic_gradient}
	\langle\nabla _f\boldsymbol{\ayche}_{\bbxi_t}(f(\bbxi_t)),\nabla \ell_{\bbtheta_t}(\mathbb{E}[\boldsymbol{\ayche}_{\bbxi}(f(\bbxi))])\rangle.
\end{align}
In vanilla stochastic gradient descent, one descends along the stochastic gradient \eqref{eq:stochastic_gradient}, i.e., performs the update given by 
\begin{align}\label{eq:vanilla_sqg_descent}
	\!\!\! f_{t+1}\!=\!(1\!-\!\lambda\alpha)f_t\!-\!\alpha\langle\nabla _f\boldsymbol{\ayche}_{\bbxi_t}({f_t(\bbxi_t)}),\!\nabla \ell_{\bbtheta_t\!}(\mathbb{E}[\boldsymbol{\ayche}_{\bbxi}(f_t(\bbxi_t))])\rangle.\!
\end{align}
However, in \eqref{eq:vanilla_sqg_descent}, the stochastic gradient at a specific random variable $\bbxi_t$, $\bbtheta_t$ is not available due to the expectation involved in the argument of $\nabla \ell_{\bbtheta_t}(\mathbb{E}[\boldsymbol{\ayche}_{\bbxi}(f(\bbxi))])$. This issue precludes use of vanilla stochastic gradient method for solving \eqref{eq:main_problem}.

Thus, we propose using a two time-scale stochastic approximation strategy called stochastic quasi-gradient method \cite{Korostelev,ermoliev1983stochastic}. To estimate the inner-expectation, one may recursively average an instantaneous approximation of the inner expectation $\boldsymbol{\ayche}_{\bbxi_t}(f(\bbxi_t))$ evaluated at $\bbxi_t$ as \cite{wang2017stochastic} \vspace{0mm}
\begin{align}\label{eq:auxiliary_update}
	\bbg_{t+1}=(1-\beta)\bbg_t+\beta\boldsymbol{\ayche}_{\bbxi_t}( f_t(\bbxi_t))
\end{align}
with the intent of estimating the expectation $\mathbb{E}_{\bbxi}\left[ \boldsymbol{\ayche}( f(\bbxi)) \right]$. \blue{By contrast, in this work, we consider a momentum analogue of this update instead \vspace{0mm}
	\begin{align}
		\bbg_{t+1} &= (1-\beta)(\bbg_t - \boldsymbol{\ayche}_{\bbxi_t}( f_{t-1}(\bbxi_t))) + \boldsymbol{\ayche}_{\bbxi_t}( f_t(\bbxi_t)). \label{track1}
	\end{align}
	Observe here that only a single observation $\xi_t$ is required per update, but two function evaluations of $\ayche$.} In \eqref{eq:auxiliary_update}, $\beta$ is a scalar learning rate chosen from the unit interval $(0,1)$. {This update may be interpreted as belonging to the family of gradient tracking schemes -- see \cite{scutari2019distributed}[p. 500-501] for illuminating discussion.} Then, we define a function sequence $f_t \in \ccalH$ initialized as null $f_0=0$ is sequentially updated using stochastic quasi-gradient descent:
\begin{align}\label{eq:sqg_descent}
	f_{t+1}=(1-\lambda\alpha)f_t-\alpha\langle\nabla _f\boldsymbol{\ayche}_{\bbxi_t}({f_t(\bbxi_t)}),\nabla \ell_{\bbtheta_t}(\bbg_{t+1})\rangle  \; ,
\end{align} 
where $\alpha$ is a step-size parameter \blue{we fix as constant throughout}. 
Further note that the term $\langle\nabla _f\boldsymbol{\ayche}_{\bbxi_t}({f_t(\bbxi_t)})$ is a function in $\ccalH$, and thus infinite dimensional. However, by applying the chain rule and the reproducing property of the kernel [cf. \eqref{eq:rkhs_properties}(i)], we obtain
\begin{align} 
	\langle\nabla _f&\boldsymbol{\ayche}_{\bbxi_t}(f(\bbxi_t)),\nabla \ell_{\bbtheta_t}(\bbg_{t+1})\rangle
	=\sum_{i=1}^{m}\nabla_f \boldsymbol{\ayche}^i_{\bbxi_t}(f(\bbxi_t))\frac{\partial \ell_{\bbtheta_t}(\bbu)}{\partial u_i}|_{\bbu=\bbg_{t+1}}\nonumber\\
	&\quad=\sum\limits_{i=1}^{m}\frac{\partial \boldsymbol{\ayche}^i_{\bbxi_t}(\omega)}{\partial \omega} |_{\omega=f(\bbxi_t)}\times{\frac{\partial \ell_{\bbtheta_t}(\bbu)}{\partial u_i}}|_{\bbu=\bbg_{t+1}}\kappa(\bbxi_t,\cdot)\label{first}\\
	&\quad=\langle\boldsymbol{\ayche}'_{\bbxi_t}( f_t(\bbxi_t)),\ell'_{\bbtheta_t}(\bbg_{t+1})\rangle \kappa(\bbxi_t,\cdot) \label{second}.
\end{align}
%
In \eqref{second}, we have used the vector inner product notation to denote the summation in \eqref{first}. Note that the kernel function $\kappa(\bbxi_t,\cdot)$ in \eqref{first} is common and therefore outside the inner product in \eqref{second}.
Utilizing this notation, the function update equation of \eqref{eq:sqg_descent} may be written as 
\begin{align}\label{eq:function_update_2}
	{f}_{t+1}=(1-\lambda\alpha)f_t-\alpha\langle\boldsymbol{\ayche}'_{\bbxi_t}( f_t(\bbxi_t)),\ell'_{\bbtheta_t}(\bbg_{t+1})\rangle \kappa(\bbxi_t,\cdot)\; . 
\end{align}
Observe that in \eqref{eq:function_update_2}, the vector $\boldsymbol{\ayche}'( {f_t(\bbxi_t)})$ associates with the sample point $\bbxi_t$ evaluated by the kernel, whereas $\ell'(\bbg_{t+1})$ associates with the {tracking parameter $\bbg_{t+1}$}. Moreover, the function sequence in \eqref{eq:function_update_2} belongs to a RKHS defined over \blue{the span of kernels $\{\kappa(\bbxi,\cdot)\}_{\bbxi\in\Xi}$}. Specifically, using the fact that $f_0=0\in\ccalH$, one may obtain through induction that the function $f_t$ at time $t$ admits an expansion via kernel evaluations of past data realizations $\bbxi_n$ and scalar weights $w_n$ for $n<t$:
\begin{align}\label{eq:kernel_expansion_t}
	f_t(\bbu)=\sum\limits_{n=1}^{t-1}w_n\kappa(\bbxi_n ,\bbu)=\bbw_t^T\mathbf{\kappa}_{\bbU_t}(\bbu)\; ,
\end{align}
where we define weight vector $\bbw_t:=[w_1,\cdots,w_{t-1}]^T$, kernel dictionary $\bbU_t = [\bbxi_1;\cdots ; \bbxi_{t-1}]$, and the empirical kernel map $\mathbf{\kappa}_{\bbU_t}(\bbu):=[\kappa(\bbxi_1,\bbu),\cdots, \kappa(\bbxi_{t-1},\bbu)]$. Thus, performing the stochastic quasi-gradient iteration in the RKHS amounts to the following parametric updates on the coefficient vector $\bbw$ and kernel dictionary $\bbU$, given by

\begin{align}\label{eq:parametric_updates}
	\bbU_{t+1}\!=&\!\left[\bbU_t, \bbxi_t \right],  \nonumber
	\\
	\bbw_{t+1}\!\!=&\!\left[(1-\alpha\lambda)\bbw_t, \!-\alpha\langle{\boldsymbol{\ayche}'_{\bbxi_t}\!( f_t(\bbxi_t))},\ell'_{\bbtheta_t}\!(\bbg_{t+1})\rangle\!\right]\!.  
\end{align}
In \eqref{eq:parametric_updates}, observe the \emph{kernel dictionary} parameterizing function $f_t$ is a matrix $\bbU_t\in \reals^{p\times (t-1)}$ which stacks past realizations of random variable $\bbxi$, and the coefficient vector $\bbw_t\in\reals^{t-1}$ as the associated scalars in the kernel expansion \eqref{eq:kernel_expansion_t} which are updated according to \eqref{eq:parametric_updates}.
Observe that the function update of \eqref{eq:function_update_2} implies that the complexity of $f_t$ is $\ccalO(t)$, due to the fact that the number of columns in $\bbU_t$, or \emph{model order} $M_t$, is $(t-1)$, and thus is unsuitable for settings where the total number of data samples is not finite, or are arriving sequentially and repeatedly. This is an inherent challenge of extending \cite{wang2017stochastic} to optimizing over nonlinear functions that belong to RKHS. To address this, we consider projections of \eqref{eq:function_update_2} onto low-dimensional subspaces, inspired by  \cite{koppel2016parsimonious}, which we detail in the following subsection.
%
%
\subsection{Subspace Projections for Complexity Control}\label{subsec:projection}
In this subsection, we turn to address the untenable growth of the function representational complexity discussed in the previous section, namely, that the model order is $M_t = (t-1)$, and grows without bound with the iteration index $t$. To do so, we adopt the idea of bias-inducing proximal projections onto low-dimensional subspaces developed in  \cite{koppel2016parsimonious}.

Specifically, we construct an approximate sequence of functions by orthogonally projecting functional stochastic gradient updates onto subspaces $\ccalH_\bbD \subseteq \ccalH$ that consist only of functions that can be represented using some dictionary $\bbD = [\bbd_1,\ \ldots,\ \bbd_M] \in {\reals^{ p \times M}}$, i.e., $\ccalH_\bbD = \{f\ :\ f(\cdot) = \sum_{n=1}^M w_n{\kappa}(\bbd_n,\cdot) = \bbw^T\boldsymbol{{\kappa}}_{\bbD}(\cdot) \}=\text{span}\{{\kappa}(\bbd_n, \cdot) \}_{n=1}^M$. Here we define ${\bbd_n\in \reals^{p}}$ as a model point which stacks exemplar realizations of $\bbxi$, i.e., $\bbd_n = \bbxi_n$. Further define $\boldsymbol{  {\kappa}}_{\bbD}(\cdot)=[{\kappa}(\bbd_1,\cdot) \ldots {\kappa}(\bbd_M,\cdot)]$, and $\bbK_{\bbD,\bbD}$ as the resulting kernel matrix from dictionary $\bbD_t$. The dictionary $\bbD_t$ is updated with the new sample $\bbxi_t$. Then to enforce projection, we replace the update in \eqref{eq:sqg_descent} with the following one which performs projection on to the subspace $\mathcal{H}_{\bbD_{t+1}}$:
%
%
%
%
%
\begin{align}\label{eq:projection_hat}
	\!\!\!\!{f}_{t+1} &\!\!= \!\argmin_{f \in \mathcal{H}_{\bbD_{t+1}}}\!
	\!\Big\lVert f \!\!-\!\left(\!{(1\!\!-\!\!\lambda\alpha)\!f_t\!-\!{\alpha\langle{\boldsymbol{\ayche}'_{\bbxi_t}\!( f_t(\bbxi_t))},\ell'_{\bbtheta_t}\!(\bbg_{t+1})\rangle \kappa(\bbxi_t,\cdot\!)} }\!\!\right)\!\!\Big\rVert_{\ccalH}^2 \nonumber \\
	&\!\!\!\!\!:=  \mathcal{P}_{\ccalH_{\bbD_{t+1}}}\!\!\left[{(1\!\!-\!\!\lambda\alpha)f_t\!-\!{\alpha\langle{\boldsymbol{\ayche}'_{\bbxi_t}( f_t(\bbxi_t))},\!\ell'_{\bbtheta_t}(\bbg_{t+1})\rangle \kappa(\bbxi_t,\cdot)}}\!\right]\!.
\end{align} 
where $\mathcal{P}_{\ccalH_{\bbD_{t+1}}}$ denotes the projection on to the subspace $\ccalH_{\bbD_{t+1}}$. 
\begin{algorithm}[h]
	\caption{Compositional Online Learning with Kernels (COLK)}
	\label{alg:colk}
	\begin{algorithmic}
		\Require $\{\bbtheta_t,\bbxi_t,\alpha,\beta,\epsilon \}_{t=0,1,2,...}$
		\State \textbf{initialize} ${f}_0(\cdot) = 0, \bbD_0 = [], \bbw_0 = []$, i.e. initial dictionary, coefficient vectors are empty
		\For{$t=0,1,2,\ldots$}
		\State \textbf{Update} auxiliary variable $\bbg_{t+1}$ according to \blue{\eqref{track1}}
		\blue{   \begin{align}\label{eq:g_update}
				\bbg_{t+1} &= (1-\beta)(\bbg_t - \boldsymbol{\ayche}_{\bbxi_t}( f_{t-1}(\bbxi_t))) + \boldsymbol{\ayche}_{\bbxi_t}( f_t(\bbxi_t))
		\end{align}}\vspace{0mm}
		\State \textbf{Update} function via a stochastic quasi-gradient step \eqref{eq:function_update_2}
		$$    \tilde{f}_{t+1}={(1-\lambda\alpha)f_t-{\alpha\langle{\boldsymbol{\ayche}'_{\bbxi_t}( f_t(\bbxi_t))},\ell'_{\bbtheta_t}(\bbg_{t+1})\rangle \kappa(\bbxi_t,\cdot)} }$$\vspace{0mm}
		\State {\textbf{Revise}  parametrization: dictionary and weights \eqref{eq:parametric_updates}
			\begin{align*}
				\tbD_{t+1} &= [\bbD_t,\;\;\bbxi_t] \\
				\tbw_{t+1} &= [(1-\alpha\lambda)\bbw_t,\;\; -{{\alpha\langle{\boldsymbol{\ayche}'_{\bbxi_t}( f_t(\bbxi_t))},\ell'_{\bbtheta_t}(\bbg_{t+1})\rangle}}]
			\end{align*}\vspace{0mm}
		}        
		\State \textbf{Compress} function representation using  KOMP as
		\begin{align}\label{eq:function_pdate}
			({f}_{t+1},\bbD_{t+1},\bbw_{t+1}) = \textbf{KOMP}(\tilde{f}_{t+1},\tbD_{t+1},\tbw_{t+1},\epsilon)
		\end{align}
		\EndFor
	\end{algorithmic}
\end{algorithm}
In order to perform the projection in \eqref{eq:projection_hat}, let us denote the unprojected version of the function as $\tilde{f}_{t+1}$, dictionary $\tilde{\bbD}_{t+1}$ and the weights $\tilde{\bbw}_{t+1}$ given by 
%
\begin{align}\label{eq:sgd_tilde}
	\tilde{f}_{t+1} =(1-\lambda\alpha)f_t-{\alpha\langle{\boldsymbol{\ayche}'_{\bbxi_t}( f_t(\bbxi_t))},\ell'_{\bbtheta_t}(\bbg_{t+1})\rangle \kappa(\bbxi_t,\cdot)}  \; .
\end{align}
This update may be represented parametrically as
\begin{align}\label{eq:param_tilde}
	\tbD_{t+1} \!=\! [\bbD_t,\;\bbxi_t], \;  \tbw_{t+1} \!=\! [(1\!-\!\lambda\alpha)\bbw_t\; ,\; -{\alpha\langle{\boldsymbol{\ayche}'_{\!\bbxi_t\!\!}( f_t(\bbxi_t))},\ell'_{\!\bbtheta_t\!\!}(\bbg_{t+1})\!\rangle}\! ] \; .
\end{align}
%
%
%
Given that we need to project $\tilde{f}_{t+1}$ onto the stochastic subspace $\ccalH_{\bbD_{t+1}}$, for a fixed dictionary $\bbD_{t+1}$, the stochastic projection in \eqref{eq:projection_hat} amounts to update the coefficient vector $\bbw_{t+1}$ as 
\begin{equation} \label{eq:hatparam_update}
	\bbw_{t+1}=  \bbK_{\bbD_{t+1} \bbD_{t+1}}^{-1} \bbK_{\bbD_{t+1} \tbD_{t+1}} \tbw_{t+1} \;,
\end{equation}
where $\bbK_{\bbD_{t+1},\bbD_{t+1}}$ and $\bbK_{\bbD_{t+1},\tbD_{t+1}}$ are the cross kernel matrices between $\{\bbD_{t+1},\bbD_{t+1}\}$ and $\{\bbD_{t+1},\tbD_{t+1}\}$, respectively. So the basic idea is to obtain a compressed version of $\tilde{f}_{t+1}$ and corresponding dictionary $\bbD_{t+1}$ and the coefficient vector $\bbw_{t+1}$. 
As previously noted, numerous approaches are possible for seeking this sparse representation. We make use of \emph{kernel orthogonal matching pursuit} (KOMP) \cite[Sec. 2.3]{Vincent2002} with allowed error tolerance $\epsilon$ to find a kernel dictionary matrix $\bbD_{t+1}$ based on the one which adds the latest sample point $\tbD_{t+1}$. This choice is due to the fact that we can tune its stopping criterion to guarantee a decrement property in expectation \cite{koppel2016parsimonious}, as well as ensure the function complexity remains finite -- see Sec. \ref{sec:convergence}. {
	We now describe the variant of KOMP we propose using, called Destructive KOMP with Pre-Fitting (see \cite{Vincent2002}, Section 2.3),  (Algorithm \ref{alg:komp}), which takes a candidate function $\tilde{f}$ of model order $\tilde{M}$ parameterized by its kernel {dictionary $\tbD\in\reals^{p\times\tilde{M}}$} and {coefficient vector $\tbw\in\reals^{\tilde{M}}$}, and approximates $\tilde{f}$ by a function $f\in \ccalH$ with a lower model order. Initially, this sparse approximation is the original function $f = \tilde{f} $ so that its dictionary is initialized with that of the original function $\bbD=\tbD$, with corresponding coefficients  $\bbw=\tbw$. Then, the algorithm sequentially removes dictionary elements from dictionary $\tbD$, yielding a sparse approximation $f$ of $\tilde{f}$, until the error threshold $\|f - \tilde{f} \|_{\ccalH} \leq \eps $ is violated, in which case it terminates. 
	
	Now we explain the execution of the KOMP algorithm in detail. At each stage of KOMP, a {pair of dictionary element associated with index} $j$ of $\bbD$ is selected to be removed which contributes the least to the Hilbert-norm approximation error $\min_{f\in\ccalH_{\bbD\setminus \{j\}}}\|\tilde{f} - f \|_{\ccalH}$ of the original function $\tilde{f}$, when dictionary $\bbD$ is used. Since at each stage, the kernel dictionary is fixed, this amounts to a computation involving weights {$\bbw \in \reals^{(M-1)}$ onl}y; that is, the error of removing dictionary point $\bbd_j$ is computed for each $j$ as 
	$\gamma_j ={\min_{\bbw_{\ccalI \setminus \{j\}}\in\reals^{({M}-1)}}} \|\tilde{f}(\cdot) - \sum_{k \in \ccalI \setminus \{j\}} w_k  {\kappa}(\bbd_k, \cdot) \|.$
	We use the notation $\bbw_{\ccalI \setminus \{j\}}$ to denote the entries {of $\bbw\in \reals^{M}$} restricted to the sub-vector associated with indices $\ccalI \setminus \{j\}$. Then, we define the dictionary element which contributes the least to the approximation error as $j^\star=\argmin_j \gamma_j$. If the error associated with removing this kernel dictionary element exceeds the given approximation budget $\gamma_{j^\star}>\eps$, the algorithm terminates. Otherwise, this dictionary elements {associated with} $\bbd_{j^\star}$ are removed, the weights $\bbw$ are revised based on the pruned dictionary as $\bbw = {\argmin_{\bbw \in \reals^{{M}}}} \lVert \tilde{f}(\cdot) - \bbw^T\boldsymbol{ {\kappa}}_{\bbD}(\cdot) \rVert_{\ccalH}$, and the process repeats as long as the current function approximation is defined by a nonempty dictionary. 
	\footnote{{We assume that the output of Algorithm \ref{alg:komp} has bounded Hilbert norm, which may be explicitly enforced, for instance, by thresholding the norm of the weight vector if it climbs above some large constant.}}
	
	\begin{algorithm}[h]
		\caption{Destructive Kernel Orthogonal Matching Pursuit (KOMP) \hspace{-2mm}}
		\label{alg:komp}
		\begin{algorithmic}
			\Require  function $\tilde{f}$ defined by dict. {$\tbD \in \reals^{p \times \tilde{M}}$}, coeffs. {$\tbw \in \reals^{\tilde{M}/2}$}, approx. budget  $\epsilon > 0$ 
			\State \textbf{initialize} $f=\tilde{f}$, dictionary $\bbD = \tbD$ with indices $\ccalI$, model order $M=\tilde{M}$, coeffs.  $\bbw = \tbw$.
			\While{candidate dictionary is non-empty $\ccalI \neq \emptyset$}
			{\For {{$j=1,\dots,\tilde{M}$}}
				\State Find minimal error with element $\bbd_j$ removed as
				$\gamma_j = \min_{\bbw_{\ccalI \setminus \{j\}}{\in\reals^{M-1}}} \|\tilde{f}(\cdot) - \sum_{k \in \ccalI \setminus \{j\}} w_k {\kappa}(\bbd_k, \cdot) \|_{\ccalH} \; .$ \vspace{-0mm}
				\EndFor}
			\State Find index minimizing error: $j^\star\! =\! \argmin_{j \in \ccalI} \gamma_j$
			\INDSTATE{{\bf{if }}  minimal error exceeds threshold $\!\gamma_{j^\star}\!\! > \! \epsilon$}
			\INDSTATED{\bf stop} 
			\INDSTATE{\bf else} 
			
			\INDSTATED Prune dictionary $\bbD\leftarrow\bbD_{\ccalI \setminus \{j^\star\}}$, {remove both the columns associated with index $j^\star$}
			\INDSTATED Revise set $\ccalI \leftarrow \ccalI \setminus \{j^\star\}$ and { order ${M} \leftarrow {M}-1$.}
			\INDSTATED Update weights $\bbw$ defined by current dictionary $\bbD$ using 
			\vspace{-0mm}$\bbw = {\argmin_{\bbw \in \reals^{{M}}}} \lVert \tilde{f}(\cdot) - \bbw^T\boldsymbol{{\kappa}}_{\bbD}(\cdot) \rVert_{\ccalH}$\vspace{-0mm}
			\INDSTATE {\bf end}
			\EndWhile
			
			\Return ${f},\bbD,\bbw$ of complexity $\!M\!\! \leq  \!\tilde{M}$ s.t. $\!\|f \!-\! \tilde{f} \|_{\ccalH}\!\leq \!\eps$
		\end{algorithmic}
\end{algorithm}} 

With the subspace projection procedure in place, we may summarize the key steps of the proposed method in Algorithm \ref{alg:colk} for solving \eqref{eq:main_problem} while maintaining a finite model order, thus breaking the ``curse of kernelization" for compositional stochastic programming over an RKHS. The method, Compositional Online Learning with Kernels (COLK), executes the stochastic subspace projection of a functional stochastic quasi-gradient step onto sparse subspaces $\ccalH_{\bbD_{t+1}}$ stated in \eqref{eq:projection_hat}. The initial function is set to null $f_0=0$, meaning that it has empty kernel dictionary $\bbD_0=[]$ and coefficient vector $\bbw_0=[]$. The notation $[]$ is used to denote the empty matrix or vector respective size {$p \times0$ or $0$}. Then, at each step, given an independent training example $\bbxi_t$ and step-size $\eta_t$, we update the auxiliary variable \blue{$\bbg_{t+1} = (1-\beta)(\bbg_t - \boldsymbol{\ayche}_{\bbxi_t}( f_{t-1}(\bbxi_t))) + \boldsymbol{\ayche}_{\bbxi_t}( f_t(\bbxi_t))$}. Then, this updated scalar is used to compute the \emph{unconstrained} functional stochastic gradient iterate $\tilde{f}_{t+1}={(1-\lambda\alpha)f_t-{\alpha\langle{\boldsymbol{\ayche}'_{\bbxi_t}( f_t(\bbxi_t))},\ell'_{\bbtheta_t}(\bbg_{t+1})\rangle \kappa(\bbxi_t,\cdot)}}$, which admits the parametric representation $\tbD_{t+1}$ and $\tbw_{t+1}$ as stated in \eqref{eq:param_tilde}. These parameters are then fed into KOMP with approximation budget $\eps$, such that
$$(f_{t+1}, \bbD_{t+1}, \bbw_{t+1})= \text{KOMP}(\tilde{f}_{t+1},\tilde{\bbD}_{t+1}, \tilde{\bbw}_{t+1},\eps)\; .$$
Next, we establish conditions and parameter selections under which Algorithm \eqref{alg:colk} solves the compositional problem in \eqref{eq:main_problem} to near-optimality while ensuring that the memory of the learned function remains finite.

%

\section{Convergence Analysis}\label{sec:convergence}
We establish the convergence of Algorithm \ref{alg:colk}, which is characterized in terms of the behavior of the function iterates $f_t$ with respect to the minimizer of \eqref{eq:main_problem}, as well as its finite-memory properties. \blue{Our analysis is inspired by the notion of coupled supermartingales illuminated in \cite{wang2017stochastic}. There are two key points of departure in our analysis: (i) the projection step introduces a finite error to the gradient estimation at each time step in pursuit of low-dimensional function representations which then determines the $\epsilon$-factor in the radius/rate of convergence; (ii) the use of a difference-of-iterates scheme in the faster time-scale \eqref{track1} yields tighter dependences on problem constants in the convex setting, and improved rates to stationarity in non-convex settings. The specific effect of (ii) is summarized in Lemma \ref{lemma:per_iterate_td_average}, which is key to deriving these results. }

To proceed with the analysis, a few quantities are first defined to simplify the exposition. Firstly, define  the  
functional stochastic quasi-gradient of the objective function (cf. \eqref{regularized_loss})
\begin{equation}\label{eq:quasi_sg}
	\!\!\!\hat{\nabla}_f R(f_t, \bbg_{t+1}; \bbxi_t, \bbtheta_t )\!=\!{{\langle{\boldsymbol{\ayche}'_{\bbxi_t}(\! f_t(\bbxi_t))},\ell'_{\!\bbtheta_t}(\bbg_{t+1})\rangle \kappa(\bbxi_t,\!\cdot\!)}}+\!\lambda f_t \; ,
\end{equation}
and corresponding compressed projected version is written as 
\begin{equation}\label{eq:projected_quasi_sg}
	\!\tilde{\nabla}_{\!f} \!R(f_t, \bbg_{t+1}; \bbxi_t, \bbtheta_t )\!=\!
	\frac{f_{t} \!-\! \ccalP_{ \ccalH_{\bbD_{t+1}}} \![ f_t \!-\! \alpha \hat{\nabla}_f R(f_t, \bbg_{t+1}; \bbxi_t, \bbtheta_t )]}{\alpha}.
\end{equation}
Utilizing this notation, the main function update step in \eqref{eq:sgd_tilde} can be written as 
\begin{equation}\label{eq:quasi_projected_fsgd_projected}
	{f}_{t+1} =f_t - \alpha \tilde{\nabla}_f R(f_t, \bbg_{t+1}; \bbxi_t, \bbtheta_t )\; .
\end{equation}
which is an update for the function $f_{t+1}$ re-written as a step involving projected gradients. These definitions are needed to clarify the technical setting, which we do next. Begin by defining as $\ccalF_t$  the filtration, i.e., the time-dependent sigma algebra containing the algorithm history $\ccalF_t \supset (\{f_u, g_u\}_{u=0}^t  \cup \{\bbtheta_s, \bbxi_s\}_{s=0}^{t-1})$. 
With this definition of algorithm history contained by $\ccalF_t$, we state the assumptions as follows. 
\begin{assumption}\label{kernel_bound}
	The reproducing kernel map is bounded as 
	\begin{align}\label{assum1}
		\sup_{\bbu\in\mathcal{U}}\kappa(\bbu,\bbu)= U^2 <\infty\; .
	\end{align}
\end{assumption}
\begin{assumption}\label{second_moments}
	At each time instant $t$, the second moment of the derivative of inner function ${{\boldsymbol{\ayche}'}_{\bbxi_t}(f(\bbxi_t))}$ and outer function ${{\ell'_{\bbtheta_t}\left(\bbg_{t+1}\right)}}$ is bounded as 
	\begin{align}\label{assum2_1}
		\!\!\!	\mathbb{E}\!\left[\!|{{{\boldsymbol{\ayche}'}_{\bbxi_t}( f(\bbxi_t))}}|^2\!\given\! \bbtheta_t\right]\!\leq G_{{\ayche}}  \;,   \mathbb{E}\Big\{|{{\ell'_{\bbtheta_t}\left(\bbg_{t+1}\right)}}|^2\!\given\! \ccalF_t \!\Big\} \!\leq G_\ell\; ,
	\end{align}
	with $G_{{\ayche}}$ and $G_\ell$  as finite constants with probability $1$ (w.p.$1$). 
\end{assumption}
\begin{assumption} \label{lips_outer}
	The instantaneous derivative of the outer function ${\ell_{\bbtheta}(\cdot)}$ is Lipschitz continuous with respect to its first scalar argument \blue{for each realization of $\bbtheta$} so that we may write
	{	\begin{align}\label{assum3}
			| {\ell'_{\bbtheta}\left(\bbu\right)}- {\ell'_{\bbtheta}\left(\bbv\right)}|\leq L_\ell\|\bbu-\bbv\|\; \blue{ \text{ w.p. } 1.}
	\end{align} }
\end{assumption}
\begin{assumption}\label{lips_inner} 
	The functions $\h$ is Lipschitz continuous and has bounded variance, i.e., 
	\begin{align}\label{assum4}
		\|\boldsymbol{\ayche}_{\bbxi}( f) - \boldsymbol{\ayche}_{\bbxi}( f')\| &\leq L_{{\ayche}} \|f-f'\|_{\ccalH}  \blue{\text{ w.p. } 1.}
	\end{align}
	for all $f, f'\in\mathcal{H}$. In addition, define ${\bbdelta_t:=\boldsymbol{\ayche}_{\bbxi_t}(f_t;\bbxi_t)}$ with $\bar{\bbdelta}_t= \mathbb{E}\left[\bbdelta_t \given \bbtheta_t\right]$.  Then $\bbdelta_t$ has finite variance as
	%
	$\mathbb{E}[\|\bbdelta_t - \bar{\bbdelta}_t\|^2\given \ccalF_t ] \leq \sigma_{\boldsymbol{\delta}}^2\; \blue{ \text{ w.p. } 1.}$
	%
\end{assumption}
\begin{assumption} 
	(Bounded Gradients)\label{assumption_gradient_bounded}
	The instantaneous gradient of the outer function {${\ell'_{\bbtheta}\left(\bbu\right)}$} is \blue{almost surely} bounded 
	$|{\ell'_{\bbtheta}\left(\bbu\right)}|\leq C_{\ell}$. 
\end{assumption}
%


Assumption \ref{kernel_bound} follows from the compactness of the feature space $\bbXi \cup \bbTheta$. Assumption \ref{second_moments} is regarding the second moments of the derivatives which limits the variance of the stochastic approximation error, and is typical in the literature \cite{nemirovski2009robust}. Assumption \ref{lips_outer} and \ref{lips_inner} regarding the Lipschitz continuity of the outer and inner function holds for most applications, and holds for most differentiable functions. \blue{Note that Assumption \ref{assumption_gradient_bounded} is automatically satisfied for most smooth convex functions when data domains are compact. } These assumptions are standard and hold for the applications in the next section. 

With the technical setting clarified, we are ready to state our theoretical results. That is, first we bound the model order (dictionary size) of the proposed algorithm in terms of the step size $\alpha$  and compression budget $\epsilon$. \blue{Specifically, we extend \cite{koppel2016parsimonious}[Theorem 3] to compositional settings, and further provide a non-asymptotic dependence of the function complexity on problem constants and the parameter dimension.  }

\begin{table}[t!]
	\centering
	{\begin{tabular}{|l|l|l|l|l|}
			\hline
			Refs. & Structure & Obj. &   Model compl. & Convergence rate\\
			\hline
			\cite{wang2012breaking} & Non-comp.& SC&   $\Omega(1)$  &  $\mathcal{O}\left({1}/{(\delta-E)^2}\right)$ \\
			\cite{dai2014scalable} & Non-comp.& SC&   $\Omega(T)$  & $\mathcal{O}\left({1}/{\delta}\right)$ \\
			%
			%
			%
			\cite{koppel2016parsimonious} & Non-comp.& SC &   $\Omega(T^p)$  & $\mathcal{O}\left(({1}/{\delta^2})\right)$  \\
			\cite{li2019towards} & Non-comp. &SC&    $\Omega(\sqrt{T}\log(d_\lambda))$  &   $\mathcal{O}\left({1}/{\delta^2}\right)$ \\ 
			\textbf{Our} & Comp. & SC&    $\Omega\left(\left({1}/{\delta}\right)^{p}\right)$  & $\mathcal{O}\left(({1}/{\delta})\log\left({1}/{\delta}\right)\right)$\\
			\textbf{Our} & Comp. & NC&    $\Omega\left(\left({1}/{\delta}\right)^{p}\right)$  & $\mathcal{O}\left({1}/{\delta^2}\right)$\\
			\hline
	\end{tabular}}
	\caption{{ Summary of recent related results of model complexity lower bounds and convergence rate upper bounds where we differentiate as follows: objective structure is compositional (comp.) or non-compositional., objective function is strongly convex (SC) or non-convex (NC). In \cite{wang2012breaking}, $E$ denotes the average error in the gradient, and in \cite{dai2014scalable}, $d_\lambda$ denotes the degree of freedom for a given regularization parameter $\lambda$. In the table, $\delta$ denotes the optimal accuracy parameter.} }
	\label{tab12}
\end{table}

\begin{theorem}\label{model_order} \blue{Consider constant step sizes $\alpha$, $\beta$, with finite compression budget $\epsilon$, and regularization parameter as $\lambda =   L_\ell U^2 G_h^2 \alpha/\beta + \lambda_0=\mathcal{O}(\alpha/\beta+1)$ with $\lambda_0<\min\{1,1/\alpha\}$. Let $M_t$ denote the model order of the function iterate $f_t$ i.e. number of columns in current dictionary $D_t$, then under Assumptions \ref{kernel_bound}-\ref{assumption_gradient_bounded},  $M_t$ satisfies
		\begin{align}\label{model_order_0}
			1\leq M_t\leq \mathcal{M} := \mathcal{O}\left(\frac{\alpha}{\epsilon}\right)^{2p}.
	\end{align}}	
\end{theorem}\vspace{-0mm}
\blue{Theorem \ref{model_order} (see See Appendix \ref{model_order_control} for proof) ensures that the number of data points in the kernel representation of $f_t$ generated from Algorithm \ref{alg:colk} is finite and dependent upon the ratio of the step-size to the compression parameter. The condition for online sparsification performed by KOMP algorithm boils down to the condition for packing number of the kernelized feature space $\phi(\mathcal{X})$ as described in \eqref{eq:min_gamma30}. The packing of kernelized feature space is inversely proportional to the radius $\sqrt{\frac{\epsilon}{\alpha C_{\ell}L_{\ayche}}}$. As this radius increases, the packing number reduces, meaning the model order required to cover feature space decreases. A larger radius may be attained by choosing a larger $\epsilon$, meaning that fewer points are required to cover the data domain, and thus yields a lower model order, at the cost of incurring greater sub-optimality.}

{As presented in Theorem \ref{model_order}, a similar gaurantees exist in the literature but for the non-compositional objective functions. The authors in \cite{wang2012breaking} used a fixed memory size for the function representation but it leads to a constant term as the gradient error in the convergence rate analysis. An online algorithm for streaming data points is proposed in \cite{dai2014scalable}  but the memory requirement is linear in $T$. Another functional gradient based technique in \cite{koppel2016parsimonious} achieves a model complexity of $\Omega(T^p)$. Further, a lower bound on the number of required features for a random feature based kernel methods is established in \cite{li2019towards} to be $\Omega(\sqrt{T}\log(d_\lambda))$. Here, $d_\lambda$ is the effective dimension if the problem and depends upon the kernel matrix. It is difficult to calculate $d_\lambda$ in the streaming data applications, and hence deciding the number of features for the problem is difficult. All the existing results are summarized in Table \ref{tab12}. A key upshot of our analysis is that it is the first to provide complexity bounds for RKHS function approximations for compositional objectives, both for strongly convex (SC) and non-convex (NC) cases.}

\blue{We begin establishing a basic property of Algorithm \ref{alg:colk}, which is inspired by, but distinct from  \cite{wang2017stochastic}. In particular, the gradient tracking scheme \eqref{track1} enables us to obtain tighter bounds on the estimation error of  the sequence $\{\bbg_t\}$ presented as follows.}
\blue{\begin{lemma}
		\label{lemma:per_iterate_td_average} 
		Consider  ${\bbdelta_t:=\boldsymbol{\ayche}_{\bbxi_t}(f_t;\bbxi_t)}$ and $\bar{\bbdelta}_t= \mathbb{E}\left[\bbdelta_t \given \bbtheta_t\right]$ as defined in Assumption \ref{lips_inner}. Then the estimation error of the auxiliary sequence {$\bbg_t$} [cf. \eqref{track1}] with respect to $\bar{\bbdelta}_t $ satisfies
		\begin{align} \label{eq:per_iterate_td_average}
			\E{\norm{\bbg_{t+1} - \bar{\boldsymbol{\delta}}_{t}}^2} \leq& (1-\beta)^2\E{\norm{\bbg_t - \bar{\boldsymbol{\delta}}_{t-1}}^2}
			\nonumber 
			\\
			&+ 2L_{{\ayche}}^2\E{\norm{f_t - f_{t-1}}_{\mathcal{H}}^2}+ 2\beta^2\sigma_{\boldsymbol{\boldsymbol{\delta}}}^2
		\end{align}
		{where $\beta\in(0,1)$ is a step-size parameter.}
\end{lemma}}
The proof of Lemma \ref{lemma:per_iterate_td_average} is provided in Appendix \ref{lemma1_proof}. 

\medskip
\blue{ {\noindent \bf Strongly Convex Objectives} Next, we characterize the convergence rate and the memory size requirements when the objective $R(f)$ [cf. \eqref{regularized_loss}] is strongly convex,  which in light of Theorem \ref{model_order}, yields a trade-off between the upper-bound on complexity $\mathcal{M}$ and convergence.}

\begin{theorem}\label{theorem:constant_stepsize_convergence0}
	\blue{Consider constant step-sizes $\alpha, \beta$, and constant compression budget $\eps$ with regularization parameter as $\lambda =   L_\ell U^2 G_h^2 \alpha/\beta + \lambda_0=\mathcal{O}(\alpha/\beta+1)$ with $\lambda_0<\min\{1,1/\alpha\}$. Then when $R(f)$ [cf. \eqref{regularized_loss}] is strongly convex,
		\begin{itemize}
			\item[(i)] the iterates of Algorithm \ref{alg:colk} converges to a neighborhood of the optimal $f^\star$ as
			\begin{align}\label{eq:constant_stepsize_parameters0}
				\mathbb{E}\left[\|f_{T+1} \!-\! f^\star\|_{\ccalH}^2  \right] 
				\leq&\left(1 - \alpha \lambda_0 \right)^{T}\|  f^\star\|_{\ccalH}^2 
				+Z_1\alpha +Z_2\frac{\beta^2}{\alpha}\nonumber 
				\\
				&\!\!\!\!\!+Z_1\alpha +Z_2\frac{\beta^2}{\alpha}+Z_3\frac{\epsilon}{\alpha}+Z_4\frac{\epsilon^2}{\alpha}
			\end{align}
			when run with constant step-sizes $\alpha>0$ and $0<\beta<1 $, with  problem-dependent constants defined as $Z_1=(\sigma_f^2+2 L_\ayche^2 L_{\ell}U^2\sigma_f^2)/\lambda_0$, $Z_2\!=\!(2 L_{\ell}U^2\sigma_{\boldsymbol{\delta}}^2)/\lambda_0$, $Z_3=(2Z)/\lambda_0$, and $Z_4\!=\!4 L_\ayche^2 L_{\ell}U^2/\lambda_0$. \label{convex_theorem_i}
			
			\item[(ii)] Moreover, to achieve $\delta$ sub-optimality, we require at least $T$ samples and an model complexity upper bound $\mathcal{M}$ as
			\begin{align}
				\!\!\!\!T\!\geq \mathcal{O\!}\left(\!\frac{1}{\delta}\log\left(\!\frac{1}{\delta}\!\right)\!\right) \;, 
				\mathcal{M}\geq  \mathcal{O\!}\left(\!\!\max\!\!\Big\{\!\!\!\left(\frac{5 Z_3}{\delta}\!\!\right)^{\!\!2p}\!\!\!\!,\left(\frac{Z_1}{Z_4}\!\right)^{\!\!p}\!\!\Big\}\!\!\!\right)
			\end{align}
			
			for step-sizes $\alpha =\frac{\delta}{5Z_1}$, $\beta =\frac{\delta}{5\sqrt{Z_1 Z_2}}$, budget $\epsilon\!\leq\!\min\Big\{\!\frac{\delta Z_4^{1/2}}{Z_1^{3/2}},\frac{\delta^2}{25 Z_1Z_3}\!\Big\}$. \label{convex_theorem_ii}
			
	\end{itemize} }
\end{theorem}
See Appendix \ref{proof_lemma_constant_step_size} for proof. 
\blue{The result in statement (i) of Theorem \ref{theorem:constant_stepsize_convergence0} is comparable to well-known finite sample analysis of stochastic gradient algorithms, which converge to a $\mathcal{O}(\alpha)$ neighborhood of the optimal value for step-size $\alpha$. } 
\blue{ However, due to the complexity of RKHS parameterizations, we restrict focus to constant compression budget $\epsilon>0$ and study how controlled is the function's representational complexity. Doing so allows us to obtain specific trade-offs between complexity and convergence, whereas $\epsilon=0$ implies the complexity of the function parameterization grows unbounded with time $T$. A further point of contrast is with approaches that fix the complexity, and may cause divergence \cite{wang2010online}. }
\begin{figure*}
	\begin{subfigure}[b]{0.5\textwidth}
		\includegraphics[width=\linewidth,height=4cm]{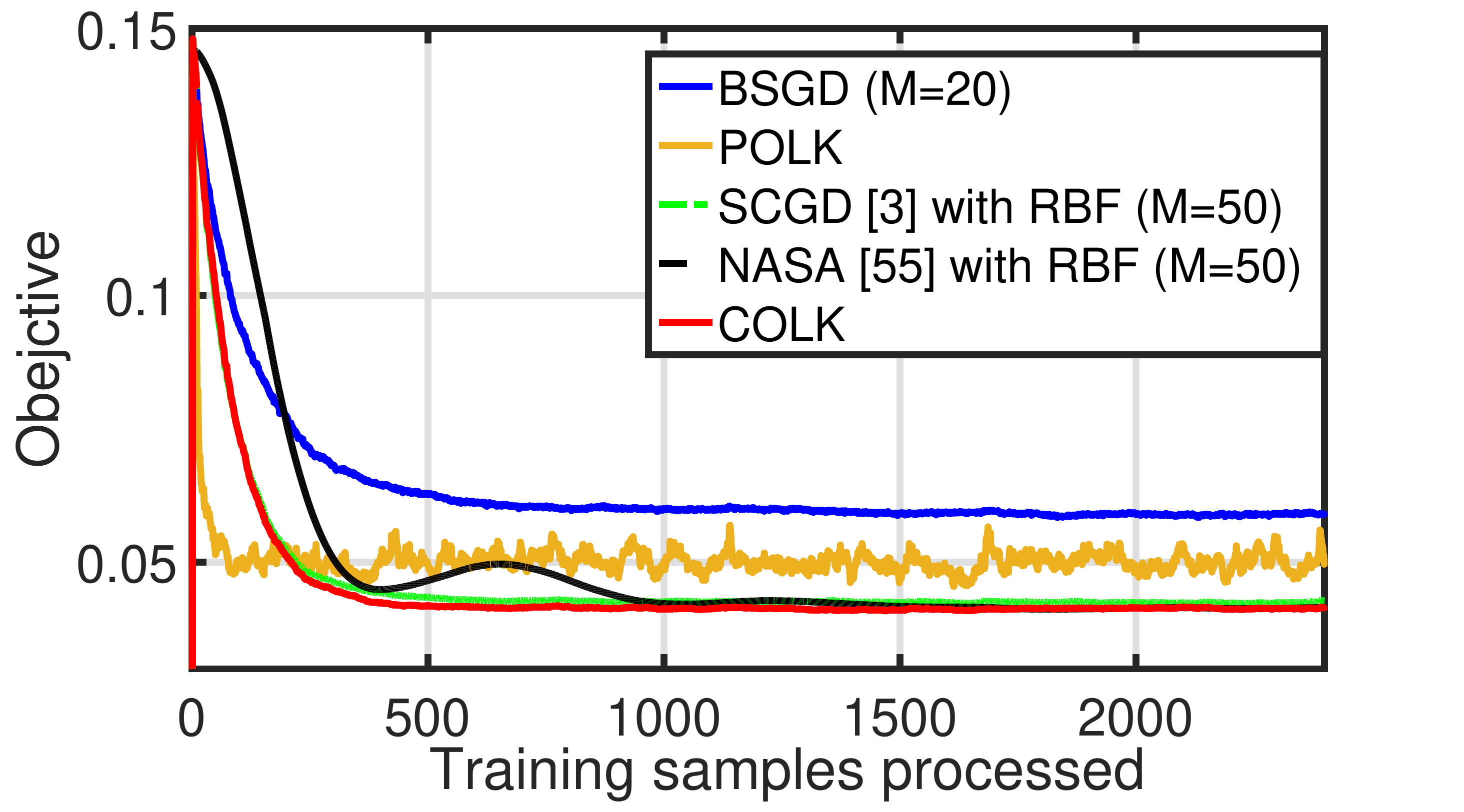}
		\caption{{Objective function}}
		\label{fig:multidist_obj}
	\end{subfigure}%
	\begin{subfigure}[b]{0.5\textwidth}
		\includegraphics[width=\linewidth,height=4cm]{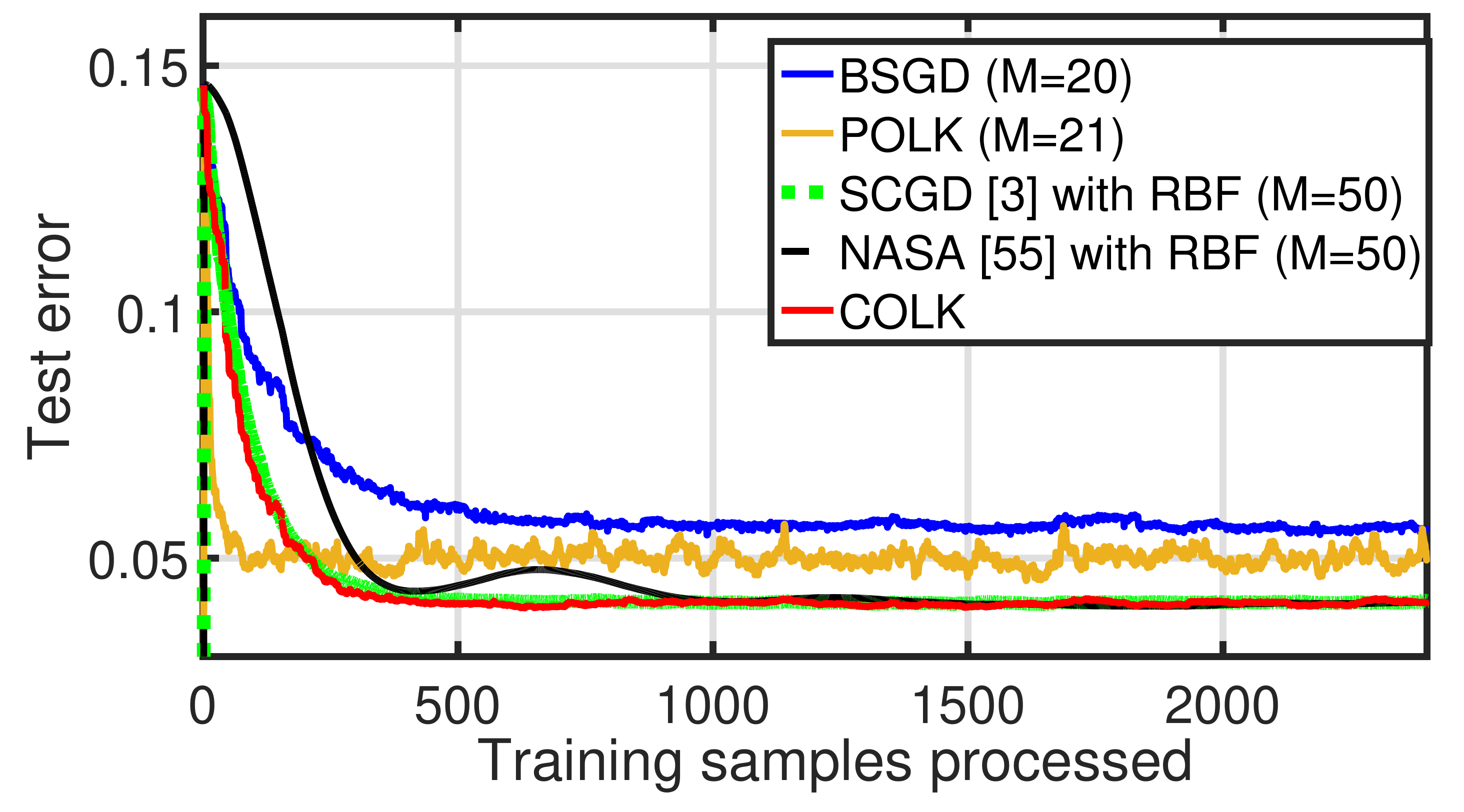}
		\caption{{Test accuracy}}
		\label{fig:multidist_err}
	\end{subfigure}\\%
	\begin{subfigure}[b]{0.5\textwidth}
		\includegraphics[width=\linewidth,height=4cm]{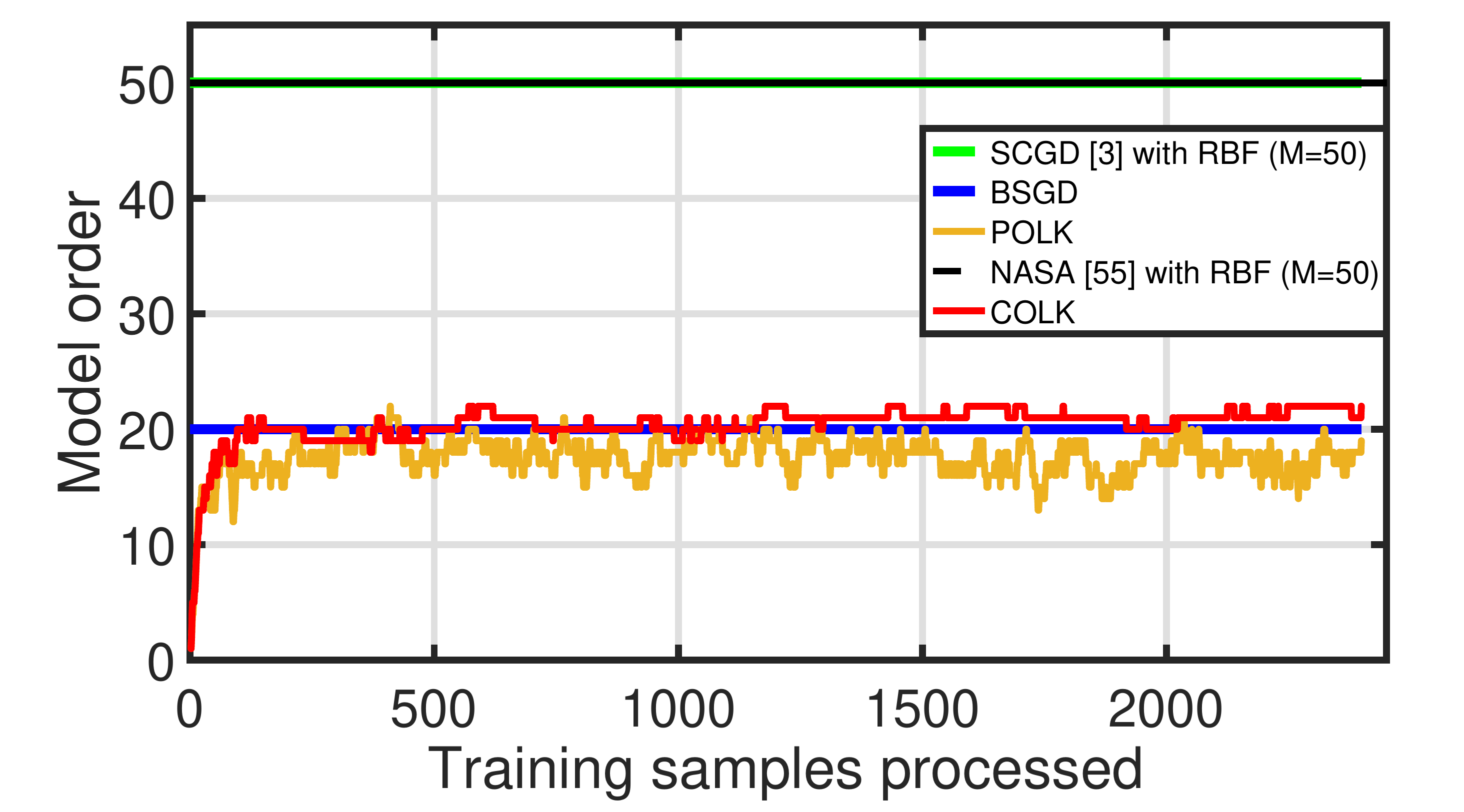}
		\caption{{Model order}}
		\label{fig:model_order}
	\end{subfigure}
	\begin{subfigure}[b]{0.5\textwidth}
		\includegraphics[width=\linewidth,height=4cm]{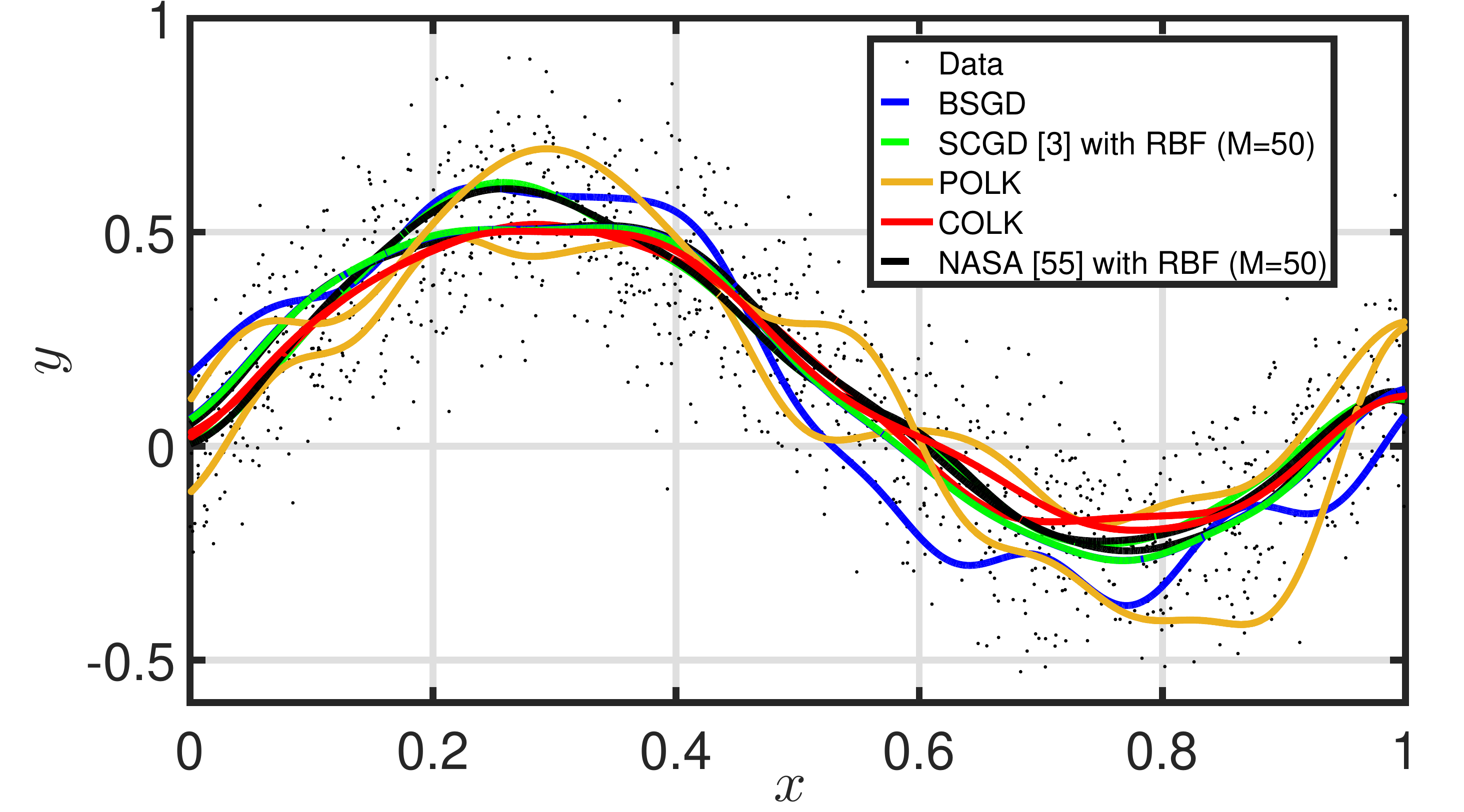}
		\caption{{Visualization of regression function plotted for two training sets}}
		\label{fig:visualization}
	\end{subfigure}%
	\vspace{-0mm}	\caption{ {COLK experimental behavior on a  regression on a synthetic data set with training outliers. Here COLK minimizes bias, variance, and third and fourth-order deviations.} }\label{fig:multidist_timeseries}\vspace{-0mm}
\end{figure*}
\blue{Instead of fixing the complexity a priori, the convergence accuracy parameter $\delta$ determines the number of iterations $T$ the algorithm must run, as well as the lower bound on the compression budget $\epsilon$, to obtain a $\delta$-suboptimal solution. Doing so then yields a specific upper bound on the model order. In doing so, we provide for the first time a non-asymptotic tradeoff between complexity and sub-optimality of solving compositional problems over an RKHS, and thus address a strictly more general setting than existing concentration bounds using random Fourier Features, whose focus is restricted to standard stochastic convex programming  \cite{sutherland2015error,sriperumbudur2015optimal,li2019towards}. We are able to address the compositional setting in contrast to these works due to the fact that our analysis only depends on function smoothness rather than tying distributional properties of the stochastic approximation error to the eigenvalue decay of the population kernel matrix, as in the aforementioned references.} 

\medskip
\blue{{\noindent \bf Non-Convex Objectives} 
	As previously noted, many practical risk measures may be non-convex -- see Example \ref{eg_risk} and \cite{ruszczynski2006optimization}. Therefore, we shift focus to analyzing algorithm performance for the setting where $R(f)$ \eqref{regularized_loss} is non-convex.

		\begin{theorem}\label{nonconvex_theorem}
			For $R(f)$ [cf. \eqref{regularized_loss}] non-convex in $f$, 
			\begin{itemize}
				\item[(i)]the function sequence $\{f_t\}$ generated by Algorithm \ref{alg:colk} run for $T$ total iterations with null initialization $f_0=0$ converges to an approximate stationary solution as
				\begin{align}\label{eq:nonconvex_theorem}
					\!\!\!\!\!\underset{1 \leq t \leq T}{\inf} \mathbb{E}\left[\left \| \nabla_f R(f_t) \right \|^2_\mathcal{H}\right] \leq \mathcal{O}\left(\frac{1}{T \alpha}+\frac{\beta^2}{\alpha} + {\alpha}+\frac{\epsilon^2}{\alpha^2}\right)
				\end{align} 
				under step-sizes $\alpha>0$, $0<\beta<1$, and budget $\epsilon$.  \label{nonconvex_theorem_i}
				
				\item[(ii)] Thus, to achieve an $\delta$-approximate stationary solution, the iteration complexity is $T\geq \mathcal{O}\left(\frac{1}{\delta^2}\right)$ and $\mathcal{M}\geq \mathcal{O}\left(\frac{1}{\delta^p}\right)$ with step size $\alpha=\mathcal{O}(\delta)$, $\beta=\mathcal{O}(\delta)$, and the compression budget $\epsilon\leq\mathcal{O}(\delta^{3/2})$. \label{nonconvex_theorem_ii}
			\end{itemize} 
		\end{theorem}

		The proof of this theorem is provided in Appendix \ref{nonconvex_proof}.
		Observe that for $\epsilon=0$, $\mathcal{M}=T$  which grows unbounded with time since we add each new sample to the kernel dictionary, and the right-hand side of \eqref{eq:nonconvex_theorem} reduces to $\mathcal{O}{(T^{-1/2})}$ because $p\geq 1$ and matches existing results \cite{ghadimi2018single}. By contrast, for fixed $\epsilon>0$, we can provide a specific tradeoff between convergence accuracy to stationarity and the required complexity of the function representation (Theorem \ref{nonconvex_theorem}(ii), which does not require $\boldsymbol{\ayche}_{\bbxi}$ to be smooth. 
	Next, we investigate the empirical validity of Algorithm \ref{alg:colk} on several practical problems.}


\section{Experiments}
\label{sec:experiments}

%
To show the efficacy of the proposed algorithm, we consider a problem of nonlinear regression (filtering)  over a $p$-dimensional parameter space.  We have again have two sets of random variables $(\bbx, \bbx') \in\ccalX\subset \reals^p$ but now the target variables are real valued $y,y'\in\reals$. The merit criterion of model fitness for a given training example $(\bbx_n, y_n)$ is the humble square loss:
\begin{equation}\label{eq:square_loss1}
	\ell(f(\bbx_n),y_n) = (f(\bbx_n) - y_n)^2.
\end{equation}

We consider as a surrogate for the approximation error the $p$-th order central moments, as distributions may be completely characterized by their moments \cite{durrett2010probability}[Ch. 3]
\begin{align}\label{eq:semi_deviation}
	\! \!\!\!\mathbb{D}[\ell(f(\bbx),\bby)] \!= \!\!\sum_{p=2}^P  \! \mathbb{E}_{\bbx, \bby\!}\Big\{\! \! \Big(\!\ell(f(\bbx),\!\bby\!)\!
	- \!{\mathbb{E}_{\bbx', \bby'}\![\ell(f\!(\bbx'),\!\bby')]}\! \Big)^p\Big\} .\!
\end{align}
which are just deviations raised to the $p$-th power \cite{ahmed2006convexity}. However, since the computational overhead scales with $P$, we truncate the upper summand index in \eqref{eq:semi_deviation} to $P=4$. Note that the dispersion measure in \eqref{eq:semi_deviation} is non-convex which is used for the experimental purposes which corresponds to the variance, skewness, and kurtosis of the loss distribution. Note that \eqref{eq:semi_deviation} may be convexified through a positive projection of $(\ell(f\!(\bbx),\!\bby)\!- \!{\mathbb{E}_{\bbx', \bby'}\![\ell\!(f\!(\bbx'),\!\bby')]} ) $, in which case the standard deviation becomes a semi-deviation, as do its higher-order analogues \cite{kalogerias2018}. This is the risk function we use for experiments in correspondence to \eqref{eq:semivariance} provided for example 1 in Sec.\ref{sec:problem}. However, for simplicity, we omit the positive projection in experiments. Next, we apply the proposed algorithm to solve the nonlinear regression problem which results in the following updates. 
\begin{align}\label{eq:algorithm_filtering}
	\hspace{-3.4cm}
	g_{t+1} =& (1-\beta)(g_t - (f_{t-1}(\bbx_t')-y_t')^2) + (f_t(\bbx_t')-y_t')^2
	\\ 
	\tilde{f}_{t+1}=&(1-\lambda\alpha)f_t\!-\!\alpha \Big\{  2 (f_t(\bbx_t)- y_t)\kappa(\bbx_{t}, \cdot) +  \eta 
	\sum_{p=2}^4  m(p,g_{t+1})  \nonumber\\
	&\qquad\quad \times \big[2(f_t(\bbx_t)\!-\! y_t) {\kappa}(\bbx_{t},\!\cdot\!)\! -\!2(f_t(\bbx_t')\!-\! y_t') {\kappa}(\bbx_{t}',\!\cdot\!)\big]\Big\},\nonumber
\end{align}
where $m(p,g_{t+1}):=\left[p ((f_t(\bbx_t)- y_t)^2  - g_{t+1})^{p-1} \right]$. Note that the gradient  of the outer function is not bounded in the above mentioned problem, but the gradient of the projected version after applying KOMP is bounded.

We evaluate COLK on synthetic and real data sets whose distributions are skewed or heavy-tailed, and compare its test accuracy against existing benchmarks that minimize only bias. Firstly, we evaluate performance on synthetic data \texttt{regression
	outliers} which has a heavier tailed distribution, i.e., more outliers are present. We inquire as to which methods overfit versus learn successfully: COLK (Algorithm \ref{alg:colk}), or methods such as  BSGD \cite{wang2012breaking},  NPBSG \cite{le2016nonparametric}, POLK \cite{koppel2016parsimonious}. \blue{{We further consider extensions of BSGD for compositional objectives, as well as a radial basis function (RBF) network of fixed size with $50$ basis functions \cite{wang2017stochastic,ghadimi2018single}. We note that \cite{ghadimi2018single} is slightly different in that it considers Jacobi rather than Gauss-Seidel-style updates for the auxiliary and main variables, in contrast to SCGD, which incurs fewer calls to a simulation oracle.}. These basis functions are  $\{\phi_i(\bbx_j)+n_j\}$ for $j=1$ to $N$ data points $\bbx_j$. Our goal in doing so is to elucidate the relationship between a fixed basis representation and the one obtained by the projection step of COLK. Of course, the nonlinearity could be defined by polynomials, sinusoidals, etc. \cite{orr1996introduction}. We consider the RBF for uniformity of experimental comparison, which takes the form  $\phi_i(x)=\exp(-(x-c_i)^2/r^2)$ where $c_i$ are the uniformly distributed centers points across the feature space $\ccalX$ with bandwidth $r=0.06$}

We generate $20$ different sets from the same data distribution and then run both POLK and COLK to learn a regression function. To generate the synthetic dataset \texttt{regression
	outliers}, we used the function $y=2x+3\text{sin}(6x)$ as the original function and target $y$'s observed after adding a zero mean Gaussian noise to $2x+3\text{sin}(6x)$. First we generate $60000$ samples of the data, and then select $20\%$ as the test data set. From the remaining $4800$ samples, we select $50\%$ at random to generate $20$ different training sets. We run COLK over these training set with the following parameter selections: a Gaussian kernel with bandwidth $\sigma=.06$, step-size parameters $\alpha=0.02$, $\beta=0.01$, $\epsilon=K\alpha^2$ with parsimony constant $K=5$, variance coefficient $\eta=0.1$, and mini-batch size of $1$. Similarity, for POLK we use $\alpha=0.5$ and  $\epsilon=K\alpha^2$ with parsimony constant $K=0.09$. We fix the kernel type and bandwidth across the different methods, and the parameters that define comparator algorithms are hand-tuned to optimize performance with the restriction that their model complexity is comparable to each other, \blue{with the exception that we choose the number of RBF features in the RBF network to be what is required to obtain comparable performance to COLK.}  We run these algorithms for different realizations of training data and evaluate their test accuracy as well as its standard deviation.

We present \blue{the convergence behavior } of Algorithm \ref{alg:colk} as well as these alternatives on the \texttt{regression outliers} data in Fig \ref{fig:multidist_timeseries}. Specifically, Fig. \ref{fig:multidist_obj} shows that the mean plus variance of the loss function is minimized as the number of samples processed increases. \blue{Observe that COLK obtains comparable sub-optimality to the RBF network, which is smaller than the comparators.}
%
Moreover, the time-series of the test-set error of COLK is given in Fig.~\ref{fig:multidist_err} which converges as the training samples increases. In Fig.~\ref{fig:model_order} we plot the model order of the function sequence defined by COLK, and observe it stabilizes over time regardless of the presence of outliers. 

\begin{figure}[t]
	\centering
	\includegraphics[width=0.7\linewidth,height=3.5cm]{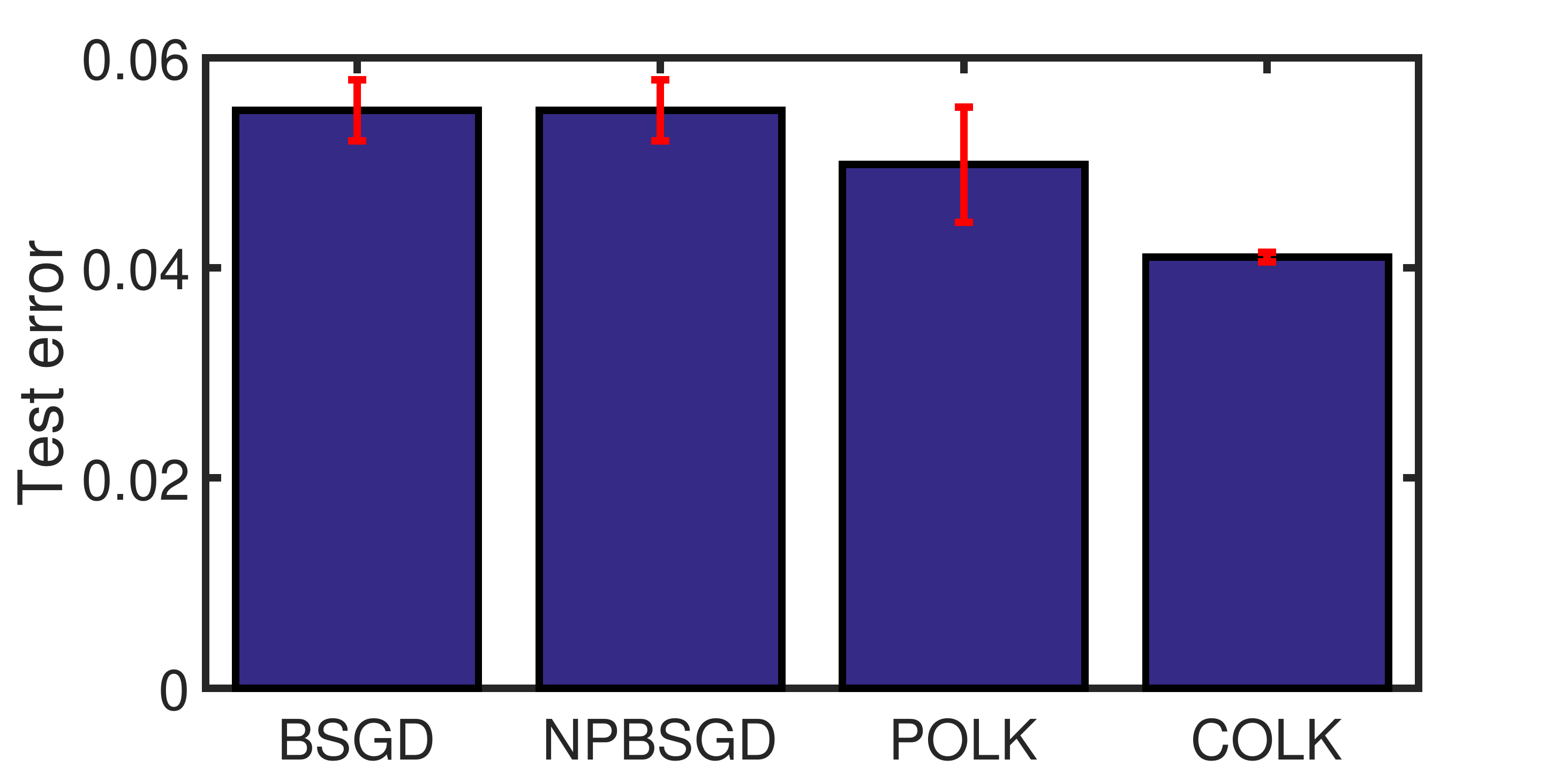}
	\caption{\small Statistical Error Comparison of COLK, with $\alpha=0.02$, $\eps=\alpha^2$, $\beta=0.01$, $K=5$, $\eta=0.1$, bandwidth $c=.06$ as compared to other methods for online learning with kernels that only minimize bias on the \texttt{regression outliers} data. This figure reports test error averages over $20$ training runs, and we report the standard deviation of test error as error bars. COLK yields both a minimal error rate and variability. The fixed model order used for the BSGD and NPBSGD is approximately equal to the converged values for COLK.\vspace{0mm}}
	\label{fig:box}
\end{figure}

\blue{Observe that in order to obtain comparable performance to COLK, one requires more than twice as many RBF features (Fig. \ref{fig:model_order}). Moreover, BSGD and other budgeted approaches yield lower accuracy estimators. } These preliminary results validate the convergence results established in Section \ref{sec:convergence}. The advantage of minimizing the bias as well as variance is depicted in Fig.~~\ref{fig:visualization} which plots the learned function for POLK and COLK for two training data sets. It can be observed that how POLK learning varies from one training set to other while COLK is robust to this change. 
%
This trend is corroborated in terms of evaluating COLK for $20$ total training runs and reporting the average test error and the standard deviation as the box and whisker plot given in Fig.~\ref{fig:box}, where we also display the result of the comparators. Observe that COLK yields the lowest error as well as the lowest standard deviation, meaning it yields inferences that are both low bias \emph{and} low variance.

To check the proposed algorithm for real data, we consider the  performance of filtering laser scans to interpolate range to a target via the \texttt{lidar} data \cite{ruppert2009semiparametric} (with added outliers) with the results shown in Fig.~\ref{fig:LIDAR}. It is clear from the figure that the proposed algorithm is robust to outliers in the data, whereas alternatives fall short in terms of either representational efficiency or overfitting.  

\section{Conclusion}\label{sec:discussion}
In this work, we addressed compositional stochastic programming in Reproducing Kernel Hilbert Space by developing a functional generalization of the stochastic quasi-gradient method operating in parallel with greedy subspace projections.
This method, Compositional Online Learning with Kernels (COLK), converges both under attenuating and constant learning rates, and yields memory-efficient parameterizations. \blue{Different from classic stochastic quasi-gradient, we introduce a momentum scheme to tighten the estimation error of the faster time-scale.}

%
We experimentally observed both with synthetic and benchmark data that COLK applied to robust supervised learning overcomes the problem of overfitting: by accounting for error variance through coherent risk, we observed consistent performance across training runs.
In future work, we hope to investigate the use of hierarchical kernels for larger parameter problems in vision or acoustics, and thus design deep learners that do not overfit. 

\appendices

\section{Proof of Theorem \ref{model_order}}\label{model_order_control}
\blue{The proof of Theorem \ref{model_order} follows from the proof of \cite[Theorem 4]{koppel2016parsimonious} which provides an asymptotic guarantee which establishes finite model order under constant step size $
	\alpha$ for the standard convex objective (non-compositional). But in this work, we are interested in characterizing the rate at which the model order (or kernel dictionary size) varies for the COLK algorithm with the number of iterations $T$. To achieve that, the initial steps of the proofs are similar to \cite[Theorem 4]{koppel2016parsimonious} which we skip here (full derivation is provided in Sec. \ref{model_order_control_supp} of supplementary for completeness) and directly start from the inequality at iteration $t$ \begin{align}\label{eq:min_gamma_optimal_weights4_0}
		\text{dist}({\kappa}(\bbxi_t,\cdot),\ccalH_{\bbD_t}) \leq \frac{\epsilon}{\alpha|V_t |} \;,
	\end{align}
	where $\text{dist}(\kappa(\bbxi, \cdot) , \ccalH_{\bbD}) =  \|  \kappa(\bbxi,\cdot) 
	- [\bbK_{\bbD_t, \bbD_t}^{-1} \bbkappa_{\bbD_t}(\bbxi)]^T
	\bbkappa_{\bbD_t}(\cdot) \|_{\ccalH}$ and $V_t=\ell'_{\bbtheta_t}\left(\bbg_{t+1}\right) {\boldsymbol{\ayche}'_{\bbxi_t}( f_t(\bbxi_t))} $. In \cite[Theorem 4]{koppel2016parsimonious}, the inequality of \eqref{eq:min_gamma_optimal_weights4_0} is obtained  for the case when the current sample $\bbxi_t$ is not added to the kernel dictionary $D_t$. This implies that when the  distance $\text{dist}({\kappa}(\bbxi_t,\cdot),\ccalH_{\bbD_t})$ for the new sample $\bbxi_t$ is less than equal to the term $\frac{\epsilon}{\alpha|V_t |}$, then the model order satisfy $M_{t+1}\leq M_t$.   
\begin{figure}[t]
	\centering
	\includegraphics[width=0.7\linewidth,height=3.5cm]{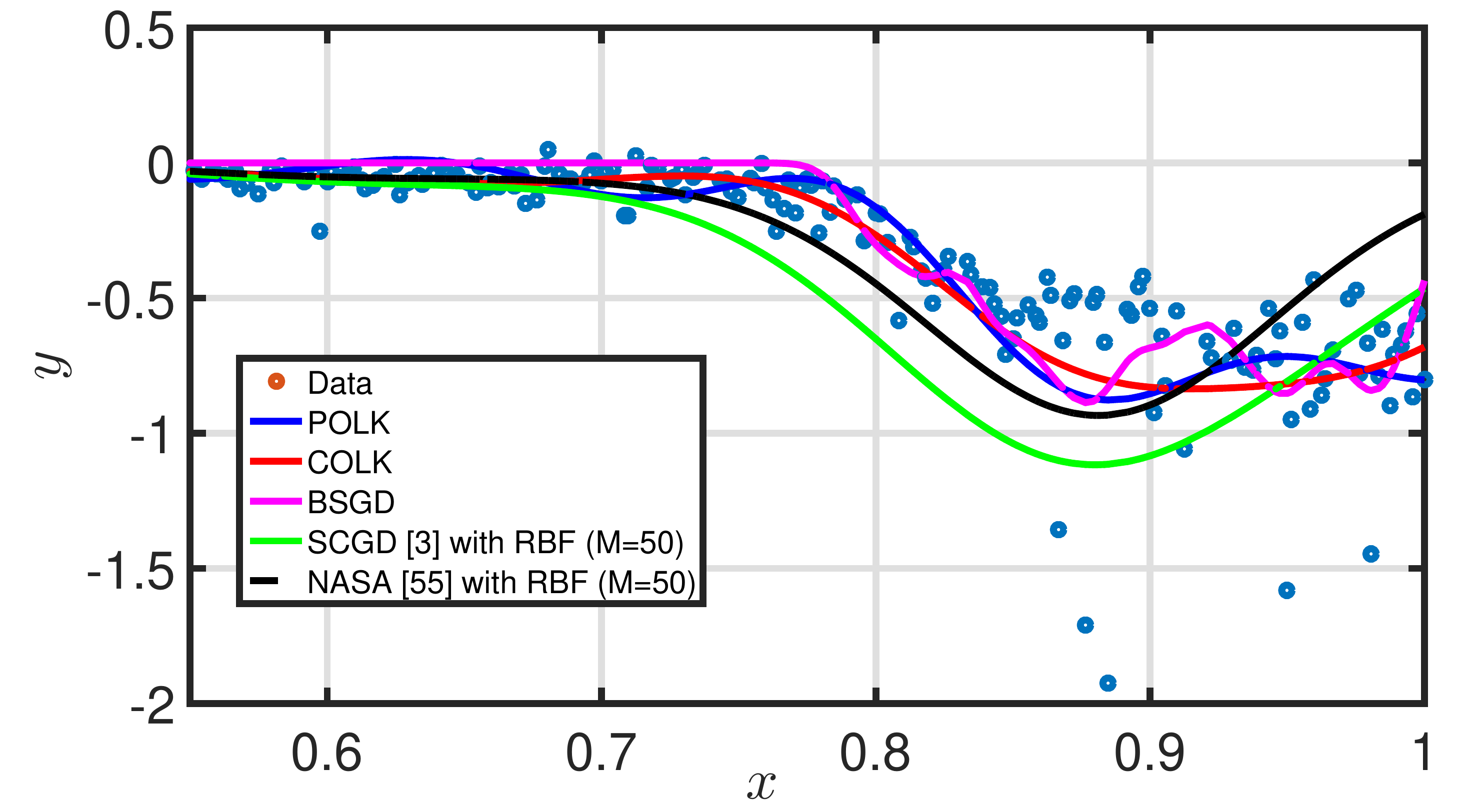}
	\caption{{Visualization of regression function plotted for LIDAR dataset}\vspace{0mm}}
	\label{fig:LIDAR}
\end{figure}

	Now, consider the contrapositive of the preceding expressions. Observe that the model order growth condition ($M_{t+1} = M_t + 1$) implies that 
	\begin{align}
		\label{eq:min_gamma2_0}
		\text{dist}({\kappa}(\bbxi_t,\cdot),\ccalH_{\bbD_t}) \geq \frac{\epsilon}{\alpha|V_t |}
	\end{align} 
	holds. This condition establishes the fact that every time a new data ($\bbxi$) is appended to kernel dictionary, then the associated product kernel is guaranteed to be at least a distance of $\blue{\frac{\epsilon}{\alpha|V_t |}}$ from every other kernel function in the current model.
	Now utilizing the Cauchy Schwartz inequality and Assumption \ref{assumption_gradient_bounded} we get $$|V_t|\leq |\ell'_{\bbtheta_t}\left(\bbg_{t+1}\right)|\| {\boldsymbol{\ayche}'_{\bbxi_t}( f_t(\bbxi_t))} \|\leq C_{\ell}L_{\ayche}.$$
	This upper bound implies that $1/|V_t| \geq 1/(C_{\ell}L_{\ayche})$, therefore we can lower bound the right hand side in \eqref{eq:min_gamma2_0} as follows
	\begin{align}\label{eq:min_gamma30}
		\frac{C{\alpha}}{|V_t|}\geq \frac{\epsilon}{\alpha C_{\ell}L_{\ayche}}\;.
	\end{align}
	From \eqref{eq:min_gamma2} (in supplementary), we obtain 
	\begin{align}
		\label{eq:min_gamma230}
		\text{dist}({\kappa}(\bbxi_t,\cdot),\ccalH_{\bbD_t}) \geq \frac{\epsilon}{\alpha C_{\ell}L_{\ayche}}.
	\end{align}
	Therefore, the KOMP stopping criterion is violated for the newest point whenever distinct dictionary points $\bbd_k$ and $\bbd_j$  for $j,k\in\{1,\dots,M_t\}$, satisfy the condition $\|\phi(\bbd_j) - \phi(\bbd_k) \|_2 > {\frac{\epsilon}{\alpha C_{\ell}L_{\ayche}}}$.  Next, we proceed in a similar manner to that of Theorem 3.1 in \cite{1315946}. Note that since the space $\mathcal{U}$ is compact and $\kappa$ is continuous, the  range $\phi(\ccalU) $ (where $\phi(\bbu)=\kappa(\bbu,\cdot)$ for $\bbu \in \ccalU$) of the kernel transformation of feature space $\ccalU$ is compact. This allows us to conclude that the number of balls of radius $\chi$ (here, $\chi = \frac{\epsilon}{\alpha C_{\ell}L_{\ayche}}$) required to completely cover the set $\phi(\ccalU)$ is finite (see, e.g., \cite{anthony2009neural}).  
	
	To prove the main result of \eqref{model_order_0}, we consider the result in \cite[Proposition 2.2]{1315946} which states that for a Lipschitz continuous Mercer kernel $\kappa$ on compact set $\mathcal{X}\subseteq\mathbb{R}^p$, for any training set $\{\bbx_t\}_{t=1}^\infty$ and any $\nu>0$, the number of elements in the dictionary is upper bounded as
	\begin{align}\label{model_order_10}
		M\leq Y\left(\frac{1}{\nu}\right)^p.
	\end{align}
	where $Y$ is a constant depends upon $\mathcal{X}$ and the kernel function. From the result in \eqref{eq:min_gamma230}, we conclude that $\sqrt{\nu}={\frac{\epsilon}{\alpha C_{\ell}L_{\ayche}}}$, which we may substitute into \eqref{model_order_10} to obtain
	\begin{align}\label{model_order_20}
		M_t\leq YC_\ell^{2p}L_\ayche^{2p}\left(\frac{\alpha}{\epsilon}\right)^{2p}=Y'\left(\frac{\alpha}{\epsilon}\right)^{2p}=\mathcal{O}\left(\frac{\alpha}{\epsilon}\right)^{2p}.
	\end{align}
	as stated in \eqref{model_order_0}. Note that in \eqref{model_order_20}, we have defined a constant $Y'=YC_\ell^{2p}L_\ayche^{2p}$.  We remark that there is a trivial lower bound of $M_t\geq 1$.  \hfill $\qed$}.

\section{Proof of Lemma \ref{lemma:per_iterate_td_average}}\label{lemma1_proof}
\blue{
	The proof employs a notion of gradient tracking for tightening the rate of convergence \cite{cutkosky2019momentum}[Lemma 2], although for an entirely different context: the error in the estimate of the inner expectation, rather than a stochastic gradient estimate error. {Whereas in \cite{cutkosky2019momentum}, the stochastic gradient sub-sampling error is the primary target hit by the auxiliary use of gradient information, instead here it is the inner-objective evaluation itself.} Begin by expanding $\bbg_{t+1} - \bar{\boldsymbol{\delta}}_t$ from \eqref{track1}, breaking the terms that multiply $1-\beta$ in two, and then add and subtract $(1-\beta)\bar{\boldsymbol{\delta}}_{t-1}$ to obtain
	\begin{align} 
		\bbg_{t+1} - \bar{\boldsymbol{\delta}}_t =& (1-\beta)(\bbg_t - \bar{\boldsymbol{\delta}}_{t-1})  - (1-\beta)(\boldsymbol{\ayche}_{\bbxi_t}( f_{t-1}(\bbxi_t))-\bar{\boldsymbol{\delta}}_{t-1}) \nonumber
		\\
		&+ \boldsymbol{\ayche}_{\bbxi_t}( f_{t}(\bbxi_t)) - \bar{\boldsymbol{\delta}}_t 
	\end{align}
	where we have gathered like terms. Now, observe that $\Et{\bar{\boldsymbol{\delta}}_{t-1} - \boldsymbol{\ayche}_{\bbxi_t}( f_{t-1}(\bbxi_t))} = \Et{\boldsymbol{\ayche}_{\bbxi_t}( f_{t}(\bbxi_t)) - \bar{\boldsymbol{\delta}}_t} = 0$ by the definition of $\delta_t$ and $\bar{\delta}_t$ in Assumption \ref{lips_inner}. Computing the norm followed by the expectation conditional on $\ccalF_t$ then yields
	\begin{align}\label{track2}
		\Et{&\norm{\bbg_{t+1} - \bar{\boldsymbol{\delta}}_t}^2}
		=(1-\beta)^2\norm{\bbg_t - \bar{\boldsymbol{\delta}}_{t-1}}^2 \nonumber
		\\
		&\!\!\!+ \Et{\|(1-\beta)(\boldsymbol{\ayche}_{\bbxi_t}( f_{t-1}(\bbxi_t))-\bar{\boldsymbol{\delta}}_{t-1})+ \bar{\boldsymbol{\delta}}_t - \boldsymbol{\ayche}_{\bbxi_t}( f_{t}(\bbxi_t)) \|^2}\nonumber
		\\
		&\!\quad-(1-\beta)\ip{\bbg_{t} - \bar{\boldsymbol{\delta}}_{t-1},(1-\beta)\Et{\boldsymbol{\ayche}_{\bbxi_t}( f_{t-1}(\bbxi_t))-\bar{\boldsymbol{\delta}}_{t-1}}\nonumber
			\\
			&\!\quad\quad+ \Et{\bar{\boldsymbol{\delta}}_t - \boldsymbol{\ayche}_{\bbxi_t}( f_{t}(\bbxi_t))}}\nonumber
	\end{align}
	where the cross term vanishes since $\bbg_t - \bar{\boldsymbol{\delta}}_{t-1}$ is independent of $\xi_t$ and the second term in the cross term is zero mean.  The second term in \eqref{track2} can again be expanded as
	\begin{align}
		&\Et{\|(1-\beta)(\boldsymbol{\ayche}_{\bbxi_t}( f_{t-1}(\bbxi_t))-\bar{\boldsymbol{\delta}}_{t-1}) + \bar{\boldsymbol{\delta}}_t - \boldsymbol{\ayche}_{\bbxi_t}( f_{t}(\bbxi_t)) \|^2} \nonumber
		\\
		&\qquad= \Et{\|(1-\beta)(\bar{\boldsymbol{\delta}}_{t-1}
			- \boldsymbol{\ayche}_{\bbxi_t}( f_{t-1}(\bbxi_t)) + \boldsymbol{\ayche}_{\bbxi_t}( f_{t}(\bbxi_t)) - \bar{\boldsymbol{\delta}}_t)\nonumber
			\\
			&\quad\qquad + \beta(\boldsymbol{\ayche}_{\bbxi_t}( f_{t}(\bbxi_t)) - \bar{\boldsymbol{\delta}}_t)\|^2} 
		\\
		&\qquad \leq 2(1-\beta)^2\Et{\|\boldsymbol{\ayche}_{\bbxi_t}( f_{t}(\bbxi_t)) - \boldsymbol{\ayche}_{\bbxi_t}( f_{t-1}(\bbxi_t))\|^2}\nonumber
		\\
		&\quad\qquad
		+ 2\beta^2\Et{\|\boldsymbol{\ayche}_{\bbxi_t}( f_{t}(\bbxi_t)) - \bar{\boldsymbol{\delta}}_t\|^2} \label{track3}
	\end{align}
	where we have used the inequality $\Ex{\norm{\mathsf{X}-\Ex{\mathsf{X}} + \mathsf{Y}}^2} \leq 2\Ex{\norm{\mathsf{X}}^2} + 2\Ex{\norm{\mathsf{Y}}^2}$ for any random variables $\mathsf{X}$ and $\mathsf{Y}$ with bounded variances. Take full expectation in \eqref{track3} and apply \eqref{assum4} to the first term on the right-hand side of the preceding expression. Now, by substituting the resulting expression into the right-hand side of \eqref{track2}, we obtain
	\begin{align}\label{track33}
		&\E{\|\bbg_{t+1} - \bar{\boldsymbol{\delta}}_{t}\|^2}  \\
		&\leq (1-\beta)^2\E{\|\bbg_t - \bar{\boldsymbol{\delta}}_{t-1}\|^2} + 2(1-\beta)^2L_{{\ayche}}^2\E{\|f_t - f_{t-1}\|_{\mathcal{H}}^2}  \nonumber\\
		&\quad + 2\beta^2\sigma_{\boldsymbol{\boldsymbol{\delta}}}^2\nonumber
		\\
		&\leq (1\!-\!\beta)^2\E{\|\bbg_t - \bar{\boldsymbol{\delta}}_{t-1}\|^2} + 2L_{{\ayche}}^2\E{\norm{f_t \!-\! f_{t-1}}_{\mathcal{H}}^2} + 2\beta^2\sigma_{\boldsymbol{\boldsymbol{\delta}}}^2 
		\label{track4}
	\end{align}
	where the last inequality \eqref{track4} follows from the fact that  $\beta\leq 1$ applied to the second term on the right-hand side of \eqref{track33}. \hfill$\blacksquare$ \\
	
}   

\section{Proof of Theorem  \ref{theorem:constant_stepsize_convergence0}}\label{proof_lemma_constant_step_size}

Before analyzing the mean convergence behavior of the distance from optimal $\mathbb{E}\left[\|f_{t+1} - f^\star\|_{\ccalH}^2  \right]$, consider the following temporal difference between {$\bbg_{t+1}$ and $\bar{\boldsymbol{\delta}}_t$} from Lemma \ref{lemma:per_iterate_td_average}
\begin{align} \label{eq:per_iterate_td_average2}
	\mathbb{E}\left[\|\bbg_{t+1} - \bar{\boldsymbol{\delta}}_t\|^2  \right] \leq& (1-\beta)^2\E{\norm{\bbg_t - \bar{\boldsymbol{\delta}}_{t-1}}^2} \nonumber
	\\
	&+ 2L_{{\ayche}}^2\E{\norm{f_t - f_{t-1}}_{\mathcal{H}}^2} + 2\beta^2\sigma_{\boldsymbol{\boldsymbol{\delta}}}^2.
\end{align}
Next, substitute the upper bound on the total expectation of $\| f_t - f_{t-1} \|_{\ccalH}^2$ as described in  Lemma \ref{lemma:iterate_relations} \eqref{lemma:value_function_difference} (\blue{stated in Appendix \ref{proof_lemma}} of the supplementary)  into \eqref{eq:per_iterate_td_average2} to obtain
\begin{align} \label{eq:per_iterate_td_average3}
	\mathbb{E}\left[\|\bbg_{t+1} - \bar{\boldsymbol{\delta}}_t\|^2  \right] \leq& (1-\beta)^2\E{\|\bbg_t - \bar{\boldsymbol{\delta}}_{t-1}\|^2} + {2 L_{\ayche}^2}[\alpha^2\sigma_f^2+ 2\epsilon^2]
	\nonumber
	\\
	&+ 2 \beta^2 \sigma_{\boldsymbol{\delta}}^2 \; ,
\end{align}
It is interesting to observe that \eqref{eq:per_iterate_td_average3} relates {$\mathbb{E}\left[\|\bbg_{t+1} - \bar{\boldsymbol{\delta}}_t\|^2  \right] $} to its previous iterate value. Utilizing this recursion, we can write the following by replacing $t+1$ by $t$
{\begin{align} \label{eq:per_iterate_td_average_previous}
		\mathbb{E}\left[\|\bbg_{t} - \bar{\boldsymbol{\delta}}_{t-1}\|^2  \right] \leq& (1-\beta)^2\E{\|\bbg_{t-1} - \bar{\boldsymbol{\delta}}_{t-2}\|^2}  \nonumber
		\\
		&+ {2 L_{\ayche}^2}[\alpha^2\sigma_f^2+ 2\epsilon^2]
		+ 2 \beta^2 \sigma_{\boldsymbol{\delta}}^2 \; ,
\end{align}}
Substituting \eqref{eq:per_iterate_td_average_previous} into the right-hand side of \eqref{eq:per_iterate_td_average3} and repeating the recursion, we can write 
%
\begin{align} \label{eq:per_iterate_td_average4}
	\mathbb{E}\left[\|\bbg_{t+1} \!-\! \bar{\boldsymbol{\delta}}_t\|^2  \right] \leq& (1\!-\beta)^{2(t+1)}{\|\bbg_{0}\! -\! \bar{\boldsymbol{\delta}}_{-1}\|^2}
	\\
	&+ \!\! \sum_{u=0}^t(1-\beta)^{2u}\Big\{\!{2 L_{\ayche}^2} [\alpha^2\sigma_f^2+ 2\epsilon^2]
	\!+\! 2 \beta^2 \sigma_{\boldsymbol{\delta}}^2 \Big\}.\nonumber
\end{align}
The first term in \eqref{eq:per_iterate_td_average4} vanishes due to the initialization $g_{0}=0$ and the convention $\boldsymbol{\delta}_{-1} = 0$. Moreover, the second term represents a geometric series and sum can be evaluated using $\sum_{u=0}^t(1-\beta)^u = [1 - (1-\beta)^t]/\beta$ provided  $\beta<1$. In this geometric sum expression, since the numerator is strictly less than unit, we can further simplify \eqref{eq:per_iterate_td_average4} to 
\begin{align} \label{eq:per_iterate_td_average_final}
	\mathbb{E}\left[\|\bbg_{t+1} - \bar{\boldsymbol{\delta}}_t\|^2  \right]\! \leq\!& \frac{2 L_{\ayche}^2}{\beta} [\alpha^2\sigma_f^2+ 2\epsilon^2]\!+\! 2 \beta \sigma_{\boldsymbol{\delta}}^2 = \ccalO\left(\!\!\frac{\alpha^2 \!+\! \eps^2}{\beta} \!+\! \!\beta\right)\; .
\end{align}
After establishing this auxiliary sequence order in terms of the order of step sizes $\alpha$ and $\beta$, we shift our focus again to the sub-optimality gap $\|f_t - f^\star\|_{\ccalH}$ in expectation. Before proceeding, note that the Hilbert-norm regularizer $(\lambda/2)\|f\|_{\ccalH}^2$ in \eqref{regularized_loss} makes the objective $R(f)$ strongly convex, i.e.
\begin{equation}\label{eq:strong_cvx}
	\frac{\lambda}{2} \|f_t - f^\star\|_{\ccalH}^2 \leq R(f_t) - R(f^\star).
\end{equation}
Using the inequality in \eqref{eq:strong_cvx} into the expression of Lemma \ref{lemma:iterate_relations}\eqref{lemma:value_function_suboptimality} (\blue{stated in Appendix \ref{proof_lemma}}), we obtain
\begin{align}\label{eq:value_function_suboptimality_strong_cvx}
	\mathbb{E}\left[\|f_{t+1} - f^\star\|_{\ccalH}^2 \given \ccalF_t \right]  &\leq \left(\!1 \!+\!  L_{\ell}U^2\frac{\alpha^2}{\beta} G_h^2-\alpha\lambda\!\right) \!\!\| f_t \!-\! f^\star\|_{\ccalH}^2  \nonumber
	\\
	&\quad\!+\! 2 \eps \|f_t \!-\!f^\star\|_{\ccalH} \!\!+\!\! \alpha^2 \sigma_f^2  \\
	&\quad\qquad\!+\!  L_{\ell}U^2\beta  \mathbb{E}\left[\|\bbg_{t+1}\!-\!\bar{\boldsymbol{\delta}}_t\|_{\ccalH}^2 \given \ccalF_t \right] \; .\nonumber
\end{align}
Now, take the total expectation of \eqref{eq:value_function_suboptimality_strong_cvx}. Considering that  $\norm{f^\star}\leq S$ which implies $\|f_t - f^\star\|_{\ccalH} \leq \|f_t\|_{\ccalH} + \|f^\star\|_{\ccalH}\leq K+S:=Z$. Now, substitute the regularizer selection $\lambda =   L_\ell U^2 G_h^2 \alpha/\beta + \lambda_0$ for $\lambda_0 < 1$, and apply \eqref{eq:per_iterate_td_average_final} to the last term on the right-hand side of the preceding expression:
\begin{align}\label{eq:mean_value_function_suboptimality}
	\mathbb{E}\left[\|f_{t+1} - f^\star\|_{\ccalH}^2  \right]  \leq &\left(1  \!-\! \alpha \lambda_0 \right)\E{ \| f_t - f^\star\|_{\ccalH}^2} 
	+ 2 \epsilon Z+ \alpha^2 \sigma_f^2     \nonumber
	\\
	&+ 2 L_{\ell}U^2\beta^2 \sigma_{\boldsymbol{\delta}}^2+   {2 L_\ayche^2 L_{\ell}U^2} [\alpha^2\sigma_f^2+ 2\epsilon^2]
	\; .
\end{align} 
Recursively substitute this expression back into itself to obtain:
\begin{align}\label{eq:mean_value_function_suboptimality222}
	\mathbb{E}&\left[\|f_{t+1} - f^\star\|_{\ccalH}^2  \right]  \leq \left(1  - \alpha \lambda_0 \right)^{t}\E{ \| f_1 - f^\star\|_{\ccalH}^2} 
	\\
	&\!\!\!\!\!\!+\!\!\! \sum\limits_{u=0}^{t-1}\!(1\!-\!\alpha\lambda_0)^u\!\Big( \!2 \epsilon Z\!\!+\!\!\alpha^2 \sigma_f^2  \!\!+\!\! 2 L_{\ell}U^2\beta^2 \sigma_{\boldsymbol{\delta}}^2 \!\!  +\!\!    2 L_\ayche^2 \! L_{\ell}U^2 [\alpha^2\sigma_f^2\!\!+\!\! 2\epsilon^2]\!\Big)\nonumber
	\; .
\end{align}
Next, since $f_1=0$ and using $\sum_{u=0}^{t-1}(1-\alpha\lambda_0)^u = [1 - (1-\alpha\lambda_0)^{t-1}]/\alpha\lambda_0\leq 1/\alpha\lambda_0$ since $\alpha\lambda_0<1$. We get the following result by substituting $t=T$ 
\begin{align}\label{eq:mean_value_function_suboptimality2222}
	\mathbb{E}\left[\|f_{T+1} - f^\star\|_{\ccalH}^2  \right] 
	\leq& \left(1  - \alpha \lambda_0 \right)^{T}\|  f^\star\|_{\ccalH}^2 
	\nonumber
	\\
	&+ \frac{1}{\alpha\lambda_0}\Big(2 \epsilon Z +\alpha^2 \sigma_f^2 + 2 L_{\ell}U^2\beta^2 \sigma_{\boldsymbol{\delta}}^2 \nonumber
	\\
	&\quad  +    2 L_\ayche^2 L_{\ell}U^2 [\alpha^2\sigma_f^2+ 2\epsilon^2]\Big)\nonumber 
	\\
	&\hspace{-3cm}= \left(1  - \alpha \lambda_0 \right)^{T}\|  f^\star\|_{\ccalH}^2 
	+Z_1\alpha +Z_2\frac{\beta^2}{\alpha}+Z_3\frac{\epsilon}{\alpha}+Z_4\frac{\epsilon^2}{\alpha}.
\end{align}
where $Z_1=(\sigma_f^2+2 L_\ayche^2 L_{\ell}U^2\sigma_f^2)/\lambda_0$, $Z_2=(2 L_{\ell}U^2\sigma_{\boldsymbol{\delta}}^2)/\lambda_0$, $Z_3=(2Z)/\lambda_0$, and $Z_4=4 L_\ayche^2 L_{\ell}U^2/\lambda_0$. \blue{From the definition of $\mathcal{M}$ in Theorem \ref{model_order}, we have $\mathcal{M}=Y'\left(\frac{\alpha}{\epsilon}\right)^{2p}$, which implies that we can write \eqref{eq:mean_value_function_suboptimality2222} as
	\begin{align}\label{eq:mean_value_function_suboptimality22222}
		\mathbb{E}\left[\|f_{T+1} - f^\star\|_{\ccalH}^2  \right] 
		\leq& \left(1  - \alpha \lambda_0 \right)^{T}\|  f^\star\|_{\ccalH}^2 +Z_1\alpha +Z_2\frac{\beta^2}{\alpha}\nonumber
		\\
		&+Z_3\left(\frac{\mathcal{M}}{Y'}\right)^{\frac{-1}{2p}}+Z_4\alpha \left(\frac{\mathcal{M}}{Y'}\right)^{\frac{-1}{p}}.
\end{align}} \blue{
	Next, in order to characterize the rates for iteration complexity as well as model complexity, it is possible to choose the step sizes $\alpha$ and $\beta$ in an optimal manner. Suppose, we require $\mathbb{E}\left[\|f_{T+1} - f^\star\|_{\ccalH}^2  \right] $ to be less than some positive value $\delta$ which describes the accuracy for the algorithm. To do so, we ensure that each terms in the right hand side of \eqref{eq:mean_value_function_suboptimality22222} is bounded by $\delta/5$. In other words, we need to satisfy the following set of inequalities 
	\begin{align}\label{ineq_first}
		&\alpha\leq \frac{\delta}{5 Z_1}, \ \beta \leq \frac{\delta}{5\sqrt{Z_1 Z_2}}, \ \ \mathcal{M}\geq Y'\left(\frac{5 Z_3}{\delta}\right)^{2p}
		\\
		& \mathcal{M}\geq Y'\left(\frac{Z_1}{Z_4}\right)^{p}\!\!\!, \ \text{and}\ -T\log(1-\alpha \lambda_0)\geq -\log\left(\frac{\delta}{5\|  f^\star\|_{\ccalH}^2}\right).\nonumber
	\end{align}
	The inequalities in \eqref{ineq_first} can be satisfied if we set $\alpha =\frac{\delta}{5Z_1}$, $\beta =\frac{\delta}{5\sqrt{Z_1 Z_2}}$ for small $\delta$ so that $-\log(1-\lambda_0 \alpha)\approx\lambda_0\alpha$ and
	\begin{align}\label{eq:possibly_incorrect_equation}
		T\geq \mathcal{O}\left(\frac{1}{\delta}\log\left(\frac{1}{\delta}\right)\right), \
		\mathcal{M}\geq \mathcal{O}(\max\Big\{\left(\frac{5 Z_3}{\delta}\right)^{2p},\left(\frac{Z_1}{Z_4}\right)^{p}\Big\})
	\end{align}
	as stated in Theorem \ref{theorem:constant_stepsize_convergence0}. From the definition of $\mathcal{M}$, we have $\mathcal{M}=Y'(\alpha/\epsilon)^{2p}$ it holds that
	\begin{align}\label{compression}
		\big(\frac{\alpha}{\epsilon}\big)^{2p}\geq \mathcal{O}\left(\max\Big\{\left(\frac{5 Z_3}{\delta}\right)^{2p},\left(\frac{Z_1}{Z_4}\right)^{p}\Big\}\right).
	\end{align}
	\blue{Using the upper bound on $\alpha$ from \eqref{eq:possibly_incorrect_equation}, we substitute $\alpha=\frac{\delta}{5 Z_1}$ into \eqref{compression}, we get 
		\begin{align}\label{compression2}
			\Big(\frac{\delta}{5\epsilon Z_1}\Big)^{2p}\geq \mathcal{O}\left(\max\Big\{\left(\frac{5 Z_3}{\delta}\right)^{2p},\left(\frac{Z_1}{Z_4}\right)^{p}\Big\}\right).
		\end{align} After rearranging the terms for $\epsilon$ in \eqref{compression2}, we obtain  $\epsilon\!\leq\!\min\Big\{\!\frac{\delta Z_4^{1/2}}{Z_1^{3/2}},\frac{\delta^2}{25 Z_1Z_3}\!\Big\}$}. 
	\hfill $\blacksquare$ }
%
%
%
\section{Proof of Theorem \ref{nonconvex_theorem}} \label{nonconvex_proof}
Note that the statement Lemma \ref{lemma:per_iterate_td_average} holds independent of the convexity of the objective. Thus, consider the statement of Lemma \ref{lemma:per_iterate_td_average} ({stated in Appendix \ref{proof_lemma} of the supplementary}) and multiply the both sides by $(1+\beta)$,  then add to the statement of Lemma \ref{nonconvex:lemma1} to obtain
\begin{align}
	\mathbb{E}\Big[& R(f_{t+1}) + \left \|\bar{\boldsymbol{\delta}}_{t}-\mathbf{g}_{t+1}\right \|^2 \mid \mathcal{F}_t\Big]\\
	&\leq \left \|\bar{\boldsymbol{\delta}}_{t-1}-\mathbf{g}_{t}\right \|^2  +2L_{{\ayche}^2}\E{\norm{f_t - f_{t-1}}_{\mathcal{H}}^2} + 2\beta^2\sigma_{\boldsymbol{\boldsymbol{\delta}}}^2 + R(f_t)\nonumber
	\\
	&- \frac{\alpha}{4}\left \| \nabla_f R(f_t) \right \|^2_\mathcal{H} + \beta \mathbb{E}[\left \|\bar{\boldsymbol{\delta}}_t-\mathbf{g}_{t+1}\right \|^2\mid\mathcal{F}_t] + \frac{\epsilon^2}{2\alpha} + \frac{L_R}{2}\sigma^2_f \alpha^2 \nonumber
\end{align}
by the observation that  $1-\beta^2 \leq 1$, $1+\beta \leq 2$. 
For the better exposition, let us define a stochastic process $\gamma_t = R(f_t) + \left \|\bar{\boldsymbol{\delta}}_{t-1}-\mathbf{g}_{t}\right \|^2$. We can write
\begin{align}
	\mathbb{E}\left[ \gamma_{t+1} \mid \mathcal{F}_t\right] \leq& \gamma_t \!+\! {2 L_\ayche}^2 \left \| f_{t}\!-\!f_{t-1} \right \|^2_\mathcal{H}\!+\!4\beta^2 \sigma^2_{\delta} - \frac{\alpha}{4}\left \| \nabla_f R(f_t) \right \|^2_\mathcal{H}\nonumber
	\\
	&+\! \beta \mathbb{E}[\left \|\bar{\boldsymbol{\delta}}_t\!-\!\mathbf{g}_{t+1}\right \|^2\mid\mathcal{F}_t] + \frac{\epsilon^2}{2\alpha} \!+\! \frac{L_R}{2}\sigma^2_f \alpha^2. 
\end{align}
Taking expectation on both sides and noting that $\mathbb{E}[\left \| f_{t}-f_{t-1} \right \|^2_\mathcal{H}] = \mathbb{E}[\alpha^2\mathbb{E}[\left \| \tilde{\nabla}_f R(f_{t-1})\right \|^2_\mathcal{H}\mid \mathcal{F}_{t-1}]]\leq \alpha^2 \sigma^2_f$, we obtain
\begin{align}
	\mathbb{E}[\gamma_{t+1}]\leq& \mathbb{E}[\gamma_{t}] + 4\beta^2\sigma^2_\delta +\frac{\epsilon^2}{2\alpha}+\ \sigma^2_f(\frac{L_R \alpha^2}{2}+{2L_\ayche^2 \alpha^2}) \nonumber
	\\
	&-\frac{\alpha}{4}\left \| \nabla_f R(f_t) \right \|^2_\mathcal{H}.
\end{align}
After rearranging the terms, we get
\begin{align}\label{bound}
	\frac{\alpha}{4}\left \| \nabla_f R(f_t) \right \|^2_\mathcal{H} \leq& \mathbb{E}[\gamma_{t}]-\mathbb{E}[\gamma_{t+1}] + 4\beta^2\sigma^2_\delta +\frac{\epsilon^2}{2\alpha}\nonumber
	\\
	&+\sigma^2_f(\frac{L_R \alpha^2}{2}+{2L_\ayche^2 \alpha^2}). 
\end{align}
Taking sum from $t=1$ to $T$, and then lower bounding the left hand side of \eqref{bound}, we obtain
%
\begin{align}
	\underset{1 \leq t \leq T}{\inf} \mathbb{E}[ \| \nabla_f R(f_t)  \|^2_\mathcal{H}] \frac{\alpha T}{4} \leq& \mathbb{E}[\gamma_{1}]-R(f^\star) + 4\beta^2\sigma^2_\delta T +\frac{\epsilon^2T}{2\alpha}\nonumber
	\\
	&+\alpha^2 \sigma^2_f T(\frac{L_R}{2}+{4L_\ayche^2}).
\end{align}
\blue{Rearranging the terms, we get
	\begin{align}
		\underset{1 \leq t \leq T}{\inf} \mathbb{E}[\left \| \nabla_f R(f_t) \right \|^2_\mathcal{H}] &\leq \frac{C_1+C_2T\beta^2+\frac{T}{2}\frac{\epsilon^2}{\alpha}+C_3T\alpha^2+C_4T{\alpha^2}}{T\frac{\alpha}{4}}\nonumber\\
		&=\mathcal{O}\left(\frac{1}{T \alpha}+\frac{\beta^2}{\alpha} + {\alpha}+\frac{\epsilon^2}{\alpha^2}\right)\label{accuracy}
	\end{align}
	where, $C_1=\mathbb{E}[\gamma_1]-R(f^\star)$, $C_2=4\sigma^2_\delta$, $C_3=\frac{\sigma^2_f L_R}{2}$, and  $C_4=4L_\ayche
	^2 \sigma^2_f$. From the definition of $\mathcal{M}$ in Theorem \ref{model_order}, we have $\mathcal{M}=Y'\left(\frac{\alpha}{\epsilon}\right)^{2p}$, which implies that we can write \eqref{accuracy} as
	\begin{align}
		\!\!\!\!\!\underset{1 \leq t \leq T}{\inf} \mathbb{E}[\left \| \nabla_f R(f_t) \right \|^2_\mathcal{H}] \leq\mathcal{O}\left(\frac{1}{T \alpha}+\frac{\beta^2}{\alpha} + {\alpha}+\left(\frac{\mathcal{M}}{Y'}\right)^{\frac{-1}{p}}\right).\label{accuracy2}
	\end{align}
	To obtain the complexity bounds presented in Theorem \ref{model_order} statement (ii), we need to consider the case when $\underset{1 \leq t \leq T}{\inf} \mathbb{E}[\left \| \nabla_f R(f_t) \right \|^2_\mathcal{H}]$ is less than the required accuracy $\delta$. From \eqref{accuracy2}, note that $\delta$ accuracy is  achieved if it holds that
	\begin{align}\label{ineq_first2}
		\frac{1}{T\alpha
		}\leq\frac{\delta}{4},\ \beta \leq \frac{\delta}{4},\ 
		\alpha\leq \delta, \ \text{and}\   \mathcal{M}\geq Y'\left(\frac{4}{\delta}\right)^{p}. 
	\end{align}
	Hence for $\alpha=\frac{\delta}{4}$, $\beta\leq \frac{\delta}{4}$, $\mathcal{M}=Y'({\alpha}/{\epsilon})^{2p} $, and from the lower bound in \eqref{ineq_first2}, it holds that 
	\begin{align}\label{compression_bound}
		\left(\frac{\alpha}{\epsilon}\right)^{2p}\geq \mathcal{O}\left(\left(\frac{4}{\delta}\right)^p\right).
	\end{align}
	Utilizing $\alpha=\frac{\delta}{4}$ and after rearranging the terms in \eqref{compression_bound}, we obtain
	$\left({\delta}/{4\epsilon}\right)^{2p}\geq \mathcal{O}\left(\left({4}/{\delta}\right)^p\right)$.
	After rearranging the terms for $\epsilon$, we obtain the bound on the compression budget $\epsilon\leq \mathcal{O}\left({\delta\sqrt{\delta}}\right)$. This implies the iteration complexity $T\geq\mathcal{O}\left({1}/{\delta^2}\right)$ and the model complexity $\mathcal{M}\geq \mathcal{O}\left({1}/{\delta^p}\right)$. }

\bibliographystyle{IEEEtran}
\bibliography{bibliography}
\onecolumn
\newpage

\section*{\centering Supplementary Material for \\Nonparametric Compositional Stochastic Optimization for Risk-Sensitive Kernel Learning}



\section{Proofs of Propositions and Technical Lemmas} \label{proof_lemma}
First, let's quantify the projection-induced error in terms of the true stochastic quasi-gradient in the following proposition.

\begin{proposition}\label{prop1}
	Given independent identical realizations $(\bbxi_t, \bbtheta_t)$ of the two associated random variables $(\bbxi, \bbtheta)$, for all $t$ it holds that  
	\begin{equation}\label{eq:prop1}
		\| \tilde{\nabla}_f R(f_t, \bbg_{t+1}; \bbxi_t, \bbtheta_t ) - \hat{\nabla}_f R(f_t, \bbg_{t+1}; \bbxi_t, \bbtheta_t )\|_{\ccalH} \leq \frac{\eps}{\alpha}\; ,
	\end{equation}
	where $\alpha>0$ denotes the step-size and $\eps>0$ is the compression parameter of Algorithm \ref{alg:komp}.
\end{proposition}

{\bf \noindent Proof:}
Consider the square-Hilbert-norm difference of $ \tilde{\nabla}_f R(f_t,\bbg_{t+1}; \bbxi_t, \bbtheta_t )  $  and $\hat{\nabla}_f R(f_t,\bbg_{t+1}; \bbxi_t, \bbtheta_t )  $ defined in \eqref{eq:quasi_sg} and \eqref{eq:projected_quasi_sg}, respectively,
\begin{align}\label{eq:norm_stoch_grad_diff}
	&\| \tilde{\nabla}_f R(f_t,\bbg_{t+1}; \bbxi_t, \bbtheta_t ) - \hat{\nabla}_f R(f_t,\bbg_{t+1}; \bbxi_t, \bbtheta_t )  \|_{\ccalH}^2  \\
	& = \Big\| \! \Big( f_{t}\! - \!\ccalP_{ \ccalH_{\bbD_{t+1}}} \Big[ f_t\! -\! \alpha \hat{\nabla}_f R(f_t,\bbg_{t+1}; \bbxi_t, \bbtheta_t ) \Big]\Big)\!  /\alpha 
	- \hat{\nabla}_f R(f_t,\bbg_{t+1}; \bbxi_t, \bbtheta_t )  \Big\|_{\ccalH}^2\; . \nonumber
\end{align}
Multiply and divide $\hat{\nabla}_f R(f_t,\bbg_{t+1}; \bbxi_t, \bbtheta_t ) $, the last term, by $\alpha$, and reorder terms to write
\begin{align}\label{eq:norm_stoch_grad_expand}
	& \Big\| \frac{\left( f_{t}  - \alpha \hat{\nabla}_f R(f_t,\bbg_{t+1}; \bbxi_t, \bbtheta_t ) \right)}{\alpha}-
	\frac{\ccalP_{ \ccalH_{\bbD_{t+1}}} \Big[ f_t - \alpha \hat{\nabla}_f R(f_t,\bbg_{t+1}; \bbxi_t, \bbtheta_t ) \Big]\Big)}{\alpha }
	\Big\|_{\ccalH}^2 \nonumber \\
	&=  \frac{1}{\alpha^2}\Big\|\left( f_{t}  - \alpha \hat{\nabla}_f R(f_t,\bbg_{t+1}; \bbxi_t, \bbtheta_t ) \right)
	- \ccalP_{ \ccalH_{\bbD_{t+1}}} \Big[\left( f_{t}  - \alpha \hat{\nabla}_f R(f_t,\bbg_{t+1}; \bbxi_t, \bbtheta_t ) \right)\Big]\Big) \Big\|_{\ccalH}^2 \nonumber \\
	& = \frac{1}{\alpha^2}\|\tilde{f}_{t+1} - f_{t+1}\|_{\ccalH}^2 \leq  \frac{\epsilon^2}{\alpha^2}\; ,
\end{align}
where we have taken the common scalar $\frac{1}{\alpha}$ out of the norm. The second equality in \eqref{eq:norm_stoch_grad_expand} holds due to the definition of $\tilde{f}_{t+1}$ and $f_{t+1}$ in \eqref{eq:sgd_tilde} and \eqref{eq:sqg_descent} and substituting here. The last inequality in \eqref{eq:norm_stoch_grad_expand} follows from the stopping criteria used for the KOMP algorithm given by $\lVert \tilde{f}_{t+1} - f_{t+1} \rVert_{\ccalH} \leq \epsilon$. $\hfill\blacksquare	$ \\
%

This result establishes that norm of the difference between the projected and unprojected stochastic quasi-gradient is the ratio of compression budget $\epsilon$ to the step-size $\alpha$.

Next we present an intermediate lemma which is vital to establishing an approximate stochastic descent relationship, and hence convergence.

\begin{lemma}\label{lemma:iterate_relations} For the given algorithm history $\ccalF_t\supset\sigma(\xi_u,\theta_u\}_{u\leq t}$ at time $t$, under the Assumptions \ref{kernel_bound} - \ref{lips_inner}, consider the sequence of iterates $f_t$ generated by Algorithm \ref{alg:colk}. Then:
	\begin{enumerate}
		\item \label{lemma:value_function_difference}The conditional expectation of the Hilbert-norm difference of the next function estimate $f_{t+1}$ and current iteration $f_t$ {$\alpha>0$} satisfies
		\begin{equation}\label{eq:value_function_difference}
			\!\!\!\!\!\mathbb{E}\left[ \|f_{t+1}\! - f_{t} \|_{\ccalH}^2 \given \ccalF_t \right] \leq  4\alpha^2U^2(G_{{\ayche}}^2G_\ell^2+\lambda^2 K^2)+ 2\epsilon^2\; .
		\end{equation}
		\item \label{lemma:value_function_suboptimality}The conditional expectation of the Hilbert-norm difference of the next function estimate $f_{t+1}$ and optimal function $f^\star$ satisfies that 
		\begin{align}\label{eq:value_function_suboptimality}
			\mathbb{E}\left[\|f_{t+1} \!-\! f^\star \|_{\ccalH}^2 \given \ccalF_t \right] &\leq \!\left(\!1 \!+ \! L_{\ell}U^2\frac{\alpha^2}{\beta} G_h^2\right) \| f_t\! -\! f^\star \|_{\ccalH}^2 + 2 \eps \|f_t \!-\!f^\star \|_{\ccalH} \nonumber \\
			&\quad- 2 \alpha \left[\! R(f_t)\! -\! R(\!f^\star\!) \right]
			\!+\! \alpha^2 \sigma_f^2 +  L_{\ell}U^2\beta  \mathbb{E}\left[\|\bar{\bbdelta}_t-\bbg_{t+1}\|_{\ccalH}^2 \given \ccalF_t \right].
		\end{align}
	\end{enumerate}
\end{lemma}
{\bf \noindent Proof of Lemma \ref{lemma:iterate_relations}\eqref{lemma:value_function_difference}}: At the current time instant $t$, consider the Hilbert-norm difference between the next function iterate $f_{t+1}$ and current estimate $f_t$ using the definition of $f_{t+1}$ in \eqref{eq:quasi_projected_fsgd_projected}, i.e.,
\begin{align}\label{eq:value_function_difference_1}
	\|f_{t+1} - f_{t} \|_{\ccalH}^2 &= \alpha^2 \|   \tilde{\nabla}_f R(f_t,\bbg_{t+1}; \bbxi_t, \bbtheta_t )  \|_{\ccalH}^2  \\
	&\quad\leq 2\alpha^2 \|   \tilde{\nabla}_f R(f_t,\bbg_{t+1}; \bbxi_t, \bbtheta_t )   \|_{\ccalH}^2 + 2\alpha^2 \|   \hat{\nabla}_f R(f_t,\bbg_{t+1}; \bbxi_t, \bbtheta_t )  - \tilde{\nabla}_f R(f_t,\bbg_{t+1}; \bbxi_t, \bbtheta_t )   \|_{\ccalH}^2\nonumber \; ,
\end{align}
where we add and subtract the functional stochastic quasi-gradient $\hat{\nabla}_f R(f_t,\bbg_{t+1}; \bbxi_t, \bbtheta_t ) $ on the first line of \eqref{eq:value_function_difference_1} and apply the inequality $(a+b)^2 \leq 2 a^2 + 2 b^2$ which holds for any $a,b$. Now, we apply the result stated in Proposition \ref{prop1} to the second term on the right hand side of \eqref{eq:value_function_difference_1}. After performing this, taking the conditional expectation  on the filtration $\ccalF_t$ yields
\begin{align}\label{eq:value_function_difference_2}
	\mathbb{E}[\|f_{t+1} - f_{t} \|_{\ccalH}^2\given \ccalF_t ]= 2\alpha^2 \mathbb{E}[\|   \hat{\nabla}_f R(f_t,\bbg_{t+1}; \bbxi_t, \bbtheta_t )  \|_{\ccalH}^2 \given \ccalF_t] + 2\epsilon^2 \; .
\end{align}
Next, utilize the definition of the stochastic quasi functional gradient provided in \eqref{eq:quasi_sg} and again using $(a+b)^2\leq 2(a^2+b^2)$, we get
\begin{align}\label{eq:middle}
	\mathbb{E}[\|f_{t+1}\! - \!f_{t} \|_{\ccalH}^2\given \ccalF_t ]\! \leq& 4\alpha^2 \mathbb{E}\Big\{\|{{\langle{\boldsymbol{\ayche}'_{\bbxi_t}( f_t(\bbxi_t))},\ell'_{\bbtheta_t}(\bbg_{t+1})\rangle \kappa(\bbxi_t,\cdot)}}\|^2_{\ccalH}  \given \ccalF_t \!\Big\} +\!4\alpha^2\lambda^2 \|f_t\|_{\ccalH}^2 \!+\! 2\epsilon^2 \; .
\end{align}
Applying the Cauchy-Schwartz inequality to the first term on right hand side of \eqref{eq:middle} yields
{\begin{align}\label{eq:value_function_difference_3}
		\nonumber\mathbb{E}[\|f_{t+1}\! - \!f_{t} \|_{\ccalH}^2\given \ccalF_t ]
		&\leq \ 4  \alpha^2 \mathbb{E}\Big\{\| {{\boldsymbol{\ayche}'_{\bbxi_t}( f_t(\bbxi_t))}}\|^2\|{\ell'_{\bbtheta_t}\left(  \bbg_{t+1}\right)}\|^2\|{{\kappa}( \bbxi_t,\cdot  )}\|^2_{\ccalH}  \given \ccalF_t \!\Big\}+\!4\alpha^2\lambda^2 \|f_t\|_{\ccalH}^2 \!+\! 2\epsilon^2\nonumber
		\\
		&\leq  \ 4\alpha^2U^2 \mathbb{E}\Big\{\mathbb{E}\left[\| {{\boldsymbol{\ayche}'_{\bbxi_t}( f_t(\bbxi_t))}}\|^2\given \bbtheta_t\right] \|{\ell'_{\bbtheta_t}\left(  \bbg_{t+1}\right)}\|^2
		\given \ccalF_t \!\Big\} +4\alpha^2\lambda^2 \|f_t\|_{\ccalH}^2 \!+\! 2\epsilon^2\; ,
\end{align}}
The second inequality in \eqref{eq:value_function_difference_3} is obtained from the Law of total expectation and using Assumption 1 [cf. \eqref{assum1}] {which implies that $ \|{\kappa}( \bbxi_t,\cdot  )\|^2_{\ccalH}\leq U^2$}. Next, due to the set projection in \eqref{eq:projection_hat} and Assumption 1,
we can conclude that $f_t$ has a bounded Hilbert norm, i.e., there
exists some $0 <K <\infty $ such that $\norm{f_t}\leq K$ for all $t$. Using the upper bound of $K$ and the upper bounds as per Assumption 2 [cf. \eqref{assum2_1}], we get the final equation
{\begin{align}\label{eq:value_function_difference_4}
		\mathbb{E}[\|f_{t+1}\! - \!f_{t} \|_{\ccalH}^2\given \ccalF_t ]\! &\leq  4\alpha^2U^2(G_\ayche^2G_\ell^2+\lambda^2 K^2)+ 2\epsilon^2\nonumber
		\\
		&=\alpha^2\sigma_f^2+2\epsilon^2 ,
\end{align}}
where we define $\sigma_f^2:=4U^2(G_\ayche^2G_\ell^2+\lambda^2 K^2)$.

\noindent {\bf Proof of Lemma \ref{lemma:iterate_relations}\eqref{lemma:value_function_suboptimality}}: This proof is a generalization of Lemma 3 in Appendix G.2 in the Supplementary Material of \cite{wang2017stochastic} to a function-valued stochastic quasi-gradient step combined with bias induced by the sparse subspace projections $\ccalP_{\ccalH_{\bbD_{t+1}}}[\cdot]$ in \eqref{eq:projection_hat}. Let us consider the square-Hilbert norm distance of $f_{t+1}$ from the optimal $f^\star$, i.e.,
\begin{align}\label{eq:projected_fsgd_suboptimality1}
	\|f_{t+1} - f^\star\|_{\ccalH}^2 &= \| f_t - \alpha \tilde{\nabla}_f R(f_t, \bbg_{t+1}; \bbxi_t, \bbtheta_t ) -f^\star\|_{\ccalH}^2 \nonumber 
	\\
	& = \| f_t - f^\star\|_{\ccalH}^2 - 2 \alpha \langle  \tilde{\nabla}_f R(f_t, \bbg_{t+1}; \bbxi_t, \bbtheta_t ), f_t -f^\star\rangle_{\ccalH} + \alpha^2 \| \tilde{\nabla}_f R(f_t, \bbg_{t+1}; \bbxi_t, \bbtheta_t )\|_{\ccalH}^2 \; ,
\end{align}
where we utilized the reformulation of the function update defined in \eqref{eq:quasi_projected_fsgd_projected} for the first equality, and expand the square in the second. Now, adding and subtracting $ \hat{\nabla}_f R(f_t, \bbg_{t+1}; \bbxi_t, \bbtheta_t )$ {(which is the (un-projected) functional stochastic quasi-gradient \eqref{eq:quasi_sg}) to first term in the inner product on right hand side of \eqref{eq:projected_fsgd_suboptimality1}} yields
\begin{align}\label{eq:projected_fsgd_suboptimality2}
	\|f_{t+1} - f^\star\|_{\ccalH}^2 & = \| f_t - f^\star\|_{\ccalH}^2 
	- 2 \alpha \langle   \hat{\nabla}_f R(f_t, \bbg_{t+1}; \bbxi_t, \bbtheta_t ), f_t -f^\star\rangle_{\ccalH} \nonumber \\
	&\qquad + 2 \alpha \langle   \hat{\nabla}_f R(f_t, \bbg_{t+1}; \bbxi_t, \bbtheta_t ) -  \tilde{\nabla}_f R(f_t, \bbg_{t+1}; \bbxi_t, \bbtheta_t ), f_t \!-\!f^\star\rangle_{\ccalH}+ \alpha^2 \|  \tilde{\nabla}_f R(f_t, \bbg_{t+1}; \bbxi_t, \bbtheta_t )\|_{\ccalH}^2 \; .
\end{align}
Applying the Cauchy-Schwartz inequality to the third inner product term on the right-hand side of \eqref{eq:projected_fsgd_suboptimality2} and then utilizing the upper bound developed in Proposition \ref{prop1}, we get 
\begin{align}\label{eq:projected_fsgd_suboptimality3}
	\|f_{t+1} - f^\star\|_{\ccalH}^2 & = \| f_t - f^\star\|_{\ccalH}^2
	- 2 \alpha \langle   \hat{\nabla}_f R(f_t, \bbg_{t+1}; \bbxi_t, \bbtheta_t ), f_t -f^\star\rangle_{\ccalH}  + 2 \eps \|f_t -f^\star\|_{\ccalH} 
	+ \alpha^2 \| \tilde{\nabla}_f R(f_t, \bbg_{t+1}; \bbxi_t, \bbtheta_t )\|_{\ccalH}^2 \; .
\end{align}
Now, as defined in Lemma \ref{lemma:iterate_relations}, let {${\boldsymbol{\delta}_t:=\boldsymbol{\ayche}_{\bbxi_t}( f(\bbxi_t))}$} with $\bar{\boldsymbol{\delta}}_t= \mathbb{E}\left[\boldsymbol{\delta}_t \given \bbtheta_t\right]$. Add and subtract $\hat{\nabla}_f R(f_t, \bar{\boldsymbol{\delta}}_t; \bbxi_t, \bbtheta_t )$, which is nothing but the stochastic quasi-gradient evaluated at $(f_t, \bar{\boldsymbol{\delta}}_t)$ rather than $(f_t, \bbg_{t+1})$, inside the inner-product term on the right-hand side of \eqref{eq:projected_fsgd_suboptimality3}, to obtain
\begin{align}\label{eq:projected_fsgd_suboptimality4}
	\|f_{t+1} \!-\! f^\star\|_{\ccalH}^2 & \!=\! \| f_t \!-\! f^\star\|_{\ccalH}^2
	- 2 \alpha \langle  \hat{\nabla}_f R(f_t, \bar{\boldsymbol{\delta}}_t; \bbxi_t, \bbtheta_t ), f_t -f^\star\rangle_{\ccalH}  + 2 \eps \|f_t -f^\star\|_{\ccalH} \nonumber\\
	&  
	\quad {+2 \alpha \langle  \langle{\boldsymbol{\ayche}'_{\bbxi_t}\!( f(\bbxi_t)},\left( {\!\ell'\!\left(\bar{\boldsymbol{\delta}}_t\right)}\!-\! {\ell'\!\left(\bbg_{t+1}\right)}\right)\rangle{{\kappa}( \bbxi_t,\cdot  )}, f_t\! -\!f^\star\rangle_{\ccalH}}+ \alpha^2 \| \tilde{\nabla}_f R(f_t, \bbg_{t+1}; \bbxi_t, \bbtheta_t )\|_{\ccalH}^2 \; ,
\end{align}
where we substitute in the definitions of $\hat{\nabla}_f R(f_t, \bar{\boldsymbol{\delta}}_{t}; \bbxi_t, \bbtheta_t )$ and $\hat{\nabla}_f R(f_t, \bbg_{t+1}; \bbxi_t, \bbtheta_t )$ [cf. \eqref{eq:quasi_sg}, \eqref{eq:quasi_sg}, respectively] in \eqref{eq:projected_fsgd_suboptimality4}, and cancel out the common regularization term $\lambda f_t$. {From the derivation in \eqref{eq:value_function_difference_1} to \eqref{eq:value_function_difference_4}, we can conclude that
	\begin{align}\label{local_bound}
		\alpha^2 \mathbb{E}\left[\|  \tilde{\nabla}_f R(f_t, \bbg_{t+1}; \bbxi_t, \bbtheta_t )\|_{\ccalH}^2 \given \ccalF_t\right]\leq&\alpha^2\sigma_f^2+2\epsilon^2
	\end{align}
	where we define $\sigma_f^2:=4U^2(G_\ayche^2G_\ell^2+\lambda^2 K^2$. }To proceed further, let us define the directional error term related to stochastic quasi-gradient and the stochastic gradient as 
\begin{equation}\label{eq:quasi_gradient_error}
	v_t={2 \alpha \langle  \langle{\boldsymbol{\ayche}'_{\bbxi_t}\!( f(\bbxi_t)},\left( {\!\ell'\!\left(\bar{\boldsymbol{\delta}}_t\right)}\!-\! {\ell'\!\left(\bbg_{t+1}\right)}\right)\rangle{{\kappa}( \bbxi_t,\cdot  )}, f_t\! -\!f^\star\rangle_{\ccalH}}\; .
\end{equation}
From here, compute the conditional expectation  on the algorithm history $\ccalF_t$:
\begin{align}\label{eq:projected_fsgd_suboptimality5}
	\mathbb{E}\left[\|f_{t+1} \!-\! f^\star\|_{\ccalH}^2 \given \ccalF_t \right] & = \| f_t - f^\star\|_{\ccalH}^2
	- 2 \alpha \langle \mathbb{E}\left[ \hat{\nabla}_f R(f_t, \bar{\boldsymbol{\delta}}_t; \bbxi_t, \bbtheta_t )\given \ccalF_t \right], f_t -f^\star\rangle_{\ccalH}  \nonumber \\
	&\quad  
	+ 2 \eps \|f_t \!-\!f^\star\|_{\ccalH}\!+\!  \mathbb{E}\left[ v_t\given \ccalF_t \right]\!\!+\! \alpha^2 \mathbb{E}\left[\|  \tilde{\nabla}_f R(f_t, \bbg_{t+1}; \bbxi_t, \bbtheta_t )\|_{\ccalH}^2 \given \ccalF_t\right].
\end{align}
Utilizing the fact that the compositional objective $R(f)$ defined in \eqref{eq:main_problem} is convex with respect to $f$ and utilizing the first order convexity condition, we have
\begin{align}\label{eq:convexity}
	\langle \mathbb{E}\left[ \hat{\nabla}_f R(f_t, \bar{\boldsymbol{\delta}}_t; \bbxi_t, \bbtheta_t )\given \ccalF_t \right], f_t -f^\star\rangle_{\ccalH}\geq R(f_t) - R(f^\star) \; .
\end{align}
Using the inequality in \eqref{eq:convexity} and the upper bound in \eqref{local_bound} into \eqref{eq:projected_fsgd_suboptimality5}, we get
\begin{align}\label{eq:projected_fsgd_suboptimality6}
	\mathbb{E}\left[\|f_{t+1} - f^\star\|_{\ccalH}^2 \given \ccalF_t \right] & = \| f_t - f^\star\|_{\ccalH}^2
	- 2 \alpha \left[ R(f_t) - R(f^\star) \right] \nonumber \\
	&\qquad+ 2 \eps \|f_t -f^\star\|_{\ccalH}+ \alpha^2 \sigma_f^2  + 2\epsilon^2+  \mathbb{E}\left[ v_t\given \ccalF_t \right] \; .
\end{align}
It remains to analyze $v_t$, the directional error associated with using stochastic quasi-gradients rather than stochastic gradients.  Proceed by applying the Cauchy-Schwartz inequality to \eqref{eq:quasi_gradient_error}, which allows us to write
{\begin{align}\label{eq:quasi_gradient_error22}
		v_t& \leq  2 \alpha (\|{{\boldsymbol{\ayche}'_{\bbxi_t}( f_t(\bbxi_t))}}\|)(\| {\ell'_{\bbtheta_t}\left(\bar{\boldsymbol{\delta}}_t\right)}\!-\! {\ell'_{\bbtheta_t}\!\!\left(\bbg_{t+1}\right)}\|)\| {{\kappa}( \bbxi_t,\cdot  )}\|_{\ccalH}\| f_t \!\!-\!\!f^\star\|_{\ccalH}
		\\
		&\leq  2 \alpha L_{\ell} U^2(\|{{\boldsymbol{\ayche}'_{\bbxi_t}( f_t(\bbxi_t))}}\|)(\|\bar{\boldsymbol{\delta}}_t-\bbg_{t+1}\|)\| f_t -f^\star\|_{\ccalH}\; ,\nonumber
\end{align}}
{where the second inequality in \eqref{eq:quasi_gradient_error22} uses Assumptions 1 [cf. \eqref{assum1}] and 3 [cf. \eqref{assum3}]. Consider Peter-Paul's inequality  $2a b \leq \rho a^2 + b^2/\rho $ for $\rho, a,b>0$, which we apply to \eqref{eq:quasi_gradient_error22} with $a=   \|\bar{\boldsymbol{\delta}}_t-\bbg_{t+1}\| $, $b=\alpha(\|{{\boldsymbol{\ayche}'_{\bbxi_t}( f_t(\bbxi_t))}})\|f_t-f^\star\|_{\ccalH}$, and $\rho = \beta_t$ so that \eqref{eq:quasi_gradient_error22} becomes}
{\begin{equation}\label{eq:quasi_gradient_error3}
		v_t\leq L_{\ell}U^2\left[\beta  \|\bar{\boldsymbol{\delta}}_t-\bbg_{t+1}\|^2 +L_{\ell}U^2  \frac{\alpha^2}{\beta} \|{{\boldsymbol{\ayche}'_{\bbxi_t}( f_t(\bbxi_t))}}\|^2\| f_t -f^\star\|_{\ccalH}^2\right]\; .
\end{equation}}
The conditional mean of $v_t$ [cf. \eqref{eq:quasi_gradient_error}], using \eqref{eq:quasi_gradient_error3}, is then
{\begin{align}\label{eq:quasi_gradient_error_mean}
		\mathbb{E}\left[v_t \given \ccalF_t \right] 
		&\leq L_{\ell}U^2\beta  \mathbb{E}\left[\|\bar{\boldsymbol{\delta}}_t-\bbg_{t+1}\|^2 \given \ccalF_t \right]+  L_{\ell}U^2\frac{\alpha^2}{\beta} \mathbb{E}\left[\|{{\boldsymbol{\ayche}'_{\bbxi_t}( f_t(\bbxi_t))}}\|^2\given \ccalF_t \right] \| f_t -f^\star\|_{\ccalH}^2.  
		\\
		& \leq L_{\ell}U^2\beta  \mathbb{E}\left[\|\bar{\boldsymbol{\delta}}_t-\bbg_{t+1}\|^2 \given \ccalF_t \right]+  L_{\ell}U^2\frac{\alpha^2}{\beta} G_h^2 \| f_t -f^\star\|_{\ccalH}^2 \nonumber\; ,
	\end{align}
}
where we have used the Assumption 1 [cf. \eqref{assum1}]. Now, substitute \eqref{eq:quasi_gradient_error_mean} into the right-hand side of \eqref{eq:projected_fsgd_suboptimality6} and gather like terms:

\begin{align}\label{eq:projected_fsgd_suboptimality_final}
	\mathbb{E}\left[\|f_{t+1} - f^\star\|_{\ccalH}^2 \given \ccalF_t \right] & \leq \left(1 +  L_{\ell}U^2\frac{\alpha^2}{\beta} G_h^2\right) \| f_t - f^\star\|_{\ccalH}^2 + 2 \eps \|f_t -f^\star\|_{\ccalH}+  {L_{\ell}U^2\beta  \mathbb{E}\left[\|\bar{\boldsymbol{\delta}}_t-\bbg_{t+1}\|^2 \given \ccalF_t \right]} \; .
\end{align}
which is as stated in Lemma \ref{lemma:iterate_relations}\eqref{lemma:value_function_suboptimality}.
\hfill$\blacksquare$ \\ 

\blue{
	
	{\begin{lemma}\label{nonconvex:lemma1}
			Consider the case that $R(f)$ [cf. \eqref{regularized_loss}] is Lipschitz smooth with parameter $L_R$. Then under Assumptions \ref{kernel_bound} - \ref{lips_inner}, and the condition that the step-size ratio satisfies  $\frac{\alpha}{\beta}L_{\ell}^2G_h^2U^2 \leq \frac{1}{4}$, i.e., $\alpha$ is sufficiently small relative to $\beta$, then the sequence of functions $f_t$ generated by the proposed algorithm satisfies 
			\begin{align}
				\mathbb{E}[R(f_{k+1})\mid \mathcal{F}_k] 
				\leq  R(f_k) - \frac{\alpha}{4} \| \nabla_f R(f_k) \|^2_\mathcal{H}
				+\beta\mathbb{E}[\left \| \mathbf{g}_{k+1}\!-\! \bar{\mathbf{\delta}}_k\right \|^2 |\mathcal{F}_k] + \frac{\epsilon^2}{2\alpha} + \frac{L_R}{2}\sigma^2_f \alpha^2. \nonumber
			\end{align}
	\end{lemma}	}
	
	{\bf \noindent Proof:}
	For a  Lipschitz smooth $R$ with parameter $L_R$, Taylor's expansion implies the quadratic upper-bound:
	%
	\begin{align}
		R(f_{t+1}) &\leq R(f_{t}) + \ip{\nabla_f R(f_{t}),f_{t+1}-f_{t}}_\mathcal{H} + \frac{L_R}{2}  \| f_{t}-f_{t-1}  \|^2_\mathcal{H}\nonumber
		\\
		&=R(f_{t}) -\alpha \ip{\nabla_f R(f_{t}),\tilde{\nabla}_f R(f_t, \bbg_{t+1}; \bbxi_t, \bbtheta_t )}_\mathcal{H} + \frac{L_R}{2}\alpha^2  \| \tilde{\nabla}_f R(f_t, \bbg_{t+1}; \bbxi_t, \bbtheta_t )\|^2_\mathcal{H}\label{first_non_convex}
	\end{align}
	where the equality in \eqref{first_non_convex} holds from the definition of the projected gradient iterate \eqref{eq:quasi_projected_fsgd_projected}. Next, we establish the upper bound on the inner product term on the right hand side of \eqref{first_non_convex} as follows
	\begin{align}
		-&\alpha\ip{\nabla_f R(f_t),\tilde{\nabla}_f R(f_t, \bbg_{t+1}; \bbxi_t, \bbtheta_t )}_\mathcal{H} \nonumber
		\\
		&\qquad= -\alpha\ip{\nabla_f R(f_t),\hat{\nabla}_f R(f_t, \bbg_{t+1}; \bbxi_t, \bbtheta_t )+\tilde{\nabla}_f R(f_t, \bbg_{t+1}; \bbxi_t, \bbtheta_t )-\hat{\nabla}_f R(f_t, \bbg_{t+1}; \bbxi_t, \bbtheta_t )}_\mathcal{H}
	\end{align}
	where we add subtract the term $-\alpha\ip{\nabla_f R(f_t),\hat{\nabla}_f R(f_t, \bbg_{t+1}; \bbxi_t, \bbtheta_t )}_\mathcal{H}$. After rearranging the terms, we get 
	\begin{align} 
		-&\alpha\ip{\nabla_f R(f_t),\tilde{\nabla}_f R(f_t, \bbg_{t+1}; \bbxi_t, \bbtheta_t )}_\mathcal{H}
		\\
		&\qquad= -\alpha\ip{\nabla_f R(f_t),\hat{\nabla}_f R(f_t, \bbg_{t+1}; \bbxi_t, \bbtheta_t )}+\alpha\ip{\nabla_f R(f_t),\hat{\nabla}_f R(f_t, \bbg_{t+1}; \bbxi_t, \bbtheta_t )-\Tilde{\nabla}_f R(f_t, \bbg_{t+1}; \bbxi_t, \bbtheta_t )}_\mathcal{H} \nonumber
	\end{align}
	where we add subtract the term $\alpha\ip{\nabla_f R(f_t),\hat{\nabla}_f R(f_t, \bbg_{t+1}; \bbxi_t, \bbtheta_t )}$. From the inequality $ab\leq(a^2+b^2)/2$, we obtain 
	\begin{align}
		-&\alpha\ip{\nabla_f R(f_t),\tilde{\nabla}_f R(f_t, \bbg_{t+1}; \bbxi_t, \bbtheta_t )}_\mathcal{H} 
		\nonumber
		\\
		&\qquad\leq -\alpha\ip{\nabla_f R(f_t),\hat{\nabla}_f R(f_t, \bbg_{t+1}; \bbxi_t, \bbtheta_t )}_\mathcal{H}  + \frac{\alpha}{2}\left[ \|\nabla_f R(f_t) \|^2_\mathcal{H}+ \| \hat{\nabla}_f R(f_t, \bbg_{t+1}; \bbxi_t, \bbtheta_t )-\Tilde{\nabla}_f R(f_t, \bbg_{t+1}; \bbxi_t, \bbtheta_t ) \|_\mathcal{H}^2\right] \nonumber
		\\
		&\qquad\leq -\alpha\ip{\nabla_f R(f_t),\hat{\nabla}_f R(f_t, \bbg_{t+1}; \bbxi_t, \bbtheta_t )}_\mathcal{H}  + \frac{\alpha}{2} \|\nabla_f R(f_t) \|^2_\mathcal{H} + \frac{\epsilon^2}{2\alpha}\label{second_non_convex}
	\end{align}
	where the inequality in \eqref{second_non_convex} holds from Proposition \ref{prop1}. Next, consider the term first term on the right hand side of \eqref{second_non_convex} and add subtract $\nabla_f R(f_{t})$ as follows
	\begin{align}\nonumber
		-\alpha \ip{\nabla_f R(f_{t}),\hat{\nabla}_f R(f_t, \bbg_{t+1}; \bbxi_t, \bbtheta_t ) }_\mathcal{H} &= -\alpha \ip{\nabla_f R(f_{t}),\nabla_f R(f_{t})+\hat{\nabla}_f R(f_t, \bbg_{t+1}; \bbxi_t, \bbtheta_t ) -\nabla_f R(f_{t}) }_\mathcal{H}\\
		&= -\alpha  \| \nabla_f R(f_t)\|^2_\mathcal{H} + \alpha\ip{\nabla_f R(f_{t}),\nabla_f R(f_{t}) - \hat{\nabla}_f R(f_t, \bbg_{t+1}; \bbxi_t, \bbtheta_t )}_\mathcal{H}.\label{here_nonconvex}
	\end{align}
	Taking the conditional expectation on the both sides and consider the term $\alpha\mathbb{E} [\ip{\nabla_f R(f_{t}),\nabla_f R(f_{t}) - \hat{\nabla}_f R(f_t, \bbg_{t+1}; \bbxi_t, \bbtheta_t )}_\mathcal{H}\mid \mathcal{F}_t ]$
	
	\begin{align}
		\alpha\mathbb{E} [&\ip{\nabla_f R(f_{t}),\nabla_f R(f_{t}) - \hat{\nabla}_f R(f_t, \bbg_{t+1}; \bbxi_t, \bbtheta_t )}_\mathcal{H}\mid \mathcal{F}_t ]\nonumber
		\\
		&= \alpha\mathbb{E}[ \ip{\nabla_f R(f_{t}),\big(\ip{\mathbf{h}'_{\xi_t}(f(\xi_t)),\ell'_{\theta_t}(\bar{\boldsymbol{\delta}}_t)} - \ip{\mathbf{h}'_{\xi_t}(f(\xi_t)),\ell'_{\theta_t}(\mathbf{g}_{t+1})}\big)\kappa(\xi_t,\cdot)}_\mathcal{H}\mid \mathcal{F}_t ]
		\\
		&\nonumber= \alpha\mathbb{E}[ \ip{\nabla_f R(f_{t}),\big(\ip{\mathbf{h}'_{\xi_t}(f(\xi_t)),\ell'_{\theta_t}(\bar{\boldsymbol{\delta}}_t) -\ell'_{\theta_t}(\mathbf{g}_{t+1})}\big)\kappa(\xi_t,\cdot)}_\mathcal{H}\mid \mathcal{F}_t ]
	\end{align}
	where the equality holds from the definition of the gradient in \eqref{eq:quasi_sg}. Using the Cauchy Schwartz inequality and the upper bounds of Assumption \ref{kernel_bound}, we get
	%
	%
	\begin{align}
		\alpha\mathbb{E} [&\ip{\nabla_f R(f_{t}),\nabla_f R(f_{t}) - \hat{\nabla}_f R(f_t, \bbg_{t+1}; \bbxi_t, \bbtheta_t )}_\mathcal{H}\mid \mathcal{F}_t ]\nonumber
		\\
		&\qquad\leq \alpha U  \| \nabla_f R(f_t)  \|_\mathcal{H} \mathbb{E}[ \|\mathbf{h}'_{\xi_t}(f(\xi_t)) \| \|\ell'_{\theta_t}(\bar{\boldsymbol{\delta}}_t)-\ell'_{\theta_t}(\mathbf{g}_{t+1}) \|\mid \mathcal{F}_t ].
	\end{align}
	Next, applying the Lipschitz smoothness of the outer function as stated in Assumption \ref{lips_outer}, we get
	\begin{align}
		\alpha\mathbb{E} [&\ip{\nabla_f R(f_{t}),\nabla_f R(f_{t}) - \hat{\nabla}_f R(f_t, \bbg_{t+1}; \bbxi_t, \bbtheta_t )}_\mathcal{H}\mid \mathcal{F}_t ]\nonumber
		&\leq \alpha U L_\ell  \| \nabla_f R(f_t) \|_\mathcal{H}\mathbb{E}[ \|\mathbf{h}'_{\xi_t}(f(\xi_t))\| \|\bar{\boldsymbol{\delta}}_t-\mathbf{g}_{t+1} \|\mid \mathcal{F}_t ].
	\end{align}
	Using the inequality $ab \leq \rho a^2 + \frac{1}{\rho}b^2$ for $0 < \rho <1$ and setting $\rho = \beta$, $a =  \|\bar{\boldsymbol{\delta}}_t-\mathbf{g}_{t+1} \|$, and $b = \alpha U L_\ell  \| \nabla_f R(f_t)  \|_\mathcal{H} \|\mathbf{h}'_{\xi_t}(f(\xi_t)) \|$, we obtain
	\begin{align}
		\alpha\mathbb{E} [&\ip{\nabla_f R(f_{t}),\nabla_f R(f_{t}) - \hat{\nabla}_f R(f_t, \bbg_{t+1}; \bbxi_t, \bbtheta_t )}_\mathcal{H}\mid \mathcal{F}_t ]\nonumber
		\\
		&\qquad\leq\beta \mathbb{E}[\left \|\bar{\boldsymbol{\delta}}_t-\mathbf{g}_{t+1}\right \|^2\mid\mathcal{F}_t] + \frac{\alpha^2}{\beta}L_\ell^2 U^2  \| \nabla_f R(f_t) \|^2_\mathcal{H} \mathbb{E}[ \|\mathbf{h}'_{\xi_t}(f(\xi_t)) \|^2\mid\mathcal{F}_t]\nonumber\\
		&\qquad\leq \beta \mathbb{E}[ \|\bar{\boldsymbol{\delta}}_t-\mathbf{g}_{t+1} \|^2\mid\mathcal{F}_t] + \frac{\alpha^2}{\beta}L_\ell^2 G^2_h U^2  \| \nabla_f R(f_t) \|^2_\mathcal{H}\label{third:non_convex}
	\end{align}
	where the second inequality in \eqref{third:non_convex} holds from  Assumption \ref{second_moments}. Substituting the result of \eqref{third:non_convex} into \eqref{here_nonconvex} after taking conditional expectation, we can write 
	\begin{align}
		-\alpha \Ex{\ip{\nabla_f R(f_{t}),\hat{\nabla}_f R(f_t, \bbg_{t+1}; \bbxi_t, \bbtheta_t ) }_\mathcal{H}~|~\mathcal{F}_t} &= -\alpha \ip{\nabla_f R(f_{t}),\nabla_f R(f_{t})+\hat{\nabla}_f R(f_t, \bbg_{t+1}; \bbxi_t, \bbtheta_t ) -\nabla_f R(f_{t}) }_\mathcal{H}\nonumber
		\\
		&= -\alpha  \| \nabla_f R(f_t)\|^2_\mathcal{H} + \beta \mathbb{E}[ \|\bar{\boldsymbol{\delta}}_t-\mathbf{g}_{t+1} \|^2\mid\mathcal{F}_t] + \frac{\alpha^2}{\beta}L_\ell^2 G^2_h U^2  \| \nabla_f R(f_t) \|^2_\mathcal{H}.\label{here_nonconvex0}
	\end{align}
	Similarly, take the conditional expectation of \eqref{second_non_convex} and utilizing the result of \eqref{here_nonconvex0}, we get 
	\begin{align}
		-\alpha \mathbb{E}[\ip{\nabla_f R(f_{t}),\tilde{\nabla}_f R(f_t, \bbg_{t+1}; \bbxi_t, \bbtheta_t )}_\mathcal{H}\mid \mathcal{F}_t] &\leq - \frac{\alpha}{2}\| \nabla_f R(f_t) \|^2_\mathcal{H} + \frac{\alpha^2}{\beta}L_\ell^2 G^2_h U^2  \| \nabla_f R(f_t) \|^2_\mathcal{H} + \beta \mathbb{E}[ \|\bar{\boldsymbol{\delta}}_t-\mathbf{g}_{t+1} \|^2\mid\mathcal{F}_t] + \frac{\epsilon^2}{2\alpha}\nonumber
		\\
		&= - \frac{\alpha}{2}\left(1-\frac{2\alpha}{\beta}L_\ell^2 G^2_h U^2\right) \| \nabla_f R(f_t) \|^2_\mathcal{H} + \beta \mathbb{E}[\left \|\bar{\boldsymbol{\delta}}_t-\mathbf{g}_{t+1}\right \|^2\mid\mathcal{F}_t] + \frac{\epsilon^2}{2\alpha}.\label{here2}
	\end{align}
	Finally, substituting \eqref{here2} into \eqref{first_non_convex} and utilizing Assumption \ref{second_moments}, we obtain the statement of Lemma \ref{nonconvex:lemma1} as follows
	\begin{align}
		\mathbb{E}[R(f_{t+1})\mid \mathcal{F}_t] \leq & R(f_t) - \frac{\alpha}{2}\left(1-\frac{2\alpha}{\beta}L_\ell^2 G^2_h U^2\right) \| \nabla_f R(f_t)  \|^2_\mathcal{H} + \beta \mathbb{E}[\left \|\bar{\boldsymbol{\delta}}_t-\mathbf{g}_{t+1}\right \|^2\mid\mathcal{F}_t] + \frac{\epsilon^2}{2\alpha} + \frac{L_R}{2}\sigma^2_f \alpha^2
		\\
		\leq& R(f_t) - \frac{\alpha}{4} \| \nabla_f R(f_t) \|^2_\mathcal{H} +\beta \mathbb{E}[\left \| \mathbf{g}_{t+1}- \bar{\mathbf{\boldsymbol{\delta}}}_t\right \|^2 \mid \mathcal{F}_t] + \frac{\epsilon^2}{2\alpha} + \frac{L_R}{2}\sigma^2_f \alpha^2
	\end{align}
	which is as stated in Lemma \ref{nonconvex:lemma1}. \hfill$\blacksquare$ \\ 

}

\section{Proof of Theorem \ref{model_order}}\label{model_order_control_supp}
Before discussing the finiteness of the model order, we discuss a lemma which helps us to relate the stopping criterion of specification in KOMP to the Hilbert subspace.
%
\begin{lemma}\label{lemma_subspace_dist}
	{	Let us define the distance of an arbitrary random feature vector  $(\bbxi)$ calculated as $\phi(\bbxi) = \kappa(\bbxi, \cdot)$ to, $\ccalH_{\bbD}=\text{span}\{\kappa(\bbd_n, \cdot) \}_{n=1}^M$, the subspace of the Hilbert space spanned by a dictionary $\bbD$ of size $M$,  as
		\begin{equation}\label{eq:hilbert_subspace_dist}
			\text{dist}( \kappa(\bbxi, \cdot) , \ccalH_{\bbD}) 
			= \min_{f\in\ccalH_{\bbD}} \| \kappa(\bbxi, \cdot) - \bbv^T \bbkappa_{\bbD}(\cdot) \|_{\ccalH} \; .
	\end{equation}}
	This set distance simplifies to following least-squares projection when $\bbD \in \reals^{p\times M}$ is fixed
	{	\begin{equation}\label{eq:hilbert_subspace_dist_ls}
			\text{dist}( \kappa(\bbxi, \cdot) , \ccalH_{\bbD}) 
			= \Big\|  \kappa(\bbxi,\cdot) 
			- [\bbK_{\bbD, \bbD}^{-1} \bbkappa_{\bbD}(\bbxi)]^T
			\bbkappa_{\bbD}(\cdot) \Big\|_{\ccalH} \; .
	\end{equation}}
\end{lemma}
\begin{IEEEproof}
	We can write the distance to a Hilbert space $\ccalH_{\bbD}$ as follows
	\begin{align}\label{eq:subspace_dist}
		\text{dist}( \kappa(\bbxi, \cdot) , \ccalH_{\bbD}) 
		=& \min_{f\in\ccalH_{\bbD}} \| \kappa(\bbxi, \cdot) - \bbv^T \bbkappa_{\bbD}(\cdot) \|_{\ccalH} \nonumber 
		\\
		=& \min_{\bbv\in \reals^{M}} \| \kappa(\bbxi, \cdot) - \bbv^T \bbkappa_{\bbD}(\cdot) \|_{\ccalH} \; ,
	\end{align}
	which is obtained since the dictionary $D$ is fixed and the only free parameter left to optimize is $\bbv$. Now similarly to \eqref{eq:proximal_hilbert_dictionary_polker} - \eqref{eq:hatparam_update}, we can obtain an optimal weight vector $\tbv^*=\bbK_{\bbD_t, \bbD_t}^{-1}\bbkappa_{\bbD_t}(\bbtheta,\bbxi)$ and then substitute back into \eqref{eq:subspace_dist}, we get
	{	\begin{align}\label{eq:subspace_dist2}
			\text{dist}(\kappa(\bbxi, \cdot) , \ccalH_{\bbD}) =  \Big\|  \kappa(\bbxi,\cdot) 
			- [\bbK_{\bbD_t, \bbD_t}^{-1} \bbkappa_{\bbD_t}(\bbxi)]^T
			\bbkappa_{\bbD_t}(\cdot) \Big\|_{\ccalH} \; .
	\end{align}}
\end{IEEEproof}
\subsection{Proof of Theorem \ref{model_order}}\label{model_order_proof}
%
{\bf \noindent Proof:}
Consider two arbitrary time instants $t$ and $t+1$, at which $f_{t}$ and $f_{t+1}$ are the iterates generated by Algorithm \ref{alg:colk}
with corresponding model order denoted by $M_t$ and $M_{t+1}$, respectively.   We consider a constant step size algorithm with fixed approximation budget $\epsilon=C\alpha^2$ for some constant $K>0$. For the sake of analysis, suppose that the model order at $t+1$ is $M_{t+1}\leq M_t$ which means that model order does not grow as we go to iterate $t+1$ from $t$. Note that the model order remains the same from $t$ to $t+1$ if the error introduced by the removal of recently appended data points ($\bbxi_t$) to dictionary $\tbD_{t+1} = [\bbD_t ; \bbxi_t ]$ [cf. \eqref{eq:param_tilde}] is less than the approximation budget $\epsilon$. In other words, the model order does not grow if the stopping criteria of KOMP (Algorithm \ref{alg:komp}), stated as $\min_{j=1,\dots,{M_t + 1}} \gamma_j > \eps$, \emph{is not} satisfied. This leads us to the conclusion that the model order remains the same from $t$ to $t+1$ if 
\begin{equation}\label{eq:komp_no_terminate}
	\min_{j=1,\dots,{M_t + 1}} \gamma_j \leq \eps \; .
\end{equation}
Further observe that the left hand side of \eqref{eq:komp_no_terminate} is a lower bound for the approximation error $\gamma_{M_t + 1}$ at $t+1$  because of the minimization over all $j=1\cdots (M_t + 1)$. This states that the error $\gamma_{M_t + 1}$ introduced by removing the recently appended pair $(\bbxi_t)$ is such that $\min_{j=1,\dots,{M_t + 1}} \gamma_j \leq\gamma_{M_t + 1} $. Therefore, the model order does not grow if $\gamma_{M_t + 1} \leq \eps$ holds because then \eqref{eq:komp_no_terminate} is satisfied. Let use analyze the error term $\gamma_{M_t + 1} $ as follows.	
The definition of $\gamma_{M_t + 1}$ with the substitution of $\tilde{f}_{t+1}$ in \eqref{eq:sgd_tilde} with the notation $$V_t:= \langle{\boldsymbol{\ayche}'_{\bbxi_t}( f_t(\bbxi_t))}, \ell'_{\bbtheta_t}\left(\bbg_{t+1}\right)\rangle, $$ allows us to write
\begin{align}\label{eq:min_gamma_expand}
	\gamma_{M_t + 1}
	&\!=\!\!\!\min_{\bbu\in\reals^{{M_t}}}\!\! \Big\|{(1\!-\!\alpha\lambda)f_t-\alpha V_t {\kappa}(\bbxi_t,\!\cdot\!) }-\!\!\!\!\!\!\!\!\sum_{k \in \ccalI \setminus \{M_t + 1\}} \hspace{-0.6cm}u_k {\kappa}(\bbxi_t,\cdot){\kappa}(\bbd_k, \cdot) \Big\|_{\ccalH}
	\\
	&=\min_{\bbu\in\reals^{{M_t}}} \Big\|(1\! -\! \alpha \lambda) \!\!\!\!\! \!\!\!\sum_{k \in \ccalI \setminus \{M_t + 1\}} \!\!\!  \!\!\!\!\!w_k {\kappa}(\bbxi_t,\cdot){\kappa}(\bbd_k, \cdot)\nonumber 
	\\
	&\hspace{1.5cm}- \alpha V_t{\kappa}(\bbxi_t,\cdot)    -\!\!\!\!\!\! \sum_{k \in \ccalI \setminus \{M_t + 1\}} \!\!\! \!\!\! u_k {\kappa}(\bbxi_t,\cdot){\kappa}(\bbd_k, \cdot)\Big\|_{\ccalH} \; 
	, \nonumber
\end{align}
where we denote the {$k^\text{th}$ data column} of $\bbD_t$ as $\bbd_k$. 
%
To obtain the value $\tbu^*$ which minimizes the above expression, we calculate the gradient of above expression with respect to $\bbu$ and set it to zero. Following the similar logic to that of \eqref{eq:proximal_hilbert_dictionary_polker} - \eqref{eq:hatparam_update}, we get the following result
\begin{align}\label{eq:minimal_weights}
	\tbu^* = (1-\alpha \lambda) \bbw - \alpha V_t \bbK_{\bbD_t, \bbD_t}^{-1} \bbkappa_{\bbD_t}(\bbxi_t) \; .
\end{align}
Next, utilizing the optimal value $\tbu^*$ from \eqref{eq:minimal_weights} into the expression in \eqref{eq:min_gamma_expand} along with the short-hand notation ${f}_{t}(\cdot)=\bbw^T \bbkappa_{\bbD_t}(\cdot)$ and $\sum_k u_k \kappa(\bbd_k, \cdot)= \bbu^T \bbkappa_{\bbD_t}(\cdot)$, we get
\begin{align}\label{eq:min_gamma_optimal_weights}
	&\Big\| (1-\alpha \lambda) \bbw^T\bbkappa_{\bbD_t}(\cdot) - \alpha V_t{\kappa}(\bbxi_t,\cdot)   -  \bbu^T \bbkappa_{\bbD_t}(\cdot) \Big\|_{\ccalH} \\
	&\quad =\! \Big\|(1\!-\!\alpha\lambda\!)\bbw^T\! \bbkappa_{\bbD_t}(\cdot) \!
	-\! \alpha V_t{\kappa}(\bbxi_t,\cdot)  \nonumber 
	\\
	&\hspace{2.2cm}- [(1-\alpha \lambda) \bbw - \alpha V_t \bbK_{\bbD_t, \bbD_t}^{-1} \bbkappa_{\bbD_t}(\bbxi_t)]^T 
	\!\!\bbkappa_{\bbD_t}(\cdot)  \Big\|_{\ccalH} \;  .\nonumber
\end{align}
Further simplifying the above expression by cancelling the similar terms  $(1-\alpha\lambda)\bbw^T\! \bbkappa_{\bbD_t}(\cdot)$ and taking taking the common term $\alpha |V_t|$ outside the norm as
\begin{align}\label{eq:min_gamma_optimal_weights2}
	\Big\| -\alpha V_t{\kappa}(\bbxi_t,\cdot) 
	&  +\alpha V_t [\bbK_{\bbD_t, \bbD_t}^{-1} \bbkappa_{\bbD_t}(\bbxi_t)]^T\bbkappa_{\bbD_t}(\cdot) \Big\|_{\ccalH}  \\
	&  = \alpha|V_t | \Big\|  {\kappa}(\bbxi_t,\cdot) 
	- [\bbK_{\bbD_t, \bbD_t}^{-1} \bbkappa_{\bbD_t}(\bbxi_t)]^T\bbkappa_{\bbD_t}(\cdot) \Big\|_{\ccalH} \nonumber\; .
\end{align}
It is remarked that the norm expression in the right hand side of \eqref{eq:min_gamma_optimal_weights2} describes the distance to the subspace $\ccalH_{\bbD_t}$ as described in \eqref{eq:hilbert_subspace_dist_ls} and defined in Lemma \ref{lemma_subspace_dist} with a scaling factor of $\alpha |V_t| $. The right hand side of \eqref{eq:min_gamma_optimal_weights2} can be written as 
\begin{align}\label{eq:min_gamma_optimal_weights3}
	&\alpha|V_t | \Big\|  {\kappa}(\bbxi_t,\cdot) 
	- [\bbK_{\bbD_t, \bbD_t}^{-1} \bbkappa_{\bbD_t}(\bbxi_t)]^T\bbkappa_{\bbD_t}(\cdot) \Big\|_{\ccalH}  \nonumber 
	\\
	&\hspace{3.5cm}= \alpha |V_t | \text{dist}({\kappa}(\bbxi_t,\cdot),\ccalH_{\bbD_t})\; ,
\end{align}
where the result in \eqref{eq:hilbert_subspace_dist_ls} on the right hand side of  \eqref{eq:min_gamma_optimal_weights3} to replace the Hilbert-norm term. Observe that when the stopping criteria of KOMP is violated, \eqref{eq:komp_no_terminate} holds and thus $\gamma_{M_t + 1} \leq \eps$. Therefore, we have that the right-hand side of \eqref{eq:min_gamma_optimal_weights3} will be upper-bounded by $\epsilon$, and we can write the following inequality
\blue{\begin{align}\label{eq:min_gamma_optimal_weights4}
		\text{dist}({\kappa}(\bbxi_t,\cdot),\ccalH_{\bbD_t}) \leq \frac{\epsilon}{\alpha|V_t |} \;,
\end{align}}
%
%
The error associated with the model order $M_{t+1}$ is denoted by $\gamma_{M_{t+1}}$. Observe that if \eqref{eq:min_gamma_optimal_weights4} holds, then $\gamma_{M_{t+1}} \leq \eps$ holds, but since $\gamma_{M_{t+1}} \geq \min_{j} \gamma_j $, we may conclude that \eqref{eq:komp_no_terminate} is satisfied. Consequently the model order at the subsequent step does not grow $M_{t+1} \leq M_t$ whenever \eqref{eq:min_gamma_optimal_weights4} is valid. 

Now, consider the contrapositive of the preceding expressions. Observe that the model order growth condition ($M_{t+1} = M_t + 1$) implies that 
\blue{	\begin{align}
		\label{eq:min_gamma2}
		\text{dist}({\kappa}(\bbxi_t,\cdot),\ccalH_{\bbD_t}) \geq \frac{\epsilon}{\alpha|V_t |}
\end{align}} 
holds. This condition establishes the fact that every time a new data ($\bbxi$) is appended to kernel dictionary, then the associated product kernel is guaranteed to be at least a distance of $\blue{\frac{\epsilon}{\alpha|V_t |}}$ from every other kernel function in the current model.

Now utilizing the Cauchy Schwartz inequality and Assumption 5 we get $$|V_t|\leq |\ell'_{\bbtheta_t}\left(\bbg_{t+1}\right)|\| {\boldsymbol{\ayche}'_{\bbxi_t}( f_t(\bbxi_t))} \|\leq C_{\ell}L_{\ayche}.$$ This upper bound implies that $1/|V_t| \geq 1/(C_{\ell}L_{\ayche})$, therefore we can lower bound the right hand side in \eqref{eq:min_gamma2} as follows
\begin{align}\label{eq:min_gamma3}
	\blue{	\frac{C{\alpha}}{|V_t|}\geq \frac{\epsilon}{\alpha C_{\ell}L_{\ayche}}\;.}
\end{align}
From \eqref{eq:min_gamma2}, we obtain 
\begin{align}
	\label{eq:min_gamma23}
	\text{dist}({\kappa}(\bbxi_t,\cdot),\ccalH_{\bbD_t}) \geq \frac{\epsilon}{\alpha C_{\ell}L_{\ayche}}.
\end{align}
Therefore, the KOMP stopping criterion is violated for the newest point whenever distinct dictionary points $\bbd_k$ and $\bbd_j$  for $j,k\in\{1,\dots,M_t\}$, satisfy the condition $\|\phi(\bbd_j) - \phi(\bbd_k) \|_2 > \blue{\frac{\epsilon}{\alpha C_{\ell}L_{\ayche}}}$.  Next, we proceed in a similar manner to that of Theorem 3.1 in \cite{1315946}. Note that since the space $\mathcal{U}$ is compact and $\kappa$ is continuous, the  range $\phi(\ccalU) $ (where $\phi(\bbu)=\kappa(\bbu,\cdot)$ for $\bbu \in \ccalU$) of the kernel transformation of feature space $\ccalU$ is compact. This allows us to conclude that the number of balls of radius $\chi$ (here, $\chi = \blue{\frac{\epsilon}{\alpha C_{\ell}L_{\ayche}}}$) required to completely cover the set $\phi(\ccalU)$ is finite (see, e.g., \cite{anthony2009neural}).  

\blue{To prove the main result of \eqref{model_order_0}, we consider the result in \cite[Proposition 2.2]{1315946} which states that for a Lipschitz continuous Mercer kernel $\kappa$ on compact set $\mathcal{X}\subseteq\mathbb{R}^p$, for any training set $\{\bbx_t\}_{t=1}^\infty$ and any $\nu>0$, the number of elements in the dictionary is upper bounded as  
	\begin{align}\label{model_order_1}
		M\leq Y\left(\frac{1}{\nu}\right)^p.
	\end{align}
	where $Y$ is a constant depends upon $\mathcal{X}$ and the kernel function. From the result in \eqref{eq:min_gamma23}, we conclude that $\sqrt{\nu}={\frac{\epsilon}{\alpha C_{\ell}L_{\ayche}}}$, which we may substitute into \eqref{model_order_1} to obtain
	\begin{align}\label{model_order_2}
		M_t\leq YC_\ell^{2p}L_\ayche^{2p}\left(\frac{\alpha}{\epsilon}\right)^{2p}=Y'\left(\frac{\alpha}{\epsilon}\right)^{2p}=\mathcal{O}\left(\frac{\alpha}{\epsilon}\right)^{2p}.
	\end{align}
	as stated in \eqref{model_order_0}. Note that in \eqref{model_order_2}, we have defined a constant $Y'=YC_\ell^{2p}L_\ayche^{2p}$.  We remark that there is a trivial lower bound of $M_t\geq 1$.  \hfill $\qed$}

\end{document}